\documentclass[11pt]{article}
\usepackage[T1]{fontenc}
\usepackage[utf8]{inputenc}
\usepackage{lmodern}
\usepackage[margin=1in]{geometry}
\usepackage{amsfonts,amsmath,amssymb,amsthm,mathtools,mathrsfs,times,color}
\usepackage{microtype}
\usepackage{enumitem,booktabs,array}
\usepackage{tikz}
\usetikzlibrary{patterns}
\usepackage[colorlinks,linkcolor=blue,anchorcolor=blue,citecolor=blue]{hyperref}
\numberwithin{equation}{section}
\newtheorem{theorem}{Theorem}[section]
\newtheorem{lemma}[theorem]{Lemma}
\newtheorem{proposition}[theorem]{Proposition}

\newtheorem{claim}[theorem]{Claim}
\newtheorem{corollary}[theorem]{Corollary}

\theoremstyle{definition}
\newtheorem{definition}[theorem]{Definition}
\newtheorem{remark}[theorem]{Remark}
\setlist[enumerate]{itemsep=4pt,topsep=5pt}

\newcommand{\thmref}[1]{Theorem~\ref{#1}}
\newcommand{\lemref}[1]{Lemma~\ref{#1}}

\newcommand{\defnref}[1]{Definition~\ref{#1}}  %added
\newcommand{\propref}[1]{Proposition~\ref{#1}}  %added

\newcommand{\Ca}{\mathcal C}
\newcommand{\Dcal}{\mathcal D}
\newcommand{\Rm}{\operatorname{Rm}}
\newcommand{\Ric}{\operatorname{Ric}}
\newcommand{\Vol}{\operatorname{Vol}}
\newcommand{\Sing}{\operatorname{Sing}}

\newcommand{\dist}{\operatorname{dist}}
\newcommand{\scr}{\operatorname{scr}}
\newcommand{\ddc}{\sqrt{-1}\partial\bar\partial}
\newcommand{\R}{\mathbb R}
\newcommand{\C}{\mathbb C}
\newcommand{\Z}{\mathbb Z}
\newcommand{\eps}{\varepsilon}

\newcommand{\norm}[1]{\left\lVert #1\right\rVert}

\hypersetup{pdftitle={Convergence of the Calabi flow on K\"ahler surfaces with bounded Sobolev constants},
 pdfauthor={Haozhao Li, Bing Wang, Kai Zheng},
 pdfsubject={}}
\title{Convergence of the Calabi flow on K\"ahler surfaces \\
with bounded Sobolev constants}
\author{Haozhao Li\footnote{Supported by NSFC grant No. 12471058 and a
research fund from University of Science and Technology of China.},
\quad Bing Wang\footnote{Supported by NSFC grant No. 12431003 and a
research fund from University of Science and Technology of China.}
\quad and \quad Kai Zheng\footnote{Supported by NSFC grant No. 12571091.}}
\date{}

\begin{document}
\maketitle
\begin{abstract}
We study global Calabi flows $g(t)$ on compact K\"ahler surfaces with
a uniform Sobolev constant. For any sequence $t_j\to\infty$, set
$g_j(t):=g(t_j+t)$. After passing to a subsequence, these flows converge
to an eternal Calabi flow on a compact extremal K\"ahler orbifold surface.
The convergence holds in the spacetime
$\hat C^\infty$-Cheeger--Gromov topology: it is uniform in the
Gromov--Hausdorff sense on compact time intervals and smooth away
from the orbifold singularities, using comparisons fixed in time.
\end{abstract}
\begingroup
\small
\setcounter{tocdepth}{2}
\tableofcontents
\endgroup

\section{Introduction}
\label{sec:intro}

Calabi introduced extremal K\"ahler metrics as critical points of the
Calabi functional in a fixed K\"ahler class \cite{Calabi82,Calabi85}.
He proposed the Calabi flow as a method for finding these metrics.
The flow decreases the Calabi functional and preserves both the
complex structure and the K\"ahler class. Chen--He established
short-time existence \cite{CHexist}. The next questions concern
long-time existence and convergence as time tends to infinity.

This paper concerns convergence. It is the second paper in our
series on the Calabi flow on K\"ahler surfaces with uniformly
bounded Sobolev constants. The first paper
\cite{arXiv:2610.05217} proves long-time existence under this
bound. Here we study the limiting geometry as $t\to\infty$.
For a global flow $g(t)$ and a sequence $t_j\to\infty$, we
consider the \emph{late-time translates}
$g_j(\tau):=g(t_j+\tau)$, obtained by moving the time origin to
$t_j$. Our aim is to find a topology in which a subsequence of
these flows converges to an extremal K\"ahler orbifold flow
defined for every $\tau\in\mathbb R$.

Let $(M,J)$ be a connected compact complex surface, and fix a real
K\"ahler class $\Omega$. Let $\omega(t)\in\Omega$ be a smooth
Calabi flow
\begin{equation}\label{eq:CF}
 \partial_t\omega=\ddc S,\qquad t\in[0,\infty),
\end{equation}
where $S$ is the complex scalar curvature. Write
$g(u,v):=\omega(u,Jv)$, $dV:=\omega^2/2$, and
\[
 V:=\int_MdV,\qquad \bar S:=V^{-1}\int_M S\,dV.
\]
Both $V$ and $\bar S$ are fixed by the class. We assume that
\begin{equation}\label{eq:Sob}
 \norm{u}_{L^4(g(t))}^{2}
 \le A\left(\int_M|\nabla u|_{g(t)}^2\,dV_{g(t)}
       +b\int_Mu^2\,dV_{g(t)}\right),\qquad b:=V^{-1/2},
\end{equation}
for every smooth $u$, with $A$ independent of time. Define
\begin{equation}\label{eq:energy}
 \Ca(t) :=\int_M(S-\bar S)^2\,dV.
\end{equation}
For a compact K\"ahler orbifold metric $k$, we also write
$\Ca(k):=\int(S_k-\bar S_k)^2\,dV_k$, where
$\bar S_k:=\Vol(k)^{-1}\int S_k\,dV_k$. Thus $\Ca(t)=\Ca(g(t))$.
%The probability-measure convention in \cite{arXiv:2610.05217}
%gives the energy $V^{-1}\Ca(t)$; its Sobolev bound implies
%\eqref{eq:Sob} with $A=C_S\max\{1,b\}$.
The Calabi functional satisfies
\begin{equation}\label{eq:intro-dissipation}
 \frac{d}{dt}\Ca(t)=-2D(t)\le0,\qquad
 D(t):=\int_M|\bar\partial\nabla^{1,0}_{g(t)}S|^2\,dV_{g(t)}.
\end{equation}

We call \eqref{eq:intro-dissipation} the \emph{dissipation identity}
and $D(t)$ the dissipation rate. Since $\Ca(t)$ has a finite
limit, the integral of $D$ over any bounded translated interval
tends to zero as the translation time tends to infinity. In complex dimension two, the Calabi bound also
controls $\int_M|\Rm|^2\,dV$. Together with the uniform Sobolev
inequality, these estimates suggest that the flow approaches
the class of compact extremal K\"ahler orbifolds. The expected
compactness thus resembles that of extremal metrics with bounded
Sobolev constants studied by Chen--Weber \cite{CW}.
This also builds on the orbifold compactness theory for critical
metrics, including constant scalar curvature K\"ahler metrics,
developed by Anderson \cite{Anderson89,Anderson05} and
Tian--Viaclovsky \cite{TVmod,TV}.

There is one essential difference: our flow evolves on a fixed
complex manifold in a fixed K\"ahler class. For this fixed
manifold and class, extremal metrics, if they exist, are unique modulo
holomorphic automorphisms
\cite[Theorem~4.16]{BermanBerndtsson}. Their moduli space modulo
these automorphisms therefore contains at most one point.
The flow slices need not be extremal, so this uniqueness does
not settle our problem. It does suggest that the compactness
theory should combine two advantages: the analytic control
available for extremal metrics, and the topological restrictions
coming from the fixed manifold and class. The analytic control comes from
elliptic estimates and orbifold compactness. The topological
restrictions enter through the integral homology lattice and
K\"ahler periods. To use these restrictions, we also need the
fixed complex structure and estimates for harmonic forms.

\par\addvspace{\baselineskip}
\noindent\textbf{Difficulties and model spaces.}
The main analytic difficulty in establishing this compactness
is regularity improvement.
Calabi flow is fourth order, so the usual parabolic maximum
principle is unavailable for its curvature evolution. This
removes a basic tool for deriving local curvature estimates.
A Sobolev bound does not give pointwise curvature control, and
the flow slices do not satisfy the extremal equation.
The compactness theory for extremal metrics therefore cannot be
applied directly to the flow slices.

We first choose suitable model spaces and study their
compactness. We then compare the evolving metrics with these
models and define a space-time canonical radius that measures
the scales on which this comparison and local regularity hold.
There are two scales in this construction. At the original
scale, the expected models are compact extremal orbifolds.
At shrinking scales, they are complete null scalar-flat
K\"ahler ALE orbifolds. Here null means that the K\"ahler form
has zero integral over every compact smooth closed oriented real
surface in the regular locus. ALE means asymptotically locally
Euclidean: the ends approach flat quotients of $\mathbb C^2$.
The scalar-flatness of the shrinking-scale limits requires a
separate noncompact argument, which we explain next.

The dissipation identity applies at both scales. Under the
Calabi rescaling
\[
 \widehat g_j(\tau):=r_j^{-2}g(t_j+r_j^4\tau),
 \qquad
 \widehat D_j(\tau):=r_j^4D(t_j+r_j^4\tau),
\]
one has, for fixed $a<b$,
\[
 \int_a^b\widehat D_j(\tau)\,d\tau
 =\int_{t_j+r_j^4a}^{t_j+r_j^4b}D(t)\,dt.
\]
The tensor $\bar\partial\nabla^{1,0}S$ is the \emph{extremal defect}.
Thus the spacetime integral of its squared norm is scaling invariant
in real dimension four. This integral tends to zero whenever the left
endpoint $t_j+r_j^4a$ tends to infinity. Once regular convergence is
established, the limiting metrics are therefore extremal.
At the original scale, they may have nonconstant scalar
curvature. At shrinking scales, the completed limits are
noncompact and inherit a homogeneous Sobolev inequality and
finite curvature energy. Decay estimates and the Jacobi equation
force the Killing field $J\nabla S$ to vanish. The scalar curvature
is therefore constant, and its decay at infinity makes this
constant zero (Lemma~\ref{lem:scalarflat}). The finite-end and
decay arguments subsequently give the ALE structure
(Proposition~\ref{prop:ALE}). Scalar-flatness is therefore a further
consequence of the noncompact geometry, not of vanishing
dissipation alone.

Choosing the shrinking-scale model requires more care.
There is initially no local volume-ratio upper bound for the
rescaled flows. Sobolev noncollapse gives the lower bound, but
does not by itself justify extraction of an orbifold limit.
For the proposed stationary models, we prove a uniform upper
bound by an orbifold version of Tian--Viaclovsky's volume-growth
argument (Proposition~\ref{prop:stationaryVolume}). In the
scale-selection argument, we distinguish comparable scales
from separated scales, addressing the issue pointed out by
Chen--Weber \cite{CW}. The resulting upper bound is included
among the canonical-radius conditions. It is
not assumed independently for the original flow.

A second obstacle is the Burns metric on the blow-up of
$\mathbb C^2$ at the origin. It is nonflat and scalar-flat,
but its tangent cone at infinity is Euclidean
\cite[Section~1.4]{JV}. Its metric blow-down is $\mathbb C^2$.
Thus the class of all scalar-flat ALE spaces allows curvature
to concentrate above a point whose tangent cone is Euclidean.
In that class, geometric regularity alone does not identify
the region of smooth convergence.

The fixed complex manifold and K\"ahler class supply the missing
restriction. Every such shrinking-scale limit $(Y,h)$ is \emph{null}:
\[
 \int_\Sigma\omega_h=0
 \quad\text{for every compact smooth closed oriented real surface }
 \Sigma\subset \mathcal R(Y).
\]
The fixed-class period observation behind this condition goes
back to the second named author. It was used in Chen--He's
blow-up analysis of toric Fano surfaces
\cite[Propositions~3.3--3.4]{CHtoric}; see also their acknowledgment.
Section~\ref{sec:fixedclass} proves the null-period statement
for arbitrary real K\"ahler classes. The fixed complex structure
and class control the $L^2$ norms of harmonic representatives.
Pairing with compactly supported Poincar\'e dual forms then bounds the
integral coordinates of transported surface classes and forces
the limiting periods to vanish. This extends the null-period argument
of Chen--Wang for K\"ahler--Ricci flow on Fano surfaces
\cite[Proposition~4.4]{ChenWangFano}.

Nullness excludes the Burns metric because its exceptional
curve has positive K\"ahler area. More generally, the mass
formula of Hein--LeBrun \cite{HL} and the orbifold rigidity
theorem of Ju--Viaclovsky \cite{JV} imply that a model with a
Euclidean end must be Euclidean space (Lemma~\ref{lem:nullAE}).
Each model has one end.
Consequently every nonflat model has tangent cone at infinity
$\mathbb C^2/\Gamma$ with $\Gamma\ne\{1\}$.
This is the gap that excludes hidden concentration over
regular limit points. It allows the regular--singular
decomposition to be determined by tangent cones, with smooth
convergence on the whole regular locus.

\par\addvspace{\baselineskip}
\noindent\textbf{The main theorem.}
At the original scale, we must also compare the compact extremal
slice limits at different times. These comparisons must identify
the limits of points held fixed in the underlying manifold;
independent identifications at each time do not suffice. The resulting convergence combines
Gromov--Hausdorff (GH) convergence uniform on compact time intervals
with smooth convergence on regular domains using comparison
embeddings independent of time. This is the spacetime
$\hat C^\infty$-Cheeger--Gromov (CG) topology developed from
Li--Li--Wang \cite{LLW} in Definition~\ref{def:topology}.
Our main result is the following.

\begin{theorem}[Long-time orbifold compactness]
\label{thm:main}\label{target:main}
Assume that the Calabi flow~\eqref{eq:CF} exists on $[0,\infty)$ and
satisfies~\eqref{eq:Sob}. For every sequence $t_j\to\infty$, define the translated flows
\[
 g_j(\tau):=g(t_j+\tau),
 \qquad \tau\in[-t_j,\infty).
\]
Then, after passage to one subsequence, the following
conclusions hold.
\begin{enumerate}[label=(\Roman*)]
\item
The translated flows $(M,g_j(\tau),J)$ converge in the spacetime
$\hat C^\infty$-Cheeger--Gromov topology of
Definition~\ref{def:topology} on every compact time interval in
$\mathbb R$. The limit is a smooth eternal Calabi flow $h(\tau)$
on a connected compact K\"ahler orbifold $(X,J_\infty)$ with
isolated quotient singularities. Write
\[
 X=\mathcal R\sqcup\mathcal S,\qquad
 \mathcal R:=\mathcal R(X),\qquad \mathcal S:=\operatorname{Sing}X.
\]
The singular set $\mathcal S$ is finite and may be empty.
The regular locus $\mathcal R$ is smooth and geodesically convex
for every $h(\tau)$: every minimizing geodesic with endpoints
in $\mathcal R$ lies entirely in $\mathcal R$.

\item
Every limiting metric $h(\tau)$ is extremal.
The orbifold vector field
\[
 W_0:=\frac12\nabla_{h(0)}S_{h(0)}
\]
is complete and real holomorphic, meaning that
$\mathcal L_{W_0}J_\infty=0$, or equivalently that its local
flow preserves the complex structure. It generates a one-parameter group $F_\tau$ of orbifold
biholomorphisms, with $F_0=\mathrm{id}_X$, such that
\begin{equation}\label{eq:mainfamily}
 h(\tau)=F_\tau^*h(0),
 \qquad
 \bar\partial_{J_\infty}
 \nabla^{1,0}_{h(\tau)}S_{h(\tau)}=0,
 \qquad \tau\in\mathbb R.
\end{equation}
Each $F_\tau$ fixes every point of $\mathcal S$.
With the time-independent comparison maps in \textup{(I)},
the singular set of the limiting spacetime is
$\mathcal S\times\mathbb R$. The maps $F_\tau$ describe the
evolution of the limiting metric and are distinct from these
comparison maps.

\item
The volume measures converge weakly, uniformly on every
compact time interval, in the sense of
Definition~\ref{volume measures converge}.
Moreover,
\[
 \operatorname{Vol}(X,h(\tau))=V,
 \qquad
 \int_X S_{h(\tau)}\,dV_{h(\tau)}=V\bar S,
 \qquad \tau\in\mathbb R.
\]

\item
The Calabi functional converges uniformly on every compact
time interval: for every $T>0$,
\begin{equation}\label{eq:maincalabiconvergence}
 \lim_{j\to\infty}\sup_{|\tau|\le T}
 \bigl|\Ca(g(t_j+\tau))-\Ca(h(\tau))\bigr|=0.
\end{equation}
\end{enumerate}
\end{theorem}

Theorem~\ref{thm:main} supplies an extremal orbifold limit without
assuming that an extremal metric already exists in $\Omega$.
Its singularities describe the possible degeneration at infinite
time. The limiting flow evolves by holomorphic automorphisms,
so all its time slices are isometric. If the limiting orbifold
has no nonzero holomorphic vector field, the flow is stationary
and its scalar curvature is constant
(Corollary~\ref{cor:noaut}). For each sequence $t_j\to\infty$,
the theorem produces a convergent subsequence of the translated
flows. Determining whether all choices give the same extremal
orbifold requires additional geometric or algebraic information.

\par\addvspace{\baselineskip}
\noindent\textbf{Outline of the proof of Theorem~\ref{thm:main}.}
\par\addvspace{\baselineskip}
\noindent\emph{Compactness of the models.}
We first establish compactness of the null scalar-flat ALE
model space. Its definition prescribes a homogeneous Sobolev
bound, a curvature-energy bound, and a bound on the number of
quotient points. Counting the latter is necessary because a
flat quotient has zero curvature energy. The volume upper
bound, one-end theorem, and finite bubble trees give preliminary
orbifold compactness. The rigidity of null models with Euclidean
ends excludes nonflat bubbles and hidden quotient points above
a regular limit point.
The comparison charts therefore cover every approximating
point sufficiently close to a compact regular region.
This gives compactness in the pointed
$\hat C^\infty$-Cheeger--Gromov topology
(Theorem~\ref{thm:stationary}).

\par\addvspace{\baselineskip}
\noindent\emph{Canonical radius and local backward pseudo-locality.}
We use these models to define the space-time canonical radius
$\scr$. At each scale $r$, the tests require model comparison
over $[t-r^4,t]$, a terminal volume bound, and backward
regularity on a fixed fraction of this interval when the
terminal curvature energy is small. If $\scr$ degenerates,
point selection gives an almost-minimal point and a rescaled
sequence with a lower canonical-radius bound on larger
backward regions. The defining estimates supply provisional
volume bounds and short regular cylinders before any
spacetime compactness theorem is available.

The selected original times tend to infinity. The
scaling-invariant time-integrated dissipation therefore tends
to zero, making the local limits extremal. Branchwise
completion initially permits several branches at a point.
The inherited model comparison and the model gaps exclude
branching and concentration at geometrically regular points.
On the completed orbifold slices, the homogeneous Sobolev
inequality and finite curvature energy give decay of the
Killing field $J\nabla S$. Its Jacobi equation along minimizing
rays forces it to vanish. The slices are therefore scalar-flat
and ALE, and the fixed-class period argument makes them null.
Extension from a common ALE end and distance continuity
identify the same points through time. The resulting stationary
model improves the defining estimates at the selected point,
contradicting its normalized canonical radius.

The positive radius bound gives local backward pseudo-locality:
a curvature bound on a terminal ball controls a smaller fixed
ball over a backward interval of length comparable to $r^4$
(Corollary~\ref{cor:backwardpseudolocality}). This is a local
regularity improvement; it requires no curvature bound on the
rest of the slice. In comparison, Li--Wang--Zheng's regularity
improvement \cite{LWZ} assumes control of the whole terminal
slice. Short backward cylinders and vanishing dissipation now
give compact extremal orbifold slice limits, including limits
at moving times.

\par\addvspace{\baselineskip}
\noindent\emph{Long-time two-sided pseudo-locality.}
To extend this control over a prescribed time interval, we
identify the completed slices by compact exterior extension
(Lemma~\ref{lem:compactextension}). The local isometry is
defined on the regular exterior of neighborhoods of all
singular points. Real analyticity permits continuation of this
isometry, and equality of total volumes rules out a nontrivial
covering. This compact argument differs from extension from an
ALE end. The resulting
extremal family is $h(\tau)=F_\tau^*h(0)$. The holomorphic flow
$F_\tau$ preserves the regular locus, and the trajectories of
a compact regular set remain in a compact regular region on
every finite interval (Lemma~\ref{lem:regularorbits}).

This preservation must still be related to the original
manifolds. Distance continuity, distances to fixed regular
points, and overlapping backward charts identify each fixed
sequence of points of the underlying manifold with one point
of the limit throughout the interval. Starting the construction
at a later terminal slice also gives forward control.
Together with the canonical-radius lower bound and a curvature
bound on a larger reference-time ball, a sufficiently small
Calabi-energy drop yields long-time two-sided pseudo-locality,
provided that the flow exists on a larger time interval
(Theorem~\ref{thm:upgrade}). This estimate is proved before
spacetime weak compactness is invoked.

A general forward pseudo-locality theorem derived directly
from local initial geometric control is not currently available
for Calabi flow. Our quantitative theorem applies to flow
segments satisfying the stated geometric bounds and small
energy-drop condition. For the global flow studied here, the
canonical-radius bound and small energy-drop condition hold on
sufficiently late translated intervals. The estimate then applies
wherever the reference-time curvature bound is satisfied.
In this late-time setting, forward pseudo-locality is obtained
together with backward pseudo-locality.

\par\addvspace{\baselineskip}
\noindent\emph{Spacetime convergence.}
Starting from the limit at time zero, we retain its regular
embeddings $\Phi_j:U\to M$ throughout each compact time interval.
Long-time two-sided pseudo-locality gives uniform mixed
derivative bounds in these fixed embeddings. A diagonal
extraction over an exhaustion of the regular locus and
increasing compact time intervals gives smooth convergence on
regular spacetime domains. The identification of fixed sequences
of points, valid also at moving times, shows that the same GH
correspondences have distortion tending to zero uniformly on
each compact time interval. These correspondences contain the
graphs of the regular embeddings.
These two requirements define the spacetime
$\hat C^\infty$-Cheeger--Gromov convergence. Because the comparison
maps are independent of time, the Calabi equation passes to the
limit without an additional Lie derivative.
The singular set is $\mathcal S\times\mathbb R$, and smooth
regular convergence together with volume control gives
convergence of the natural measures. This proves
Theorem~\ref{thm:main}\textup{(I)--(III)}. A fixed finite group
of holomorphic isometries also induces an action on the limit,
as proved in Subsection~\ref{sec:equivariantcompactness}.

\par\addvspace{\baselineskip}
\noindent\emph{Continuity of the Calabi functional.}
Finally, Subsection~\ref{sec:calabicontinuity} proves \textup{(IV)}.
At nearby times of small dissipation, scalar-flatness controls
the shrinking-scale regular cores. A length-independent
estimate controls the scalar energy on the conical regions
between them. Together these exclude scalar-energy loss.
Monotonicity and isometry of the limiting slices then give
Calabi-energy convergence uniformly on compact time intervals.\qed
\begingroup
\widowpenalty=10000
\par\addvspace{\baselineskip}
\noindent\textbf{Relation to earlier work.}
The strategy of the proof originates in Chen--Wang's \emph{Space of
Ricci flows (I)} and \emph{(II)}
\cite{ChenWangI,ChenWangA,ChenWangB}. Part~(I) develops the
interaction between blow-up analysis, backward pseudo-locality,
and distance continuity. Part~(II) organizes the argument
through a compact model space, canonical radius, and long-time
two-sided pseudo-locality. The canonical-radius terminology
comes from Part~(II) \cite[Definition~3.5]{ChenWangA}. Its
space-time formulation follows
\cite[Definition~2.7]{ChenWangRemarks}. See also the survey
\cite[Section~2, p.~385]{ChenWangSurvey}.
\par\endgroup

In the K\"ahler--Ricci flow theory, the shrinking-scale models
are Calabi--Yau spaces with mild singularities: their regular
loci carry Ricci-flat K\"ahler metrics. The original-scale limits in
the Fano setting are K\"ahler--Ricci solitons. Here the models
are null scalar-flat ALE orbifolds and compact extremal
orbifolds. For K\"ahler--Ricci flow, the scalar-curvature evolution
forces regular rescaled limits with zero scalar curvature to
be Ricci-flat. For Calabi flow, vanishing scaling-invariant
spacetime dissipation gives extremality, and noncompact decay
gives scalar-flatness. Nullness supplies the rigidity needed
to exclude Euclidean-end bubbles. The time scale is $r^4$
rather than $r^2$, because Calabi flow is fourth order whereas
K\"ahler--Ricci flow is second order.
Both theories pass from slice to spacetime compactness using
long-time regularity and fixed maps.

The compactness theorem provides a starting point for studying
the final geometry of the Calabi flow. Further information can
be used to exclude orbifold singularities, identify the complex
surface, and prove uniqueness of the limit. The theorem also
supplies geometric input for the companion
papers on polarized surfaces \cite{Polarized} and on the
controlled K\"ahler cone \cite{ControlledCone}.

\par\addvspace{\baselineskip}
\noindent\textbf{Organization of the paper.}
The paper is organized as follows. Section~\ref{sec:preliminaries}
collects the basic estimates, tools for passing to the limit, and
shrinking-scale arguments for a fixed K\"ahler class.
Section~\ref{sec:stationary} establishes compactness of the null
scalar-flat ALE models, including volume growth, one end, and
the regularity gap. Section~\ref{sec:scr} develops the canonical
radius, interior Calabi smoothing, conditional slice compactness,
and local backward pseudo-locality. 
Section~\ref{sec:upgrade} proves long-time two-sided
pseudo-locality. Section~\ref{sec:mainproof} defines the spacetime
topology and proves weak compactness and continuity of the Calabi
functional, completing Theorem~\ref{thm:main}. It ends with the
finite-group equivariant version.

\par\addvspace{\baselineskip}
\noindent\textbf{On the usage of AI.}
The authors used ChatGPT to assist in preparing this paper
and take full responsibility for its contents.

\par\addvspace{\baselineskip}
\noindent\textbf{Acknowledgement.}
The authors are deeply grateful to Professor Xiuxiong Chen
for his encouragement and constant support.

\section{Preliminaries}
\label{sec:preliminaries}
\label{sec:bounds}\label{sec:transport}\label{sec:local}
\subsection*{Notation, terminology and constants}
\addcontentsline{toc}{subsection}{Notation, terminology and constants}
\label{index:notation}
\noindent\emph{Flow and convergence.}
\begin{description}[style=nextline,leftmargin=3.5cm,labelwidth=3.3cm,labelsep=0.3cm,align=left,font=\normalfont,itemsep=0pt,parsep=0pt,topsep=4pt]
\item[$\Omega$] The fixed real K\"ahler class.
\item[$V$] The fixed total volume.
\item[$S,\bar S$] Complex scalar curvature and its average.
\item[$A$] The uniform Sobolev constant in \eqref{eq:Sob}.
\item[$b=V^{-1/2}$] The lower-order coefficient in \eqref{eq:Sob}.
\item[$\Ca$] The Calabi functional; see \eqref{eq:energy}.
\item[$\Dcal_g$] The operator $\bar\partial\nabla_g^{1,0}$; see \eqref{eq:dissipation}.
\item[Extremal defect] The tensor $\Dcal_gS$; it vanishes precisely when the metric is extremal.
\item[Late-time translates] The flows $g_j(\tau):=g(t_j+\tau)$ with $t_j\to\infty$.
\item[Late-time slices] The metrics $g(u_j)$ with $u_j\to\infty$.
\item[$D(t)$] The dissipation integral $\int_M|\Dcal_{g(t)}S|^2\,dV$.
\item[$E_0$] The uniform bound for $\int_M|\Rm|^2\,dV$ (Lemma~\ref{lem:energy}).
\item[$\kappa$] The volume noncollapse constant in Lemma~\ref{lem:noncollapse}.
\item[$d_0$] The diameter bound in Lemma~\ref{lem:noncollapse}.
\item[$\mathcal R,\mathcal S$] The regular and singular loci: $\mathcal R=\mathcal R(X)$, $\mathcal S=\Sing X$.
\item[$\mathcal E$] The provisional exceptional set in Proposition~\ref{prop:slice}.
\item[$W_0,F_\tau$] The extremal vector field and its holomorphic flow; see \eqref{eq:mainfamily}.
\item[$\hat C^\infty$-CG] Compatible GH and smooth convergence on the regular part (Definition~\ref{def:modelhat}). Its spacetime version uses time-independent comparisons and convergence on compact time intervals (Definition~\ref{def:topology}).
\item[$\delta_I$] The distance \eqref{eq:familyGH} between compact metric families modulo a single time-independent isometry, using one correspondence for all times in $I$.
\item[Volume convergence] Weak convergence against continuous test functions, uniformly on compact time intervals through fixed GH maps; see Definition~\ref{volume measures converge}.
\item[Coverage] For $K\Subset U\Subset\mathcal R(X)$ and comparison charts $\Phi_j:U\to X_j$, every sequence $x_j$ with GH limit in $K$ eventually lies in $\Phi_j(U)$. The same requirement applies along subsequences; see Section~\ref{Convergence on the regular part}.
\end{description}

\par\addvspace{\baselineskip}
\noindent\emph{Models.}
\begin{description}[style=nextline,leftmargin=3.5cm,labelwidth=3.3cm,labelsep=0.3cm,align=left,font=\normalfont,itemsep=0pt,parsep=0pt,topsep=4pt]
\item[ALE, AE] Asymptotically locally Euclidean and asymptotically Euclidean; see \eqref{eq:modelALEdefinition}.
\item[Nullness] The condition $\int_\Sigma\omega=0$ for every compact smooth closed oriented real surface $\Sigma$ in the regular locus; see \eqref{eq:modelNull}.
\item[$\mathcal A(A,\Lambda,\alpha)$] The model moduli: pointed null scalar-flat K\"ahler ALE orbifolds, up to pointed biholomorphic isometry, with Sobolev and weighted-energy bounds (Definition~\ref{def:modelmoduli}).
\item[$\mathscr E_\alpha$] Weighted energy: $\int_X|\Rm|^2dV+\alpha\,\#\Sing X$; see \eqref{eq:modelComplexity}.
\item[$\mu_j$] Marked measure: $|\Rm_{h_j}|^2dV_{h_j}+\alpha\sum_{q\in\Sing X_j}\delta_q$, where $\delta_q$ is unit point mass at $q$; see \eqref{eq:modelMarkedMeasure}.
\item[Labels and marks] Indices for selected nonflat bubbles and original singular points. A mark records an original singular point, assigned scale zero.
\item[Cluster] Labels converging to the same attachment point of a parent bubble, with scales negligible relative to the parent; see Lemma~\ref{lem:modelhierarchy}.
\item[$\eps_{\rm reg}$] The stationary regularity threshold in \eqref{eq:modelEpsRegularity}.
\item[$\omega_4$] The volume of the Euclidean unit ball in real dimension four.
\item[$G(A)$] The bound for orbifold group orders in \eqref{eq:modelGroup}.
\end{description}

\par\addvspace{\baselineskip}
\noindent\emph{Canonical radius.}
The model and test constants are fixed in Section~\ref{Choice of constants}.
\begin{description}[style=nextline,leftmargin=3.5cm,labelwidth=3.3cm,labelsep=0.3cm,align=left,font=\normalfont,itemsep=0pt,parsep=0pt,topsep=4pt]
\item[$A_*$] The homogeneous Sobolev bound for the model class, with $A_*>2A$.
\item[$\Lambda$] The model weighted-energy bound, chosen as $2E_0+1$.
\item[$\alpha$] The weight of each singular point in $\mathscr E_\alpha$, chosen as $\eps_R/4$.
\item[$N$] The derivative order $12$ used in the canonical-radius tests.
\item[$\theta$] The backward time fraction $1/16$ in test \textup{(R)}.
\item[$\eps_{\rm stat}$] The small-energy threshold for the stationary estimate \eqref{eq:modelinterior}.
\item[$H_m$] The stationary curvature-derivative constants in \eqref{eq:modelinterior}.
\item[$\eps_R$] The small-energy threshold in test \textup{(R)}.
\item[$C_m$] The constants $4H_m+1$ in test \textup{(R)}.
\item[$V_*$] The uniform upper volume ratio of the model class.
\item[$B_0$] The larger volume bound used in test \textup{(V)}.
\item[$a_h$] The model scale defined by injectivity radius and curvature derivatives.
\item[$\rho_h$] The $1$-Lipschitz regularity scale obtained from $a_h$ in \eqref{eq:modelrho}.
\item[$a_E$] The regularity lower bound for models sufficiently close to Euclidean space.
\item[$\zeta_E$] The GH closeness threshold for Euclidean regularity.
\item[$\zeta_B$] The GH separation from flat cones with several branches.
\item[$G_*$] The upper bound for the quotient-group order of each provisional branch.
\item[$N_B$] The upper bound for the number of provisional branches.
\item[$\eta$] The model regularity threshold selecting the domains for smooth comparison.
\item[$\delta$] The required accuracy of the model comparisons; see \eqref{eq:tolerances}.
\item[$G(\tau)$] The normalized flow in Section~\ref{sec:canonicaltests}.
\item[$W_G(R)$] The union of normalized radius-$R$ balls over $\tau\in[-1,0]$.
\item[$K\subset K^+\subset U$] The nested regular model domains in test \textup{(P)}.
\item[$\mathscr C,\Phi$] A compact relation matching model points with points of the manifold, and a regular comparison embedding. Both are fixed in time; see test \textup{(P)}.
\item[Capture condition] The ball inclusion \eqref{eq:scrcapture}: small ambient balls lie in the comparison chart image.
\item[\textup{(P), (R), (V)}] Model comparison, backward regularity, and terminal volume control.
\item[$\scr(p,t)$] The space-time canonical radius; each scale $s$ tests $[t-s^4,t]$ (Definition~\ref{def:scr}).
\item[$r_0$] The positive lower radius bound on $M\times[t_0,\infty)$; it may depend on the fixed flow and $t_0$ (Theorem~\ref{thm:radiuslower}).
\end{description}

\par\addvspace{\baselineskip}
\begin{samepage}
\noindent\emph{Long-time estimates.}
\begin{description}[style=nextline,leftmargin=3.5cm,labelwidth=3.3cm,labelsep=0.3cm,align=left,font=\normalfont,itemsep=0pt,parsep=0pt,topsep=4pt]
\item[\textup{(B1)--(B3)}] Slice compactness, short backward charts, and distance continuity (Corollary~\ref{cor:geometricinputs}, Lemma~\ref{lem:shortdrop}).
\item[$\eps$] The energy-drop threshold in Theorem~\ref{thm:upgrade}.
\item[$C_m$] The curvature-derivative bounds in Theorem~\ref{thm:upgrade}, distinct from those in \textup{(R)}.
\end{description}
\end{samepage}

\par\addvspace{\baselineskip}
\noindent\emph{Terms used in the proofs.}
\begin{description}[style=nextline,leftmargin=3.5cm,labelwidth=3.3cm,labelsep=0.3cm,align=left,font=\normalfont,itemsep=0pt,parsep=0pt,topsep=4pt]
\item[Anchors] A chosen countable dense set of reference points in an open regular part of the common ALE end.
\item[End-normalized identification] The global isometry extending the prescribed identification on the common end.
\item[Regular tube] A connected open set with compact closure in the regular locus, joining specified regular regions.
\item[Regular-time set] The set $\{\tau:\gamma(\tau)\in\mathcal R(X)\}$ for the limiting path $\gamma$ of a fixed sequence of points.
\item[Spatial margin] A lower bound $cr$ for the ambient distance between the smaller chart image and the complement of the larger image, uniform in the sequence and time.
\item[Past-time margin] An additional backward interval of length at least $cr^4$ on which estimates are available, with $c>0$ uniform in the sequence.
\item[Cone atlas] Cone annular charts with the nested domains, bounded overlap, and uniform $C^2$ estimates specified in Lemma~\ref{lem:atlasDneck}.
\end{description}

\subsection{Basic estimates and rescaling}
\label{sec:basic-estimates}

With the conventions of \eqref{eq:CF}, the variation identity is
\begin{equation}\label{eq:dissipation}
\frac{d}{dt}\Ca(t)=-2D(t),\qquad
D(t)=\int_M|\Dcal_{g(t)}S|^2\,dV,\qquad
\Dcal_g f:=\bar\partial\nabla_g^{1,0}f.
\end{equation}
See \cite[Section~3]{MR2953047}.
In particular, $\Dcal_g f=0$ means that the complex gradient of $f$
is a holomorphic vector field.

The four-dimensional Chern--Weil formulas and K\"ahler curvature
decomposition \cite{Besse,CH} bound the curvature energy in terms of
the Calabi functional; see \cite[Lemma 3.12]{arXiv:2610.05217}.

\begin{lemma}[Uniform curvature energy]\label{lem:energy}
Along the Calabi flow,
\[
\int_M|\Rm|^2\,dV\le E_0,
\]
where $E_0$ depends only on the topology, the fixed K\"ahler class, and $\Ca(0)$.
\end{lemma}

The Sobolev volume lower bound is due to Y. D. Wang \cite{WangYD}
and G. Carron \cite{Carron}.

\begin{lemma}[Noncollapse and diameter]\label{lem:noncollapse}
The Sobolev inequality \eqref{eq:Sob} implies constants $\kappa>0$ and $d_0<\infty$, depending only on $A,V$, such that
\[
\Vol_{g(t)}B(x,r)\ge\kappa r^4,\quad 0<r\le V^{1/4},
\qquad \operatorname{diam}(M,g(t))\le d_0.
\]
Consequently the individual time slices form a GH precompact family of compact metric spaces.
\end{lemma}
\begin{proof}
Write $v(r):=\Vol B(x,r)$ and $f(r):=r^{-4}v(r)$.
A cutoff supported in $B(x,r)$, equal to one on $B(x,r/2)$,
with $|\nabla\eta|\le2/r$, gives
$v(r/2)^{1/2}\le5Ar^{-2}v(r)$ for $r\le V^{1/4}$.
Thus $f(r)\ge(20A)^{-1}f(r/2)^{1/2}$.
Iteration and $f(s)\to\omega_4$ at a smooth centre give
$\kappa:=(20A)^{-2}$.
Disjoint balls of radius $R:=V^{1/4}$, centred $2R$ apart
along a diameter-realizing geodesic, give
$d_0:=2\kappa^{-1}V^{1/4}$.
Likewise, disjoint $\varepsilon/2$ balls centred at a maximal
$\varepsilon$-separated set give
\begin{equation}\label{eq:coveringnumber}
 N(\varepsilon)\le\frac{16V}{\kappa\varepsilon^4},
 \qquad 0<\varepsilon\le2V^{1/4}.
\end{equation}
These covering and diameter bounds imply GH precompactness
\cite{BBI}.
\end{proof}

\begin{lemma}[Vanishing dissipation]\label{lem:tail}
For every $a<b$ and every $t_j\to\infty$,
\begin{equation}\label{eq:tail}
\int_a^b D(t_j+\tau)\,d\tau
=\tfrac12\bigl(\Ca(t_j+a)-\Ca(t_j+b)\bigr)\longrightarrow0.
\end{equation}
\end{lemma}
\begin{proof}
The nonnegative decreasing function $\Ca$ has a limit
$\Ca_\infty$. Integrating \eqref{eq:dissipation} gives
\begin{equation}\label{eq:totaldissipation}
 \Ca(t_2)+2\int_{t_1}^{t_2}D(t)\,dt=\Ca(t_1),
 \qquad 0\le t_1<t_2.
\end{equation}
In particular, $\int_0^\infty D(t)\,dt
=\frac12(\Ca(0)-\Ca_\infty)<\infty$.
Taking $t_1=t_j+a$ and $t_2=t_j+b$ proves \eqref{eq:tail}.
\end{proof}

We also record the rescaling relevant to point selection. We will use the Calabi rescaling:
\begin{equation}\label{Calabi rescaling}
g_j(\tau):=r_j^{-2}g(t_j+r_j^4\tau).
\end{equation}
Then the scalar curvature scales by $r_j^2$, the
scalar-curvature gradient by $r_j^4$, and volume by $r_j^{-4}$.
Consequently,
\begin{equation}\label{eq:scaleddissipation}
 \begin{aligned}
 D_j(\tau)&:=\int_M|\Dcal_{g_j(\tau)}S_{g_j(\tau)}|^2\,dV_{g_j(\tau)}
   =r_j^4D(t_j+r_j^4\tau),\\
 \int_{-L}^0D_j(\tau)\,d\tau
   &=\int_{t_j-Lr_j^4}^{t_j}D(t)\,dt.
 \end{aligned}
\end{equation}
Thus the spacetime dissipation is invariant under the
rescaling in real dimension four. It tends to zero whenever
$t_j-Lr_j^4\to\infty$, by \eqref{eq:totaldissipation}. For shrinking intervals at finite times, the same conclusion follows
from absolute continuity of the integral in \eqref{eq:totaldissipation};
see Section~\ref{sec:fixedclass}.

\subsection{Pointed \texorpdfstring{$\hat C^\infty$}{hat C-infinity}-Cheeger--Gromov convergence}
\label{sec:pointedconvergence}

We recall the pointed-$\hat C^\infty$-Cheeger--Gromov convergence of
Li--Li--Wang \cite[Definition~7.1 and Remarks~7.3, 7.6]{LLW},
in the K\"ahler orbifold setting used here. We abbreviate it as
\emph{pointed $\hat C^\infty$-CG convergence}.
The hat allows singular limits: the whole spaces converge in the
pointed GH sense, while the metrics converge smoothly on the regular
part through compatible comparison maps. We also require convergence
of the complex structures. The approximating spaces and the
basepoints may be singular.

\begin{definition}[Pointed-$\hat C^\infty$-Cheeger--Gromov convergence]
\label{def:modelhat}
Let $d_j,d$ be the length distances of complete pointed K\"ahler
orbifolds $(X_j,h_j,J_j,p_j)$ and $(X,h,J,p)$. We say that the
former converge to the latter in the pointed
$\hat C^\infty$-Cheeger--Gromov topology if, for every $L>0$ and
compact $K\Subset\mathcal R(X)\cap B_d(p,L)$, there are numbers
$\delta_j\to0$ and maps $\psi_j:B_d(p,L)\to X_j$ satisfying
\begin{align}
 d_j(\psi_j(p),p_j)&<\delta_j,\label{eq:modelGHbase}\\
 |d_j(\psi_j(x),\psi_j(y))-d(x,y)|&<\delta_j,
 \qquad x,y\in B_d(p,L),\label{eq:modelGHdist}\\
 B_{d_j}(p_j,L-\delta_j)&\subset
 N_{\delta_j}^{d_j}\bigl(\psi_j(B_d(p,L))\bigr),\label{eq:modelGHcover}
\end{align}
and smooth embeddings $\Phi_j:U\to\mathcal R(X_j)$, on a
neighborhood $U$ of $K$, such that
\begin{equation}\label{eq:modelhatcompatible}
 \psi_j|_K=\Phi_j|_K
\end{equation}
and
\begin{equation}\label{eq:modelhatsmooth}
 \Phi_j^*h_j\to h,\qquad \Phi_j^*J_j\to J
 \quad\hbox{in }C^\infty(K).
\end{equation}
Here $N_{\delta_j}^{d_j}$ denotes the metric neighborhood of the
indicated set. The maps and embeddings are required for all
sufficiently large indices of the given sequence. In a proof of
subsequential compactness, a countable exhaustion and diagonal
selection produce one subsequence on which these requirements hold
for every compact regular domain.
\end{definition}

Conditions \eqref{eq:modelGHbase}--\eqref{eq:modelGHcover} express
pointed GH convergence. The maps $\psi_j$ need not
be continuous outside $K$. Condition~\eqref{eq:modelhatcompatible}
requires that their restriction to $K$ be the same identification
used for smooth tensor convergence. The embeddings $\Phi_j$ need
not be holomorphic for the original complex structures; the
convergence of $\Phi_j^*J_j$ records the complex geometry in a fixed
smooth domain.

Finite order versions of these conditions, with embeddings defined near
$\overline U$ and estimates on $K\Subset U\Subset\mathcal R(X)$,
specify the neighbourhoods of this topology.

This compatibility is the starting point for the spacetime version
in Definition~\ref{def:topology}. There the same points must be
compared throughout a time interval, using comparison maps
independent of time.

\subsection{Extremal metrics and their flows}
\label{sec:extremal-limits}

All orbifolds below are connected and effective, with finitely many
isolated quotient singularities. A neighbourhood of $q$ is
$B^4/\Gamma_q$, where $\Gamma_q$ acts freely on $S^3$ and the
metric has a smooth invariant lift. Set
\[
 \mathcal S(X):=\operatorname{Sing}X=\{q:\Gamma_q\ne\{1\}\},
 \qquad \mathcal R(X):=X\setminus\mathcal S(X),
\]
and abbreviate these sets by $\mathcal S,\mathcal R$ when $X$ is fixed.
Thus a point is regular precisely when its tangent cone is Euclidean.
If $X$ is complete, $\mathcal R$ is geodesically convex:
a minimizing geodesic through a singular point would give two
orbits at distance $\pi$ in $S^3/\Gamma_q$, forcing every group
element to fix a unit vector. This contradicts the nontrivial free
action. The intrinsic metric completion of $\mathcal R$ is therefore $X$.

For a multifold, a punctured neighbourhood of $q$ is a disjoint
union of finitely many branches, each completing to $B^4/\Gamma_{q,\ell}$.
A point with several branches is a joining point; a multifold with
one branch everywhere is an orbifold. Its singular set includes
joining points and nontrivial quotient points. Choose a finite
$\mathcal E\supset\operatorname{Sing}X$; for convergence limits,
include all curvature-concentration points, even those smooth
after completion. The normalization is
\[
 X\setminus\mathcal E=\bigsqcup_{a=1}^NC_a,
 \qquad \widetilde X:=\bigsqcup_{a=1}^N\widetilde X_a,
\]
where $\widetilde X_a$ is the intrinsic length-metric completion
of $C_a$, a complete connected orbifold.

We next describe the evolution of an extremal limit. Set
\[
 d^c:=\frac{\sqrt{-1}}2(\bar\partial-\partial),
 \qquad dd^c=\ddc,
 \qquad d^cf=-\frac12df\circ J.
\]
For a K\"ahler metric $h$ with K\"ahler form $\omega$, set
$W(h):=\frac12\nabla_hS_h$.
For every real tangent vector $Y$,
\[
(\iota_{W(h)}\omega)(Y)
=-h(W(h),JY)
=-\frac12dS_h(JY)
=d^cS_h(Y).
\]
Since $d\omega=0$, Cartan's formula gives
$
\mathcal L_{W(h)}\omega=dd^cS_h.
$
If $h$ is extremal, then
$\bar\partial\nabla_h^{1,0}S_h=0$.
Thus $W(h)$ is real holomorphic and
$\mathcal L_{W(h)}J=0$.

Consequently, along a Calabi flow $h(\tau)$ on a fixed
complex manifold whose slices are all extremal, the
evolution equations take the form
\begin{equation}\label{eq:Lie}
 \partial_\tau\omega=dd^cS=\mathcal L_W\omega,
 \qquad
 \partial_\tau h=\mathcal L_Wh.
\end{equation}

\begin{lemma}[Preservation of the regular locus]\label{lem:regularorbits}
Let $(X,h_0,J)$ be a compact extremal K\"ahler orbifold surface
with isolated quotient singularities $\mathcal S$. The real holomorphic
field $W:=\tfrac12\nabla_{h_0}S_{h_0}$ generates a complete flow
$F_\tau$ of orbifold biholomorphisms. Each $F_\tau$ preserves
$\mathcal R(X)$ and fixes every point of $\mathcal S$.
For every nonempty compact $K\subset \mathcal R(X)$ and $0<L<\infty$,
one has
\[
 \{F_\tau(x):x\in K,\ |\tau|\le L\}\Subset \mathcal R(X).
\]
In particular, if $\mathcal S\ne\varnothing$, then
\[
 \min_{\substack{x\in K\\ |\tau|\le L}}
 d_{h_0}(F_\tau(x),\mathcal S)>0.
\]
This positive separation may depend on $h_0$, $K$, and $L$.
For $h_\tau:=F_\tau^*h_0$, the intrinsic curvature derivatives
are bounded uniformly on $X\times\R$, with bounds depending
only on $h_0$ and the derivative order.
\end{lemma}
\begin{proof}
Extremality makes $W$ a smooth real holomorphic orbifold vector
field. Compactness of $X$ makes it complete. Its flow consists
of orbifold biholomorphisms, which preserve the isotropy groups
up to isomorphism and hence preserve $\mathcal R(X)$ and $\mathcal S$.
Since $\mathcal S$ is finite and $F_0$ is the identity, continuity
implies that $F_\tau$ fixes each singular point.

The displayed orbit set is the continuous image of the compact
set $[-L,L]\times K$ and lies in $\mathcal R(X)$. Its distance
from the closed singular set is therefore positive when that
set is nonempty. Finally,
$F_\tau:(X,h_\tau)\to(X,h_0)$ is an isometry, so
\[
 |\nabla^m\Rm(h_\tau)|_{h_\tau}(x)
 =|\nabla^m\Rm(h_0)|_{h_0}(F_\tau(x))
 \le\sup_X|\nabla^m\Rm(h_0)|_{h_0}<\infty.
\]
The last bound follows from orbifold smoothness and compactness,
for every $m\ge0$.
\end{proof}

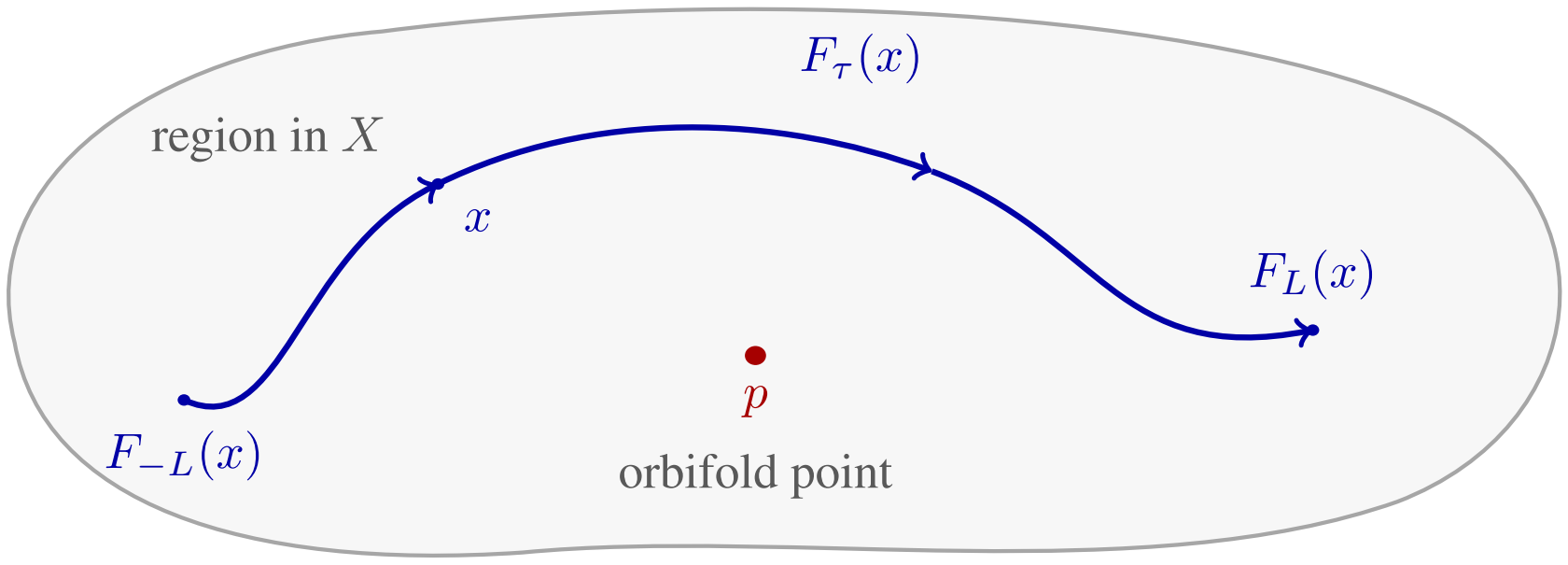
\begin{figure}[htbp]
\centering
\begin{tikzpicture}[x=1cm,y=0.9cm,font=\small]
 \path[fill=black!3,draw=black!35,thick]
 (0.4,1.7) .. controls (0.1,3.0) and (1.4,4.0) .. (3.0,4.15)
 .. controls (5.4,4.5) and (8.9,4.35) .. (10.5,3.5)
 .. controls (11.75,2.8) and (11.6,1.05) .. (10.2,0.45)
 .. controls (8.4,-0.25) and (5.9,0.25) .. (4.0,0.05)
 .. controls (2.0,-0.1) and (0.6,0.45) .. (0.4,1.7) -- cycle;
 \node[black!65] at (2.2,3.3) {region in $X$};
 \draw[blue!65!black,very thick,->]
 (1.6,1.25) .. controls (2.3,0.9) and (2.4,2.4) .. (3.4,2.95);
 \draw[blue!65!black,very thick,->]
 (3.4,2.95) .. controls (4.5,3.55) and (5.8,3.5) .. (6.9,3.05);
 \draw[blue!65!black,very thick,->]
 (6.9,3.05) .. controls (8.1,2.55) and (8.2,1.5) .. (9.6,1.8);
 \foreach \x/\y in {1.6/1.25,3.4/2.95,9.6/1.8} {
   \fill[blue!65!black] (\x,\y) circle (0.045);
 }
 \node[blue!65!black,below,inner sep=6pt] at (1.6,1.25) {$F_{-L}(x)$};
 \node[blue!65!black,below right,inner sep=5pt] at (3.4,2.95) {$x$};
 \node[blue!65!black,above,inner sep=6pt] at (6.4,3.5) {$F_\tau(x)$};
 \node[blue!65!black,above,inner sep=6pt] at (9.6,1.8) {$F_L(x)$};
 \fill[red!65!black] (5.65,1.6) circle (0.075);
 \node[red!65!black,below,inner sep=6pt] at (5.65,1.6) {$p$};
 \node[black!65] at (5.65,0.65) {orbifold point};
\end{tikzpicture}
\caption{A regular flow line avoids the orbifold point $p$.}
\label{fig:orbifoldtrajectory}
\end{figure}

Figure~\ref{fig:orbifoldtrajectory} gives a spatial illustration
of Lemma~\ref{lem:regularorbits}. The blue curve is a finite
segment of a regular trajectory in $X$; the red point is a fixed
orbifold point. The drawing is schematic.

\begin{lemma}[Evolution of an extremal limit]\label{lem:transport}
Let $h(\tau)$ be a smooth Calabi flow on a complex manifold,
defined on an interval containing zero, and suppose every slice is
extremal. Then the vector field
$W(h(\tau))=\frac12\nabla_{h(\tau)}S_{h(\tau)}$ is independent
of time. If $F_\tau$ is its local holomorphic flow, then
\begin{equation}\label{eq:localextremalpullback}
 h(\tau)=F_\tau^*h(0)
\end{equation}
wherever the relevant integral curve segments remain in the domain
of the given solution.

If this domain is the entire regular locus of a compact complex
orbifold surface $X$ with isolated quotient singularities and $h(0)$
extends smoothly to $X$, formula~\eqref{eq:localextremalpullback}
holds on the entire regular locus throughout the given time interval
and defines an eternal Calabi flow on $X$.
\end{lemma}
\begin{proof}
If $F$ is a biholomorphism, then
$S_{F^*h}=S_h\circ F$ and
$
 \nabla_{F^*h}(S_h\circ F)=F^*(\nabla_hS_h).
$
Hence
\begin{equation}\label{eq:Wnaturality}
 W(F^*h)=F^*W(h),
 \qquad
 (F^*W(h))_p=(dF_p)^{-1}W(h)_{F(p)}.
\end{equation}
This gives naturality under biholomorphic
pullbacks.

Fix an arbitrary time $\tau_0$ and temporarily write
$h:=h(\tau_0)$ and $W:=W(h(\tau_0))$. Let $\Phi_s$ be the local
holomorphic flow of this fixed vector field. Differentiating
\eqref{eq:Wnaturality} at $s=0$ gives
\[
 \mathrm dW_h[\mathcal L_Wh]
 =\left.\frac{d}{ds}\right|_{s=0}W(\Phi_s^*h)
 =\left.\frac{d}{ds}\right|_{s=0}\Phi_s^*W
 =\mathcal L_WW=[W,W]=0.
\]
By \eqref{eq:Lie},
\[
 \left.\partial_\tau W(h(\tau))\right|_{\tau=\tau_0}
 =\mathrm dW_h[\partial_\tau h(\tau_0)]
 =\mathrm dW_h[\mathcal L_Wh]=0.
\]
Thus $W(h(\tau))$ is independent of time; denote it by $W_0$.
It follows that
\begin{equation}\label{eq:transportcancellation}
 \frac{d}{d\tau}\bigl(F_{-\tau}^*h(\tau)\bigr)
 =F_{-\tau}^*\bigl(\partial_\tau h-\mathcal L_{W_0}h\bigr)=0.
\end{equation}
Evaluating at $\tau=0$, we get
\eqref{eq:localextremalpullback}.
%More explicitly, on $U\Subset U^+$ it holds for all sufficiently small $|\tau|$ for which the integral-curve segments beginning in $\overline U$ remain in $U^+$, where the solution is defined. This local version will be used in Section~\ref{sec:upgrade}.

Suppose next that $h(0)$ extends smoothly to $X$ as in the
statement. Extremality extends by continuity. Lemma~\ref{lem:regularorbits}
gives a complete holomorphic flow preserving $\mathcal R(X)$.
Thus every finite characteristic segment remains in the domain
of the solution, and \eqref{eq:transportcancellation} proves
\eqref{eq:localextremalpullback} globally there. By
\eqref{eq:Wnaturality}, $F_\tau^*W_0=W_0$, so
$F_\tau^*h(0)$ satisfies the Calabi flow for all time.
\end{proof}

%%%%%%%%%%%%%%%%%%%%%%%%%%%%%%%%%%%%%%%%%%%%%%%%%%%%%%%%%%%%%%%%%%%%%%%%%%%%%%%%%%%%%%%%%%%%%%%%%%%%%%%%%%%%%%%%%%%%%%%
\subsection{Convergence on the regular part}
\label{Convergence on the regular part}

Here convergence is assumed, rather than proved. We explain
how Sobolev inequalities, volume measures, and parallel complex
structures pass to the limit. These facts are used for the
stationary models in Section~\ref{sec:stationary} and for
the flow slices in Section~\ref{sec:slicecompactness}.

Throughout this subsection, $(X_j,h_j,p_j)$ are complete pointed
four-manifolds or orbifolds converging in pointed GH sense to a
proper length space $(X,d,p)$. Off a finite $\mathcal E$, the limit
has a smooth metric $h$ whose length distance agrees locally with $d$.
Choose Borel GH approximations
$f_j:B_{h_j}(p_j,R_j)\to X$, with $R_j\to\infty$, distortion
$\varepsilon_j\to0$, $d(f_j(p_j),p)\le\varepsilon_j$ and
$\varepsilon_j$-dense image in $B_d(p,R_j-1)$.
For $K\Subset U\Subset X\setminus\mathcal E$, assume smooth embeddings
$\Phi_j:U\to X_j$ satisfy
\[
 \Phi_j^*h_j\to h\quad\hbox{smoothly},\qquad
 \sup_{x\in K}d(f_j(\Phi_j(x)),x)\to0.
\]
We use \emph{coverage} to mean that, for every $K\Subset U$,
any sequence $x_j\in X_j$ with GH limit in $K$ eventually lies in
$\Phi_j(U)$, also after passage to a subsequence.
Under the preceding assumptions, the charts have this property.
Indeed, choose $K\Subset U'\Subset U$.
Smooth metric convergence and a first-exit argument give a uniform
positive radius for balls about $\Phi_j(K)$ contained in $\Phi_j(U')$.
If $x_j$ has GH limit $x\in K$, compatibility and vanishing distortion
give $d_{h_j}(x_j,\Phi_j(x))\to0$. Hence $x_j\in\Phi_j(U')$ for large $j$.
In particular, $f_j^{-1}(K)\subset\Phi_j(U)$ for large $j$.

Cutoffs of vanishing capacity extend Sobolev inequalities across
isolated singularities; see \cite[Proof of Theorem 5.1]{MR1001844}
\cite[Lemma~2.4]{CW}. We retain the lower-order coefficient
and apply the same argument on each branch.
\begin{lemma}
\label{lem:sobolev-passage}
Suppose \eqref{eq:Sob} holds on $(X_j,h_j)$ with $(A,b_j)$,
where $0\le b_j\to b<\infty$, and
$\Vol_h B_h(q,r)\le C_qr^4$ for $q\in\mathcal E$ and sufficiently
small $r$. Then:
\begin{enumerate}[label=\textup{(\roman*)}]
\item Each $q\in\mathcal E$ has zero $W^{1,2}$ capacity. A cutoff
$\chi_\rho$ vanishing on $B(q,\rho)$ and equal to one outside
$B(q,2\rho)$ satisfies $|d\chi_\rho|\le C\rho^{-1}$ and
\begin{equation}\label{eq:pointcapacity}
 \int|d\chi_\rho|^2\,dV_h
 \le C\rho^{-2}\Vol_h B_h(q,2\rho)\le C\rho^2\longrightarrow0.
\end{equation}
\item Every compactly supported
$u\in W^{1,2}_{\rm loc}(X\setminus\mathcal E)$ with finite
$\int(u^2+|du|^2)\,dV_h$ satisfies
\begin{equation}\label{singular sobolev}
 \|u\|_{L^4(h)}^2
 \le A\left(\int|du|^2\,dV_h+b\int u^2\,dV_h\right).
\end{equation}
\item If $\operatorname{supp}u\subset B_h(x,r)$ and $\Vol_h B_h(x,r)\le Cr^4$,
then
\begin{equation}\label{eq:localhomogeneoussobolev}
 \|u\|_{L^4(h)}^2\le2A\int|du|^2\,dV_h
 \quad\hbox{when }AbC^{1/2}r^2\le\tfrac12.
\end{equation}
For $b=0$, the constant remains $A$.
\item The same statements hold on each completed branch, for
functions supported away from its outer boundary, by zero extension
to the other branches. No bound on the number of branches is needed.
\end{enumerate}
\end{lemma}
\begin{proof}
The cutoff estimate \eqref{eq:pointcapacity} follows from the volume
bound. Smooth convergence first gives \eqref{singular sobolev}
for smooth functions supported away from $\mathcal E$.
For a general $u$, truncate it to
$u_m:=\max\{-m,\min\{u,m\}\}$ and let $\chi_\rho$ be the product
of the cutoffs at the relevant exceptional points. Then
\[
 d(\chi_\rho u_m)-du_m
   =(\chi_\rho-1)du_m+u_m\,d\chi_\rho\longrightarrow0
   \quad\hbox{in }L^2
\]
by absolute continuity and \eqref{eq:pointcapacity}.
Also $\chi_\rho u_m\to u_m$ in $L^2$ and $L^4$.
Approximate $\chi_\rho u_m$ smoothly away from $\mathcal E$,
apply the inequality, then let $\rho\downarrow0$ and $m\to\infty$.
Monotone convergence and $|du_m|\le|du|$ prove
\eqref{singular sobolev}.
The local estimate follows by absorbing
$b\int u^2\le bC^{1/2}r^2\|u\|_4^2$.
For a single branch, apply the same argument to the zero extension
to all other branches.
\end{proof}

\begin{lemma}[Convergence of volume measures]
\label{lem:volume-passage}
Set $\mu_j:=(f_j)_*(dV_{h_j}|_{B_{h_j}(p_j,R_j)})$.
Assume that for each $R<\infty$ there are $C_R,r_R>0$ such that
\[
 \Vol_{h_j}B_{h_j}(z,r)\le C_Rr^4,
 \qquad z\in B_{h_j}(p_j,R),\quad 0<r<r_R,
\]
for all large $j$. Then $dV_h$, extended by zero on $\mathcal E$,
is Radon and $\mu_j\rightharpoonup dV_h$ on bounded regions.
Thus $\mu_j(B)\to dV_h(B)$ for every bounded Borel continuity set $B$.
If $f_j(x_j)\to x$ and $r_j\to r>0$, then
$\Vol_{h_j}B_{h_j}(x_j,r_j)\to dV_h(B_d(x,r))$ whenever
$dV_h(\partial B_d(x,r))=0$.
\end{lemma}

\begin{proof}
Set $\nu:=dV_h$ off $\mathcal E$ and $\nu(\mathcal E):=0$.
Coverage, compatibility and smooth convergence give, for
$\varphi\in C_c(X\setminus\mathcal E)$ supported in a regular chart,
\begin{equation}\label{eq:regular-volume-passage}
 \int_X\varphi\,d\mu_j
 =\int_U\varphi\circ f_j\circ\Phi_j\,dV_{\Phi_j^*h_j}
 \longrightarrow\int_X\varphi\,d\nu.
\end{equation}
A partition of unity gives the same conclusion for arbitrary
compact regular support. Choose $q_j$ with $f_j(q_j)\to q\in\mathcal E$.
For small fixed $r$, the GH distortion and source upper volume bound give
\begin{equation}\label{eq:exceptional-source-volume}
 \mu_j(B_d(q,2r))
 \le\Vol_{h_j}B_{h_j}(q_j,3r)\le C_R(3r)^4.
\end{equation}
Taking the supremum in \eqref{eq:regular-volume-passage} over
$0\le\varphi\le1$ supported compactly in
$B_d(q,2r)\setminus\mathcal E$ yields
\begin{equation}\label{eq:exceptional-limit-volume}
 \nu(B_d(q,2r))\le C_R(3r)^4.
\end{equation}
Thus $\nu$ is locally finite. For $\varphi\in C_c(X)$, remove
the $2r$ neighbourhood of $\mathcal E$ by a cutoff and use
\eqref{eq:regular-volume-passage} on its complement.
The omitted integrals are bounded by
$C\|\varphi\|_\infty\#\mathcal E\,r^4$ on both sides.
Letting $r\downarrow0$ proves $\mu_j\rightharpoonup\nu$.
Portmanteau gives convergence on bounded continuity sets
\cite{Billingsley}. Finally, if $f_j(x_j)\to x$ and $r_j\to r$,
GH distortion gives $\delta_j\downarrow0$ with
\[
 f_j^{-1}(B_d(x,r-\delta_j))
 \subset B_{h_j}(x_j,r_j)
 \subset f_j^{-1}(B_d(x,r+\delta_j)).
\]
The two bounds converge to $\nu(B_d(x,r))$ at a continuity radius.
\end{proof}

\begin{lemma}[Convergence of complex structures]
\label{lem:parallel-complex}
Let $h_j(\tau)\to h(\tau)$ smoothly locally on $U\times I$,
where $U$ is a smooth four-manifold and $I$ is a compact interval,
possibly a single time. Suppose $J_j$ is independent of $\tau$,
$J_j^2=-\mathrm{Id}$, and $J_j$ is $h_j(\tau)$-orthogonal and
parallel for every $\tau$. A subsequence converges smoothly locally
to a time-independent $J$, and every $(h(\tau),J)$ is K\"ahler.
This applies in particular to time-independent comparison charts.

If $h$ is smooth on $B^4$, every orthogonal parallel complex
structure on $B^4\setminus\{0\}$ extends uniquely and smoothly
across $0$. A finite group preserving $h$ and the punctured
complex structure preserves the extension; its tangent action
at $0$ is conjugate into $U(2)$.
\end{lemma}
\begin{proof}
At a fixed time, orthogonality bounds $J_j$ and
$\partial J_j=\Gamma(h_j)*J_j$ bounds every derivative inductively.
A diagonal subsequence converges; all algebraic and parallel
identities persist at every time, so the limit is K\"ahler.

On the punctured ball, the same equation bounds $J$ and $\partial J$.
Joining points approaching $0$ by paths of vanishing length gives
a unique limit $J(0)$. If $P_x$ is radial parallel transport, then
$J(x)=P_xJ(0)P_x^{-1}$, which is smooth also at $0$.
Parallelness and the algebraic identities persist by continuity,
and uniqueness preserves the group action.
\end{proof}

%%%%%%%%%%%%%%%%%%%%%%%%%%%%%%%%%%%%%%%%%%%%%%%%%%%%%%%%%%%%%%%%%%%%%%%%%%%%%%%%%%%%%%%%%%%%%%%%%%%%%%%%%%%%%%%%%%%%%%%

%%%%%%%%%%%%%%%%%%%%%%%%%%%%%%%%%%%%%%%%%%%%%%%%%%%%%%%%%%%%%%%%%%%%%%%%%%%%%%%%%%%%%%%%%%%%%%%%%%%%%%%%%%%%%%%%%%%%%%%%%%%%%%%%%%%%%%%%%%%%%%%%%%%%%%%%%%%%%%%%%%%%%%%%%%%%%%%%%%%%%%%%%%%%%%%%%%%%%%%%%%%%%%%%%%%%%%%%%%%%%%
% BEGIN proofs/fixed_class_bubble_appendix.tex
\subsection{Shrinking-scale limits in a fixed K\"ahler class}
\label{sec:fixedclass}
\label{sec:bubbles}\label{sec:hodgeperiods}

We record geometric and cohomological restrictions on shrinking-scale
limits. The scalar-flatness and ALE statements assume a completed
extremal orbifold; the period argument needs only regular convergence
near a cycle. Existence of the limits is established in
Section~\ref{sec:slicecompactness}.

For $r_j\downarrow0$, set
\begin{equation}\label{eq:scaling}
\widehat g_j(\tau):=r_j^{-2}g(t_j+r_j^4\tau).
\end{equation}
In dimension four, curvature energy and centered Calabi energy are invariant under this scaling. The dissipation rescales by
\[
\widehat D_j(\tau)=r_j^4D(t_j+r_j^4\tau).
\]
If both endpoints of a rescaled fixed interval approach the same endpoint of the original flow, where its monotone energy has a limit, then
\begin{equation}\label{eq:scaled-drop}
\int_{-L}^0\widehat D_j(\tau)\,d\tau
=\tfrac12\bigl[\Ca(t_j-Lr_j^4)-\Ca(t_j)\bigr]\longrightarrow0.
\end{equation}
The rescaled Sobolev inequality becomes
\begin{equation}\label{eq:scaledSob}
\norm{u}_{L^4(\widehat g_j)}^2
\le A\left(\int|\nabla u|^2\,dV_{\widehat g_j}
+r_j^2V^{-1/2}\int u^2\,dV_{\widehat g_j}\right).
\end{equation}
Thus a regular bubble inherits a homogeneous Sobolev inequality.
The compact limit of unrescaled time translates retains the lower-order term.

The estimates for extremal metrics of
\cite[Theorem~3.1]{CW} give decay of the scalar-curvature Killing field
at infinity. Its Jacobi equation along a ray then forces the field
to vanish.

\begin{lemma}[Scalar-flatness of complete extremal bubbles]\label{lem:scalarflat}
Let $(Y,h,J)$ be a connected complete noncompact extremal K\"ahler orbifold of real dimension four with finitely many isolated quotient points. Assume
\[
\norm{u}_{L^4(h)}^2\le A\int_Y|\nabla u|^2\,dV_h
\quad(u\in C_c^\infty(Y)),\qquad
\int_Y|\Rm_h|^2\,dV_h<\infty.
\]
Then $S_h\equiv0$.
\end{lemma}
\begin{proof}
We fix $o\in Y$. For $r:=d(o,y)$ large, a ball of radius $r/4$ centered at $y$ lies outside a compact set containing all orbifold points. Its curvature energy tends uniformly to zero. The interior estimate for extremal metrics \cite[Theorem~3.1]{CW} gives
\begin{equation}\label{eq:taildecay}
|\nabla^k\Rm|(y)=o(r^{-2-k})
\end{equation}
for every fixed $k$. In particular the Killing field $K:=J\nabla S$ satisfies $|K|=o(r^{-3})$ and $|\nabla K|=o(r^{-4})$.

We take a minimizing ray $\gamma$ from any regular point $x$. It cannot pass through a nontrivial isolated quotient point. A minimizing path through the cone vertex would require two directions at distance $\pi$ in its link, whereas a nontrivial free spherical quotient has diameter strictly less than $\pi$. Along $\gamma$, in a parallel frame, $y(s):=K(\gamma(s))$ solves
\[
y''(s)=-A(s)y(s),\qquad |A(s)|=o(s^{-2}).
\]
The decay of $y,y'$ permits integration twice from infinity:
\[
y(s)=-\int_s^\infty(u-s)A(u)y(u)\,du.
\]
For large $R$, we put $M_R:=\sup_{s\ge R}s^3|y(s)|<\infty$. If $|A(u)|\le\eps u^{-2}$ for $u\ge R$, then
\[
s^3|y(s)|\le\eps M_Rs^3\int_s^\infty(u-s)u^{-5}\,du
=\frac{\eps}{12}M_R.
\]
We choose $\eps<12$ and take the supremum. Then $M_R=0$, and ODE uniqueness gives $K(x)=0$. Since the initial point was arbitrary, the scalar curvature is constant. Equation~\eqref{eq:taildecay} makes this constant zero.
\end{proof}

If a multifold occurs, its normalisation separates the local branches. By Lemma~\ref{lem:sobolev-passage}, the homogeneous
Sobolev inequality passes with the same constant to every
normalised component. No such component is compact. Approximating its constant function by these cutoffs would contradict the homogeneous inequality. Thus Lemma~\ref{lem:scalarflat} applies to every normalised component once smooth orbifold removal is justified.

We apply the finite-end theorem and the exterior decay improvement of
Tian--Viaclovsky \cite[Theorems~1.1, 2.5 and 5.1]{TV}.
The orbifold points lie in a compact core, so the exterior argument
applies on smooth ends.

\begin{proposition}[ALE conclusion after completion]\label{prop:ALE}
Under the hypotheses of Lemma~\ref{lem:scalarflat}, $Y$ has finitely many ALE ends. The ends have positive differentiable ALE order. In the scalar-flat K\"ahler case the usual exterior estimates give any order $\tau<2$.
\end{proposition}
\begin{proof}
Equation~\eqref{eq:taildecay} gives $|\Rm|=o(r^{-2})$, and finite curvature energy implies $\int|\Ric_-|^2<\infty$. The homogeneous Sobolev inequality and these bounds are the hypotheses of the finite-end theorem, including its finite-orbifold extension \cite[Theorems 1.1 and 2.5]{TV}. It gives ALE order zero. By Lemma~\ref{lem:scalarflat}, the metric is scalar-flat. With its complex orientation it is anti-self-dual. The exterior decay improvement \cite[Theorem 5.1]{TV} applies on each end. For an orbifold, choose cutoffs outside a compact core containing all
quotient points. The equations and coordinate estimates on each end
are unchanged, so the decay improvement applies there. Parallelness of $J$ gives the corresponding complex-structure decay in the ALE coordinates.
\end{proof}

Even after one-endedness is established
(Lemma~\ref{lem:oneend}), several normalized components could
still meet at a multifold point. Excluding such branching requires
a separate argument.

% BEGIN proofs/null_periods_hodge_section.tex

The K\"ahler Hodge decomposition \cite[Section~2.G]{Besse} identifies the
self-dual and anti-self-dual cohomology spaces from the complex
structure and the K\"ahler class. The resulting expression for the
harmonic norm is independent of the representative metric.

\begin{lemma}[Fixed cohomological Hodge norm]\label{lem:fixedHodgeNorm}
Let $(M,J)$ be a closed connected complex surface and let $\Omega$ be a fixed real
K\"ahler class.  For $a\in H^2(M;\R)$, let $\eta_{a,g}$ denote its harmonic
representative for a K\"ahler metric $g$ whose K\"ahler form represents
$\Omega$.  Then $\|\eta_{a,g}\|_{L^2(g)}$ is independent of $g$.
\end{lemma}
\begin{proof}
We write $Q(a,b):=\int_M a\wedge b$ for the intersection pairing.  The K\"ahler
identities and Lefschetz decomposition identify the self-dual and
anti-self-dual harmonic cohomology subspaces with the fixed spaces
(see \cite[Section~2, Example, p.~654]{LeBrunPolarized}):
\begin{align*}
P&:=\R\Omega\oplus
  \left(\bigl(H^{2,0}(M)\oplus H^{0,2}(M)\bigr)\cap H^2(M;\R)\right),\\
N&:=\{a\in H^{1,1}(M;\R):Q(a,\Omega)=0\}.
\end{align*}
Contraction of a harmonic $(1,1)$ form with the K\"ahler
form is a harmonic function, hence constant.  If its cohomology class is
$Q$-orthogonal to $\Omega$, that constant is zero.  The form is therefore
primitive and, on a complex surface, anti-self-dual.  Real forms of type
$(2,0)+(0,2)$ and the K\"ahler form are self-dual.  Thus the decomposition depends only on $(J,\Omega)$.

We write $a=a_++a_-$ in the $Q$-orthogonal splitting $P\oplus N$.
The Hodge star identities \cite[13.5 and 13.7]{Besse} give
\begin{equation}\label{eq:fixedHodgeNorm}
\|\eta_{a,g}\|_{L^2(g)}^2
 =Q(a_+,a_+)-Q(a_-,a_-),
\end{equation}
which is independent of the representative metric.  The harmonic
representative itself need not be fixed.
\end{proof}

Compactly supported Poincar\'e dual forms convert the global
$L^2$ bounds into bounds for periods. Integrality then restricts
the transported homology classes to a finite set.

\begin{lemma}[Finitely many transported homology classes]
\label{lem:transportedLattice}
Let $g_j$ be metrics on a fixed closed oriented four-manifold $M$.
We choose a basis $e_1,\ldots,e_b$ of $H^2(M;\Z)/\mathrm{torsion}$, and let
$\alpha_{a,j}$ be its $g_j$-harmonic representatives.  Suppose
\[
\sup_j\|\alpha_{a,j}\|_{L^2(g_j)}<\infty
\qquad (1\leq a\leq b).
\]
Let $U$ be a smooth four-manifold, $\lambda_j>0$, and let embeddings
$\Phi_j:U\to M$ satisfy
$h_j:=\lambda_j\Phi_j^*g_j\to h$ smoothly on compact subsets.
For every fixed compact integral two-cycle $\Sigma$ in $U$, the classes
$\Phi_{j*}[\Sigma]$ lie in a finite subset of
$H_2(M;\Z)/\mathrm{torsion}$.
\end{lemma}
\begin{proof}
Choose an orientation on $U$ and a closed two-form
$\theta\in\Omega_c^2(U)$ representing the compactly supported
Poincar\'e dual of $[\Sigma]\in H_2(U;\R)$. Thus, for every closed
two-form $\beta$ on $U$,
\[
 \int_\Sigma\beta=\int_U\theta\wedge\beta.
\]
Set $\beta_{a,j}:=\Phi_j^*\alpha_{a,j}$. In real dimension four,
the $L^2$ norm of a two-form is invariant under constant metric
scaling, so
\[
 \|\beta_{a,j}\|_{L^2(U,h_j)}
 =\|\alpha_{a,j}\|_{L^2(\Phi_j(U),g_j)}
 \le\|\alpha_{a,j}\|_{L^2(M,g_j)}\le C_a.
\]
Here $C_a:=\sup_j\|\alpha_{a,j}\|_{L^2(M,g_j)}<\infty$ by hypothesis.
In our application to K\"ahler metrics on a fixed complex surface
in a fixed K\"ahler class, Lemma~\ref{lem:fixedHodgeNorm} supplies this
bound; in fact, the global norm is independent of $j$.
The support of $\theta$ is fixed and compact. Convergence of $h_j$
there gives $\|\theta\|_{L^2(U,h_j)}\le C_\Sigma$. Since
$\beta_{a,j}$ is closed, Cauchy--Schwarz yields
\[
 \left|\langle e_a,\Phi_{j*}[\Sigma]\rangle\right|
 =\left|\int_U\theta\wedge\beta_{a,j}\right|
 \le \|\theta\|_{L^2(U,h_j)}
       \|\beta_{a,j}\|_{L^2(U,h_j)}
 \le C_\Sigma C_a.
\]
These pairings are integers and, by the universal coefficient
theorem, determine the free homology class. Only finitely many
bounded integer coordinate vectors occur.
\end{proof}

Chen--Wang \cite[Proposition~4.4]{ChenWangFano} prove null periods
for K\"ahler--Ricci flow bubbles on Fano surfaces using the integral
class $c_1(M)$. See also Chen--He \cite{CH} for Calabi-flow bubbles.
The preceding Hodge and lattice lemmas extend this period argument
to an arbitrary fixed real K\"ahler class.

\begin{lemma}[Null periods on fixed-class bubbles]\label{lem:period}
Let $\omega_j$ be K\"ahler forms on a fixed closed connected complex surface $(M,J)$,
all representing one fixed real K\"ahler class $\Omega$.
Suppose $\lambda_j\to\infty$ and embeddings give smooth K\"ahler
convergence
\[
(\lambda_j\Phi_j^*g_j,\Phi_j^*J,\lambda_j\Phi_j^*\omega_j)
 \longrightarrow(h,I,\omega_h)
\]
on a smooth region $U$.  Then every compact integral two-cycle $\Sigma$
in $U$, in particular every compact closed oriented real surface, satisfies
\[
\int_\Sigma\omega_h=0.
\]
\end{lemma}
\begin{proof}
Lemma~\ref{lem:fixedHodgeNorm} supplies the harmonic $L^2$ bounds of
Lemma~\ref{lem:transportedLattice}.  Hence
$c_j:=\Phi_{j*}[\Sigma]$ ranges over only finitely many free integral
homology classes.  Since $[\omega_j]=\Omega$,
\[
\langle\Omega,c_j\rangle
 =\int_\Sigma\Phi_j^*\omega_j
 =\lambda_j^{-1}\int_\Sigma\lambda_j\Phi_j^*\omega_j
 \longrightarrow0.
\]
These numbers belong to a finite set, so they vanish for all sufficiently
large $j$.  Therefore
\[
\int_\Sigma\omega_h
 =\lim_j\int_\Sigma\lambda_j\Phi_j^*\omega_j
 =\lim_j\lambda_j\langle\Omega,c_j\rangle=0.
\]
Torsion has no effect on real periods.
\end{proof}

The convergence hypotheses in Lemma~\ref{lem:period} can be weakened.
In fact, with $(M,J)$ and $\Omega$ fixed, the proof only needs
uniform equivalence of the rescaled metrics to a fixed smooth metric
on the support of a compactly supported Poincar\'e dual form for each cycle, together
with convergence of the rescaled K\"ahler-form integrals on that cycle.
The first condition gives finitely many transported homology classes
by the proof of Lemma~\ref{lem:transportedLattice}; the second gives null periods.
The finite-order regular charts used later supply both properties.
The argument applies to real surfaces without restrictions on
holomorphicity or self-intersection, and uses neither scalar-flatness
nor ALE asymptotics.

The null-period conclusion also extends across isolated quotient
points. For an orbifold with finitely many such points, the links
are rational homology three-spheres. Their
real homology in degrees one and two vanishes.  Excision and the long
exact sequence therefore show that inclusion of the regular locus
induces an isomorphism on real degree-two homology.  All periods of the
completed orbifold K\"ahler form vanish as well.  We apply this argument
component by component if a multifold is first normalized. Any finitely
many provisionally removed smooth energy-concentration points can also
be avoided, since their links are three-spheres. Thus the argument does
not require smooth convergence across those points in advance.

% END proofs/null_periods_hodge_section.tex

% END proofs/fixed_class_bubble_appendix.tex

\section{Moduli of null scalar-flat ALE K\"ahler orbifolds}
\label{sec:stationary}

We prove compactness of the ALE models used to define the canonical
radius in Section~\ref{sec:scr}. We first bound volume growth and
exclude branching in the limits. We then construct finite bubble
trees, retaining both curvature bubbles and original orbifold points.
Nullness and rigidity at an asymptotically Euclidean end exclude
concentration at a regular limit point. This gives smooth convergence
on the entire regular locus, with charts covering the corresponding
regions of the approximating spaces.

%\paragraph{ALE space}
We initially allow an ALE space to have finitely many ends. On each
end $E$, there are coordinates
$\Psi: (\R^4\setminus B)/\Gamma\rightarrow E$, where $\Gamma\subset U(2)$ acts freely
on $S^3$, such that
\begin{equation}\label{eq:modelALEdefinition}
 |\partial^k(h-e)|+|\partial^k(J-J_0)|=O(r^{-\tau-k}),
 \qquad k\ge0,
\end{equation}
for some $\tau>0$.
We call $\tau$ the differentiable decay order of the end $E$.
Here $e$ and $J_0$ are the standard Euclidean
metric and complex structure, and the estimates are taken in the
lifted end coordinates. The coordinate systems, their starting
radii, and the constants in these estimates may depend on the
space. An end with trivial group is called asymptotically Euclidean (AE).
When $(h,J)$ is K\"ahler, we also write $\omega:=h(J\cdot,\cdot)$.

%Uniform volume estimates will be established independently of these choices. Lemma~\ref{lem:oneend} will show that the K\"ahler orbifolds under consideration have exactly one end.

We fix $A,\Lambda>0$ and $0<\alpha<\eps_{\rm reg}(A)$. The constant
$\eps_{\rm reg}$ is chosen below the critical regularity thresholds.
\begin{definition}[The model moduli]\label{def:modelmoduli}
The space $\mathcal A(A,\Lambda,\alpha)$ consists of pointed
biholomorphic isometry classes $(X,h,J,p)$ with the following
properties.
\begin{enumerate}[label=\textup{(\roman*)}]
\item $(X,h,J)$ is a complete noncompact scalar-flat K\"ahler ALE
orbifold. The base point $p$ may be singular.
\item For every smooth compactly supported orbifold function $u$,
\begin{equation}\label{eq:modelSob}
 \|u\|_{L^4(h)}^2\le A\int_X|du|_h^2\,dV_h.
\end{equation}
\item Its weighted energy, defined by the curvature integral plus
the weighted number of singular points, satisfies
\begin{equation}\label{eq:modelComplexity}
 \mathscr E_\alpha(X,h)
 :=\int_X|\Rm_h|^2\,dV_h+\alpha\,\#\Sing X\le\Lambda.
\end{equation}
\item Its K\"ahler form is \emph{null}, in the sense that
\begin{equation}\label{eq:modelNull}
 \int_\Sigma\omega=0
\end{equation}
for every compact smooth closed oriented real surface
$\Sigma\subset\mathcal R(X)$.
\end{enumerate}
\end{definition}
The number of singular points in \eqref{eq:modelComplexity} cannot
be replaced by the curvature integral: a flat cone $\R^4/\Gamma$ has zero curvature
energy, but its vertex must still be recorded when $\Gamma$ is
nontrivial.

The class
contains every constant rescaling $h\mapsto c^2h$ of each of its members, including
the flat quotients satisfying the stated bounds. Indeed,
all defining bounds are invariant under a constant rescaling with $J$ fixed: $\omega$ scales by $c^2$,
the two sides of \eqref{eq:modelSob} both scale by $c^2$, and
$\int|\Rm|^2dV$ is unchanged in dimension four. Nullness and the
number of quotient points are unchanged.

\begin{lemma}
\label{nullness condition}
Let $X$ be an oriented orbifold of real dimension four with
finitely many isolated quotient singularities, and let
$\mathcal F\subset X$ be a finite set containing
$\operatorname{Sing}X$. Then inclusion induces an isomorphism
\[
 H_2(X\setminus\mathcal F;\R)
 \xrightarrow{\;\simeq\;}
 H_2(X;\R).
\]
Consequently, a smooth closed orbifold two-form $\omega$ on
$X$ satisfies
\begin{align}\label{class vanishes}
 [\omega]=0\in H^2_{\mathrm{dR}}(X;\R)
\end{align}
if and only if the nullness condition \eqref{eq:modelNull} holds
for every compact smooth closed oriented surface
$\Sigma\subset X\setminus\mathcal F$.
In this case there is a smooth orbifold one-form $a$ on $X$
such that $\omega=da$.

In particular, \eqref{eq:modelNull}
is equivalent to exactness of the K\"ahler form, and it is
enough to test its periods on surfaces avoiding any
additional finite set of smooth points.
\end{lemma}

\begin{proof}
For each $q\in\mathcal F$, choose pairwise disjoint small orbifold
coordinate neighborhoods $U_q$, with trivial group when $q$ is smooth.
Each $U_q$ is contractible, and its punctured neighborhood
retracts onto $S^3/\Gamma_q$.
Transfer for the finite covering
$S^3\to S^3/\Gamma_q$ gives
$
 H_1(S^3/\Gamma_q;\R)
 =
 H_2(S^3/\Gamma_q;\R)
 =
 0.
$
It follows that
\[
 H_2(U;\R)=0,
 \qquad
 H_1(U\setminus\mathcal F;\R)
 =
 H_2(U\setminus\mathcal F;\R)
 =
 0,\qquad U:=\bigcup_{q\in\mathcal F}U_q.
\]
The isomorphism follows from the Mayer--Vietoris sequence for
$
 X=(X\setminus\mathcal F)\cup U
$
in degree two:
\[
 0\longrightarrow
 H_2(X\setminus\mathcal F;\R)
 \longrightarrow
 H_2(X;\R)
 \longrightarrow0.
\]
Thus \eqref{eq:modelNull} implies that $[\omega]$ pairs trivially with every class in
$H_2(X;\R)$.
The orbifold de Rham theorem and the universal coefficient
theorem over $\R$ give \eqref{class vanishes}.
The remaining assertions follow immediately.
\end{proof}

We equip the model moduli with the pointed $\hat C^\infty$-CG
topology of Definition~\ref{def:modelhat}.

\begin{theorem}[Compactness of the model moduli]\label{thm:stationary}
Every sequence $
(X_j,h_j,J_j,p_j)\in\mathcal A(A,\Lambda,\alpha)
$ has a subsequence
converging in Definition~\ref{def:modelhat} to a member of the same
space. In particular, the underlying pointed metric family is
compact in pointed GH topology.
\begin{enumerate}[label=\textup{(\roman*)}]
\item
Every model and every limit has one ALE end, with every decay order
$\tau<2$. There are constants $\kappa=\kappa(A)>0$ and
$B=B(A,\Lambda)<\infty$ such that
\begin{equation}\label{eq:modelVolume}
 \kappa r^4\le\Vol_h B_h(x,r)\le Br^4,
 \qquad x\in X,\ r>0.
\end{equation}
\item
The natural volume measures converge weakly on bounded regions;
in particular, the volumes of corresponding bounded domains
converge when the limit boundary has zero volume.
\item
At a regular point there is neither curvature-energy
concentration nor accumulation of orbifold singular points.
\item \label{the convergence charts}
 For every $K\Subset U\Subset\mathcal R(X)$, the smooth embeddings
$\Phi_j$ in Definition \ref{def:modelhat} satisfy
the coverage property: if $x_j\in X_j$ converges
in the GH sense to a point of $K$, then $x_j\in\Phi_j(U)$
for all sufficiently large $j$.
See the beginning of Section~\ref{Convergence on the regular part}.
\end{enumerate}
\end{theorem}

The volume and preliminary compactness arguments use the critical
equations, the homogeneous Sobolev inequality, and the curvature
bound. The topology of K\"ahler ALE surfaces excludes branching.
Nullness is then used to exclude concentration over regular points.

%%%%%%%%%%%%%%%%%%%%%%%%%%%%%%%%%%%%%%%%%%%%%%%%%%%%%%%%%%%%%%%%%%%%%%%%%%%%%%%%%%%%%%%%%%%%%%%%%%%%%%%%%%%%%%%%%%%%%%%
\subsection{Uniform volume growth}
\label{Uniform volume growth}
Scalar-flat K\"ahler metrics satisfy
$\Delta\Ric=\Rm*\Ric$ and
$\Delta\Rm=\nabla^2\Ric+\Rm*\Rm$, where $*$ denotes universal
contractions. We use the epsilon-regularity theorem
\cite[Theorem~4.1]{TVmod}, also formulated under a lower volume
bound in \cite[Theorem~1.2]{TV}:
\begin{equation}\label{eq:modelEpsRegularity}
 \int_{B(x,2r)}|\Rm|^2<\eps_{\rm reg}
 \quad\Longrightarrow\quad
 \sup_{B(x,r)}|\nabla^k\Rm|
 \le C_k(A)r^{-2-k}
       \left(\int_{B(x,2r)}|\Rm|^2\right)^{1/2}.
\end{equation}
The estimate holds in smooth orbifold charts at isolated quotient
points; see \cite[proof of Proposition~4.3, pp.~903--904]{TV}.
For limiting punctures, \cite[Theorem~6.4]{TVmod} requires a connected
smooth punctured neighborhood, finite energy, a local Sobolev
inequality and an upper volume bound. We apply it branchwise only
after the critical annular construction of \cite[Section~7]{TVmod}
and \cite[Section~4]{TV} separates the branches.

The uniform upper bound \cite[Theorem~1.3]{TV} is stated for smooth
manifolds. We adapt it to orbifolds by counting the original quotient
points and the branches of provisional limits. The metric-dependent
orbifold bounds \cite[Theorems~1.1 and~2.5]{TV} enter only as
intermediate estimates.

\begin{proposition}\label{prop:stationaryVolume}
A complete scalar-flat K\"ahler ALE orbifold surface with the homogeneous Sobolev
\eqref{eq:modelSob} and $\int|\Rm|^2\le\Lambda$ satisfies the uniform volume growth
\eqref{eq:modelVolume}, with constants depending only on $A$ and
$\Lambda$.
\end{proposition}
\begin{proof}
\noindent\textit{Step 1: Sobolev constant controls orbifold order.}
The cutoff iteration in Lemma~\ref{lem:noncollapse}, applied to
\eqref{eq:modelSob}, gives
\begin{align}\label{the lower volume bound}
 \Vol B(x,r)\ge\kappa r^4,\qquad \kappa:=(16A)^{-2}.
\end{align}
At a quotient point the volume ratio tends to
$\omega_4/|\Gamma_x|$, hence
\begin{equation}\label{eq:modelGroup}
 |\Gamma_x|\le G(A):=\omega_4/\kappa.
\end{equation}
The same bound holds for every end group. Indeed, rescaling that end
passes the homogeneous Sobolev inequality with constant $A$ to its
punctured flat cone. Lemma~\ref{lem:sobolev-passage}, with $b=0$,
extends it across the vertex, where the preceding iteration applies.

\noindent\textit{Step 2: \texorpdfstring{$b_2$}{b2} and number of quotient singularities.}
The boundary links of a compact core $K$ are rational homology
three-spheres. The exact sequence gives
$H_c^2(X;\R)\simeq H^2(X;\R)$ with nondegenerate intersection
pairing. Thus compactly supported cohomology injects into reduced
$L^2$ cohomology: if $d\beta_j\to\zeta$ in $L^2$, pairing with any
compactly supported closed two-form shows $[\zeta]=0$.
The $L^2$ Hodge theorem and Carron's estimate \cite{Carron99}, in
the form following \cite[Theorem~2.1]{TV}, yield
\begin{equation}\label{eq:modelBetti}
 b_2(X)\le\dim\mathcal H^2_{(2)}(X)
       \le C(A)\int_X|\Rm|^2\,dV.
\end{equation}
The orbifold extension is explained in
\cite[proof of Theorem~2.5]{TV}. Orbifold Gauss--Bonnet on an ALE
exhaustion gives
\begin{equation}\label{eq:modelGB}
 \chi_{\rm top}(X)=\int_Xe(h)
 +\sum_{q\in\Sing X}(1-|\Gamma_q|^{-1})
 +\sum_{\rm ends}|\Gamma_\infty|^{-1}.
\end{equation}
Since $\chi_{\rm top}(X)\le1+b_2(X)$,
$|e(h)|\le C|\Rm|^2dV_h$, and every singular point contributes at
least $1/2$, we obtain
\begin{equation}\label{eq:modelOldCount}
 \#\Sing X\le2+C(A)\Lambda.
\end{equation}

\noindent\textit{Step 3: a bound for the number of branches.}
Suppose temporarily that a sequence of these orbifolds $(X_j,h_j,p_j)$ satisfies
\begin{align}\label{temp upper volume}
 \Vol B(x,r)\le C_0r^4,\qquad0<r\le r_0,
\end{align}
where $C_0,r_0>0$ are fixed for the sequence.
Together with \eqref{the lower volume bound}, this gives uniform local
packing bounds. Completeness and chains of intersecting balls give
covering bounds on every fixed larger ball. We therefore extract a
complete proper pointed GH limit $(X,d,p)$.
Let $N_{\mathrm{old}}:=\lceil2+C(A)\Lambda\rceil$ be the bound
in \eqref{eq:modelOldCount}, and fix
$\varepsilon_*:=\frac12\varepsilon_{\rm reg}(A)$.
Let $\mathcal E$ be the union of the limits of $\operatorname{Sing}X_j$
and the points with curvature-energy concentration at least
$\varepsilon_*$. Then
\[
 \#\mathcal E\le N_{\mathrm{old}}+\Lambda/\varepsilon_*.
\] Away from
$\mathcal E$, we have smooth convergence $h_i\to h_\infty$. The critical
annular construction under \eqref{temp upper volume} in \cite[Section 4]{TV} gives
finitely many punctured branches, and
\cite[Theorem~6.4]{TVmod} removes the puncture on each branch.

The smooth convergence away from $\mathcal E$ has coverage
from the critical annular construction. Together with the
temporary upper bound \eqref{temp upper volume},
Lemma~\ref{lem:volume-passage} gives convergence of the natural
volume measures on bounded regions.
In particular, volumes of corresponding bounded domains
converge whenever the limiting boundary has zero volume.
\begin{claim}The number of these branches is bounded independently of $C_0$.
\end{claim}
\begin{proof}
Let $q$ be a provisional multifold point with cone groups
$\Gamma_1,\ldots,\Gamma_m$. Take small spherical quotient links
in its smooth branches. Their source counterparts bound a connected
source neighbourhood $U_j$ containing the region that shrinks to \(q\). The boundary lies in the converging cone collars and avoids
all singular points of the underlying orbifolds. Orbifold Gauss--Bonnet for \(U_j\) gives
\[
 \chi_{\rm top}(U_j)
 =\int_{U_j}e(h_j)
  +\sum_{q'\in U_j\cap\Sing X_j}(1-|\Gamma_{q'}|^{-1})
  +\sum_{\ell=1}^m|\Gamma_\ell|^{-1}+o(1).
\]
Here the error tends to zero by first passing along the convergent
sequence and then shrinking the cone collars. Since $U_j$ is connected and has
nonempty boundary, $\chi_{\rm top}(U_j)\le1+b_2(U_j)$.
Mayer--Vietoris for $U_j$ and its exterior, whose common boundary
has zero real $H_1,H_2$, gives $b_2(U_j)\le b_2(X_j)$. By
$
1-\frac1{|\Gamma_{q'}|}\geq 0
$, we get
\begin{equation}\label{eq:modelConeCount}
 \sum_{\ell=1}^m|\Gamma_\ell|^{-1}
 \le1+b_2(U_j)+C\int_{U_j}|\Rm|^2+o(1)
 \le1+C(A)\Lambda+o(1).
\end{equation}
It remains to bound the groups $\Gamma_\ell$.
The temporary upper bound \eqref{temp upper volume} and smooth convergence away
from $\mathcal E$ allow us to apply Lemma~\ref{lem:sobolev-passage} with
$b=0$. Thus every completed local branch inherits the homogeneous
Sobolev inequality \eqref{eq:modelSob} with the original constant $A$.
The same is therefore true on its tangent cone. By Step~1,
$
 |\Gamma_\ell|\le G(A).
$
Combining this with \eqref{eq:modelConeCount}, we get
\[
 m\le G(A)\bigl(1+C(A)\Lambda\bigr).
\]
\end{proof}

\noindent\textit{Step 4: an upper volume bound for each limit.}
The preceding bounds allow us to apply
\cite[Proposition~4.3]{TV} to every provisional smooth multifold
limit. There are constants $\eps>0,A_0>\omega_4$, depending only on
$A,\Lambda$, such that
\begin{equation}\label{eq:modelSmallEnergyVolume}
 \int_{B(x,2)}|\Rm|^2<\eps
 \quad\Longrightarrow\quad
 \Vol B(x,1)\le A_0.
\end{equation}
Here smooth orbifold removal in Step~3 precedes the use of
epsilon-regularity and Bishop comparison on each branch; see also
\cite[Remark on p.~907]{TV}. For a multifold one counts the branches
through the exceptional points. Their number is bounded in terms of
$A,\Lambda$, since $\varepsilon_*$ was fixed in terms of $A$ alone.
Thus neither $\eps$ nor $A_0$ depends on the temporary constant $C_0$.

\begin{claim}\label{ALE ends}
Every such complete provisional limit $X$ also has finite
volume growth upper bound
\begin{align}\label{A1 upper bound}
\Vol_{h_\infty}B(p,r)\le A_1r^4
\quad\text{for all }r>0,
\end{align}
where \(A_1\) may depend on the particular complete provisional limit.
\end{claim}
\begin{proof}
By the branchwise conclusion of
Lemma~\ref{lem:sobolev-passage}, every normalized component
$\widetilde X_a$ inherits
\[
 \|u\|_{L^4(\widetilde X_a)}^2
 \le
 A\int_{\widetilde X_a}|du|^2\,dV.
\]
In particular, $\widetilde X_a$ cannot be compact.
Finite curvature energy and
\eqref{eq:modelEpsRegularity} on exterior balls give
$|\Rm|=o(r^{-2})$. Also $\int|\Ric_-|^2<\infty$.
Theorems~1.1 and~2.5 of \cite{TV} therefore give finitely many ALE
ends of order zero and finite upper-volume growth on each
$\widetilde X_a$. Gluing $\widetilde X_a$ at
finitely many points preserves such a bound, with a constant
depending on that limit. This proves \eqref{A1 upper bound}.
\end{proof}

\noindent\textit{Step 5: a uniform upper volume bound.}
Suppose the required uniform upper bound fails. After selecting a
sequence we have
\[
 \sup_{x\in X_j,\ r>0}r^{-4}\Vol B(x,r)\longrightarrow\infty.
\]
For a fixed threshold $H>A_0$ define the first violating radius
\begin{align}\label{aj threshold}
 a_j(H):=\inf\{r>0:\ \Vol B(x,r)>Hr^4
                         \text{ for some }x\in X_j\}.
\end{align}
For large $j$ this is finite and positive: each fixed smooth ALE
orbifold has bounded curvature and uniform small-radius upper ratios.
Choose $a_j(H)\le r_j<2a_j(H)$ with
$\Vol B(x_j,r_j)>Hr_j^4$, and put $\tilde h_j:=r_j^{-2}h_j$.
Then $\Vol_{\tilde h_j}B(x_j,1)>H$ and, by the definition of $a_j(H)$,
$\Vol_{\tilde h_j}B(y,\rho)\le H\rho^4$ for every $y$ and
$0<\rho\le1/2$. Step~3 therefore gives a complete pointed multifold
limit with $\Vol B(x_\infty,1)\ge H>A_0$. By
\eqref{eq:modelSmallEnergyVolume},
\[
 \int_{B(x_\infty,2)}|\Rm_{h_\infty}|^2\,dV_{h_\infty}\ge\eps.
\]

Let \(\mathcal E_\infty\) denote the finite exceptional set of convergence.
We choose $
K\Subset B(x_\infty,2)\setminus \mathcal E_\infty$ such that
$
\int_K|\operatorname{Rm}_{h_\infty}|^2\,dV_{h_\infty}
>\frac{3\varepsilon}{4}
$. Then smooth convergence provides \(K_j\subset X_j\) satisfying
\[
\int_{K_j}|\operatorname{Rm}_{h_j}|^2\,dV_{h_j}
=
\int_{K_j}|\operatorname{Rm}_{\tilde h_j}|^2\,dV_{\tilde  h_j}
\geq\frac{\varepsilon}{2}.
\]
 We repeat this construction,
increasing the fixed threshold on successive subsequences.
Suppose finitely many cores $
K_j^{(l)},1\leq l\leq m$ have already been chosen. Write
$B(x_j^{(\ell)},r_j^{(\ell)})$ for their selecting balls. Let
$A_\ell<\infty$ be the constant in \eqref{A1 upper bound}, as supplied by Step~4. Choose the next threshold larger
than $10000\max\{A_0,A_1,\ldots,A_{m}\}$ and all preceding
thresholds. Because $a_j(H)$ is nondecreasing in $H$, its selected
radius $r_j$ satisfies $r_j>r_j^{(\ell)}/2$ for every earlier
choice.

The distinction between bounded and unbounded ratios of selecting
radii is essential. The latter case was omitted in the earlier
Tian--Viaclovsky argument; see their correction
\cite[Section~4, pp.~906--907]{TV} and Chen--Weber's treatment
\cite[proof of Theorem~4.3]{CW}.
An earlier pointed limit controls volume only at bounded radii
in its own scale. We count disjoint curvature cores, not selecting
balls, which may be nested.

After passing to a subsequence, each ratio
$\lambda_j:=r_j/r_j^{(\ell)}$ tends either to a finite
$\lambda\ge\frac12$ or to infinity.
In the finite case the two radius-two selecting balls must be
disjoint for large $j$. Otherwise
$
 B(x_j,r_j)\subset
 B(x_j^{(\ell)},8r_j).
$
In the earlier rescaling, the radius $8\lambda_j$ tends to
$8\lambda$. Distance spheres in a smooth multifold have zero
natural volume, branchwise as for a smooth metric. Sandwiching
these balls between fixed radii $8\lambda\pm\delta$ and then
letting $\delta\downarrow0$ justifies volume convergence at these
moving radii. Thus \eqref{A1 upper bound} gives, for large $j$,
\[
Hr_j^4
<
\Vol_{h_j}B(x_j,r_j)\leq
\Vol_{h_j}B(x_j^{(\ell)},8r_j)
\leq 2A_\ell(8r_j)^4,
\]
contradicting $H>10000A_\ell$.

If $r_j/r_j^{(\ell)}\to\infty$, an earlier compact core has
diameter tending to zero in the new scale. It either leaves every
bounded pointed ball or converges to a single point. There are
only finitely many earlier cores. In the new limit, choose the
next compact regular core away from all of their limit points
and from the exceptional set. Removing these finitely many points
does not change the integral of the smooth limiting curvature
density, so the core can still have energy exceeding $3\eps/4$.
This uses the curvature measure of the smooth limit, not the weak
limit of the approximating curvature measures, which may have atoms.
Its source representatives are disjoint from every earlier core.

Thus each stage provides a new disjoint source region with energy
at least $\eps/2$. More than $2\Lambda/\eps$ stages contradict
the total curvature energy bound. Under the supposed unboundedness
of maximal volume ratios, each of the finitely many successive
fixed thresholds could have been selected. This contradiction
proves the uniform upper bound.
\end{proof}

\subsection{Preliminary compactness and one end}
\label{Preliminary compactness and the number of ends}
We apply the preceding argument to obtain pointed compactness of the
models. We record each original quotient point, including those where
the curvature vanishes.

\begin{proposition}[Preliminary pointed compactness]
\label{prop:modelprecompact}
Every sequence
$
(X_j,h_j,J_j,p_j)\in\mathcal A(A,\Lambda,\alpha)
$
has a subsequence converging in the pointed GH sense to a scalar-flat K\"ahler multifold \((X,d,p)\), which is a connected complete proper metric space.
There is a finite exceptional set $\mathcal E$ such that
\[
\operatorname{Sing}X\subseteq\mathcal E\subset X,
\qquad \#\mathcal E\le\Lambda/\alpha,
\]
and the metrics and complex structures converge smoothly on compact subsets of \(X\setminus\mathcal E\).

Every component \((\widetilde X_a,h_a,J_a)\) is a complete noncompact scalar-flat K\"ahler ALE orbifold with finitely many ends. Each end has differentiable decay of every order \(0<\tau<2\).
\end{proposition}
\begin{proof}
Step~3 of Proposition~\ref{prop:stationaryVolume}, now with its uniform
upper volume bound, gives the preliminary extraction. To record old
quotient points, including flat ones, use the marked measures
\begin{equation}\label{eq:modelMarkedMeasure}
 \mu_j:=|\Rm_{h_j}|^2dV_{h_j}
       +\alpha\sum_{q\in\Sing X_j}\delta_q.
\end{equation}
A diagonal subsequence converges weakly on bounded regions to a
Radon measure $\mu$ with $\mu(X)\le\Lambda$. Set
$\mathcal E:=\{x:\mu(\{x\})\ge\alpha\}$, so
$\#\mathcal E\le\Lambda/\alpha$.

At $x\notin\mathcal E$, choose a neighborhood of marked measure
less than $\alpha<\eps_{\rm reg}$. The corresponding source
neighborhoods contain no orbifold points and have curvature energy
below $\eps_{\rm reg}$. The estimate \eqref{eq:modelEpsRegularity},
noncollapse and injectivity-radius control give smooth convergence
on smaller balls. Lemma~\ref{lem:parallel-complex} gives simultaneous
smooth convergence of $J_j$. The critical annular construction and
\cite[Theorem~6.4]{TVmod} complete every punctured branch to a smooth
orbifold; Lemma~\ref{lem:parallel-complex} extends its parallel complex
structure equivariantly. Multiple branches are still allowed.

The extraction includes coverage. Together with the uniform upper
volume bound, Lemma~\ref{lem:volume-passage} therefore gives natural
volume convergence on bounded regions, and convergence of volumes
for domains with zero-volume limiting boundary.

Claim~\ref{ALE ends} shows that each normalized component is complete,
noncompact and ALE of order zero, with finite topological type and
homogeneous Sobolev constant $A$. The decay theorems
\cite[Theorem~6.1]{TVmod} and \cite[Theorem~5.1]{TV} improve this to
$|\partial^k(h-e)|=O(r^{-\tau-k})$ for every $0<\tau<2$.
The parallel equation for $J$ gives
$|\partial^k(J-J_0)|=O(r^{-\tau-k})$ for a constant orthogonal
complex structure $J_0$. The end group preserves $J_0$ and is
therefore conjugate into $U(2)$.
\end{proof}

Hein--LeBrun prove one-endedness for smooth ALE K\"ahler surfaces
\cite[Proposition~4.2]{HL}. The positive-order decay needed here is
also discussed in \cite[Section~2]{HRS}. We give the harmonic-function
argument, including the integrations at isolated quotient points.

\begin{lemma}[One end]\label{lem:oneend}
Let $X$ be a connected complete K\"ahler ALE orbifold surface, with finitely
many isolated quotient points. Suppose each end has positive $\tau$. Then $X$ has exactly one end.
\end{lemma}
\begin{proof}
Suppose there are at least two ends $E_a$. Choose end values
$c_a\in\{0,1\}$, not all equal. Harmonic Dirichlet solutions on a
compact exhaustion, with these boundary values, satisfy $0\le u_R\le1$
and have uniformly bounded Dirichlet energy by comparison with one
fixed function. A subsequence gives a finite-energy harmonic function
$u$ on $X$. For $1<p<2$, the barrier
$\Delta_h r^{-p}=[p(p-2)+O(r^{-\tau})]r^{-p-2}<0$ gives
\begin{align}\label{uca}
 u-c_a=O(r^{-p}),\qquad |du|=O(r^{-p-1}) ,\quad
\text{ on }E_a,
\end{align}
so $u$ is nonconstant.

Put $\beta:=-du\circ J=2d^cu$. On a K\"ahler surface,
$d\beta$ is primitive of type $(1,1)$, hence anti-self-dual.
Integration by parts with exhaustion cutoffs gives
$\int\chi_R^2|d\beta|^2\le CR^{-2}\int|du|^2\to0$.
Thus $d\beta=0$. All integrations hold in orbifold charts; excising
small balls about quotient points introduces boundary terms tending
to zero by smoothness of the lifts.

Since $H^1(S^3/\Gamma_a;\R)=0$, write $\beta=dv_a$ on each end.
The K\"ahler identity $*du=\beta\wedge\omega$ yields
\begin{align}\label{the individual zero-flux identities}
 \int_{S_{a,R}}\partial_\nu u\,dA
 =\int_{S_{a,R}}dv_a\wedge\omega
 =\int_{S_{a,R}}d(v_a\omega)=0.
\end{align}
Green's identity and these individual zero fluxes give
\[
 \int_{X_R}|du|^2
 =\sum_a\int_{S_{a,R}}(u-c_a)\partial_\nu u\,dA
 =O(R^{2-2p})\longrightarrow0,
\]
where \eqref{uca} and $\operatorname{Area}(S_{a,R})=O(R^3)$ were used.
This contradicts the nonconstancy of $u$.
\end{proof}

\subsection{One branch and bubble trees}

We first exclude multiple local branches in the preliminary limits.
We then select finitely many nonflat blow-up scales, retain the original
singular points as marked points, and organize these data into bubble trees.
The construction follows the bubble arguments of
\cite{Bando,BandoCorrection,TVmod,TV}; the topological input is the
K\"ahler resolution construction in \cite[proof of Lemma~2.11]{HV}
and the vanishing theorem \cite[Proposition~3.6]{HVdef}.

\begin{lemma}[Finite bubble trees]\label{lem:modelhierarchy}
Let \((X,d,p)\) be the preliminary limit and \(\mathcal E\) its exceptional set from Proposition~\ref{prop:modelprecompact}. Then every point of \(X\) has exactly one local orbifold branch; hence \(X\) is a scalar-flat K\"ahler orbifold.

After passing to a further subsequence, the degeneration over \(\mathcal E\) is described by finitely many bubble trees of scalar-flat K\"ahler ALE orbifolds \((Y,h,J)\), possibly flat. The following properties hold.
\begin{enumerate}[label=\textup{(\roman*)}]
\item Every source quotient point converging to $\mathcal E$ is retained
in an orbifold chart on a component of these trees or in a terminal
quotient chart.
\item For each parent--child pair, the corresponding source regions are
joined by a smooth annulus diffeomorphic to
$(S^3/\Gamma)\times(0,1)$. The group at infinity of the child and the
local group at the attachment point of the parent agree up to unitary
conjugacy of their linear actions.
\item Each connecting annulus admits compatible cone charts. For fixed
$0<a<b$, consider subannuli $\{a\rho<r<b\rho\}$ whose larger
subannuli $\{a\rho/2<r<2b\rho\}$ lie between the two end collars;
here $r$ is the radial coordinate. On the smaller subannuli, the
rescaled metrics and complex structures converge smoothly to those of
the corresponding flat quotient, uniformly over all such scales $\rho$.
\end{enumerate}
The number of components and terminal quotient charts is bounded by
a constant depending only on $A$, $\Lambda$, and $\alpha$.
\end{lemma}
\begin{proof}
We use the source sequence $(X_j,h_j,J_j,p_j)$ from
Proposition~\ref{prop:modelprecompact} and write $d_j:=d_{h_j}$ for its
length distances. The constants $A,\Lambda,\alpha$, the volume bounds
\eqref{eq:modelVolume}, and the regularity estimate
\eqref{eq:modelEpsRegularity} remain fixed. A blow-up $$
(X_j,r_j^{-2}h_j,x_j)\longrightarrow(Y,h)$$ is described by a center and scale
$
(x_j,r_j).
$
Two choices are equivalent if their scale ratio tends to a positive
finite number and their center distance is bounded relative to that
scale. Their limits then differ only by a constant rescaling and a
change of base point.
Parent and child refer to consecutive larger and smaller blow-up
scales.

We prove the lemma in the following steps.

\noindent\textit{Step 1: vanishing \texorpdfstring{$H^1$}{H1}.}
By Lemma~\ref{lem:oneend}, every $X_j$ has one end.
\begin{claim}\label{eq:modelb1zero}
$
 H^1(X_j;\R)=0.
$
\end{claim}
\begin{proof}
Resolve the quotient singularities to a smooth complex surface
$\widetilde X_j$. The construction in
\cite[proof of Lemma~2.11]{HV} gives a K\"ahler metric agreeing with
$h_j$ outside the quotient neighborhoods. Each fixed source already
has finite ALE topology, so the decay theorem used in
Section~\ref{Preliminary compactness and the number of ends} applies
before any uniform Betti bound. The auxiliary metric has ALE order
greater than one, and \cite[Proposition~3.6]{HVdef} gives
$H^1(\widetilde X_j;\R)=0$.

Replacing a cone on $S^3/\Gamma$ by its resolution does not change
first real cohomology. Indeed, the resolved filling retracts onto a
tree of rational curves; both fillings and their common boundary have
zero $H^1$, and the $H^0$ restriction map is surjective. Mayer--Vietoris
therefore gives $H^1(X_j;\R)\simeq H^1(\widetilde X_j;\R)=0$.
\end{proof}

\noindent\textit{Step 2: the provisional limit is orbifold.}
Suppose that a provisional limit has at least two local branches at
$q$. Choose a small regular ball in each of two branches, away from
$q$ and inside a neighborhood where the branches are disjoint. Fix
these balls once and for all. Their source images have volumes bounded below by a fixed number
$
v_0>0.
$

For a sufficiently small $r$, take a spherical quotient cross-section
of radius $r$ in the first branch. Smooth convergence on the radial
collar between radii $r/2$ and $2r$ produces a connected cooriented hypersurface
\[
 L_{j,r}\simeq S^3/\Gamma\subset\mathcal R(X_j).
\]
Because \(H^1(X_j;\mathbb R)=0\), this hypersurface separates \(X_j\). Indeed, a nonseparating cooriented hypersurface
would admit a closed curve with algebraic intersection number one,
contradicting Claim \ref{eq:modelb1zero}. Since $X_j$ is proper and one-ended, exactly one
of the two complementary sides is relatively compact.

\begin{claim}
The two fixed balls lie on opposite sides, even if the limiting
branches reconnect outside the chosen neighbourhood.
\end{claim}
\begin{proof}
we choose an inner point at distance approximately $r/4$ from $q$ in
each branch. Their source representatives have distance at most
$r/2+o(1)$. A minimizing path between them stays within distance
$3r/4+o(1)$ of a source representative of $q$, and therefore misses
$L_{j,r}$, whose distance from that representative is $r+o(1)$.
Join each inner point to its fixed ball within its limiting branch.
The first path crosses the selected cross-section once, whereas the
second misses it. Transferring these paths to the sources and
concatenating them with the inner minimizing path gives a path
between the fixed balls with exactly one transverse crossing of
$L_{j,r}$. Since $L_{j,r}$ separates, the balls are on opposite sides.
\end{proof}

Choose a compactly supported function $u_{j,r}$ equal to one on the
relatively compact side and zero on the other, with
$|du_{j,r}|\le C/r$ in the collar. The collar has volume $O(r^4)$,
whereas one fixed ball lies in $\{u_{j,r}=1\}$. Thus

\begin{equation}\label{eq:modelnobranchcutoff}
 v_0^{1/2}\le \|u_{j,r}\|_{L^4(h_j)}^2
 \le A\int_{X_j}|du_{j,r}|_{h_j}^2\,dV_{h_j}
 \le CAr^2.
\end{equation}
For each fixed $r$, take $j$ sufficiently large, and then let
$r\downarrow0$. The contradiction excludes multiple branches.
This is the separating-collar argument of
\cite[Proposition~7.2]{TVmod}, with one-endedness providing the
compact side.

The same argument applies to every complete rescaled limit used later. Indeed,
Claim \ref{eq:modelb1zero} and the number of ends are unchanged by rescaling,
and the homogeneous Sobolev constant and four-dimensional curvature
energy are scale invariant. Every provisional limit used below is
therefore an orbifold. Proposition~\ref{prop:modelprecompact} and
Lemma~\ref{lem:oneend} make its completion a one-ended scalar-flat
K\"ahler ALE orbifold.

\noindent\textit{Step 3: bubbles are nonflat.}
We first record local orbifold compactness under a curvature bound,
including old source quotient points. If $|\Rm|\le K$ on $B(x,2r)$,
choose $2a<r$ with $K^2B(2a)^4<\eps_{\rm reg}$. The upper volume
bound and \eqref{eq:modelEpsRegularity} give all curvature derivative
bounds on smaller balls. Noncollapse then permits the frame-bundle
construction of \cite[Proposition~4.1]{KLorb}: smooth
$O(4)$-equivariant comparison maps descend to orbifold maps covering
the smaller source regions and identifying their isotropy actions.
Consequently, isolated quotient points cannot coalesce, appear or
disappear there while curvature and noncollapse remain controlled.

If curvature is unbounded on smaller balls $B_{h_j}(x_j,r)$,
take $B_{h_j}(x_j,2r)$ as the outer region and select points
$z_j$ by maximising $$
|\operatorname{Rm}|(z)\,
d(z,\partial\text{outer region})^2
.$$ With $Q_j:=|\Rm_{h_j}|(z_j)$, the point-selection argument gives, after rescaling,
\[
\tilde h_j:=Q_jh_j,
\qquad
|\operatorname{Rm}_{\tilde h_j}|(z_j)=1,
\]
with uniformly bounded curvature on balls whose rescaled radii tend to infinity.
The preceding local argument gives smooth orbifold convergence on
each fixed ball. The resulting limit is complete and nonflat. It
inherits the Sobolev and volume bounds and finite curvature energy;
by Section~\ref{Preliminary compactness and the number of ends} and Step~2, it is a one-ended scalar-flat K\"ahler
ALE orbifold. Thus nonflat blow-ups can be extracted before the
bubble trees have been constructed.

\noindent\textit{Step 4: inequivalent nonflat scales are finite.}
Every nonflat complete bubble has energy at least
$\eps_{\rm gap}(A)>0$. Otherwise \eqref{eq:modelEpsRegularity} on
arbitrarily large balls gives
$|\Rm|(x)\le C(A)R^{-2}\eps_{\rm gap}^{1/2}\to0$, contradicting
nonflatness.

The smooth curvature density of a completed bubble has no atoms.
For any finite collection of inequivalent nonflat limits, choose
compact regular cores $K$ avoiding the exceptional points and the
limiting centers of finer bubbles, with energy greater than
$3\eps_{\rm gap}/4$. Smooth convergence and scale invariance give
\begin{align}\label{local j energy}
 \int_{K_j}|\Rm_{h_j}|^2\,dV_{h_j}
 \ge\tfrac12\eps_{\rm gap},\qquad K_j:=\Phi_j(K), \qquad \text{ for large }j.
\end{align}

Choose any finite collection of complete nonflat blow-ups, with \((x_j^{(a)},r_j^{(a)})\), \(a=1,\ldots,N\), which are pairwise inequivalent.
The selected limiting cores $
K_j^{(a)}$ can then be made disjoint. At comparable scales, inequivalent centers move infinitely far apart in normalised distance. At separated scales, the finer core lies in a small neighbourhood
omitted from the coarser core. Summing \eqref{local j energy} over these disjoint cores yields
\begin{equation}\label{eq:modelessentialcount}
 N\le\frac{2\Lambda}{\eps_{\rm gap}}.
\end{equation}

Hence the attainable cardinalities have a largest value, denoted again
by $N$. Choose a family attaining $N$, retain its common
subsequence, and call its members the selected nonflat labels.

\noindent\textit{Step 5: construct the bubble trees.}
We also assign a label, called a \emph{mark}, to each original
singular point converging to $\mathcal E$. Their number satisfies $N_0\le\Lambda/\alpha$. Write
$(x_j^{(a)},r_j^{(a)})$ for the center and length scale of a label.
For a label representing an original singular point, set $r_j^{(a)}:=0$.

For each point of $\mathcal E$, we start with the original limit as
parent with scale $s_j=1$. At each later stage, a cluster
consists of the remaining labels whose centers converge to the
same parent attachment; their positive scales are $o(s_j)$.
A singleton nonflat label gives its selected component, while a
singleton label of an original singular point terminates in a quotient chart.

For a cluster with at least two labels, set
\[
 r_j:=\max\left\{
       \max_a r_j^{(a)},\;
       \max_{a,b}d_j(x_j^{(a)},x_j^{(b)})
       \right\},
\]
where the maxima range over that cluster. Choose one label center
as the base point. Then $r_j/s_j\to0$, and all cluster centers
remain at bounded distance in the $r_j$-rescaling. By
Proposition~\ref{prop:modelprecompact} and Step~2, a subsequence
converges to a complete ALE orbifold limit. Pass to a further
subsequence so that all normalised label scales and centers
converge.

If some $r_j^{(a)}/r_j$ has a positive limit, the new limit is the
selected nonflat component of that label, up to rescaling and
recentering. There is at most one such label, since two would be
equivalent. Retain this component, remove its label from the list,
and group the remaining labels by their limiting attachment points.
Their scales are now $o(r_j)$.

Otherwise, every positive label scale is $o(r_j)$. The maximum is
then determined by center distances, so at least two centers have
distinct limits. The complete limit is flat: a nonflat limit would,
by Step~4, be equivalent to a selected label of comparable scale
in this cluster, a contradiction. Retain this flat component and
group the labels by their distinct limiting attachment points.
Repeat the construction on each resulting cluster.

Every unresolved attachment is represented by a remaining label.
If its normalized curvature is unbounded, Step~3 gives a nonflat
blow-up and Step~4 supplies its label. If curvature is bounded,
the comparison maps of Step~3 identify each singular point with an
original quotient point.

Every nonflat level removes one label. Every retained flat level
splits a cluster into at least two nonempty smaller clusters.
Hence the procedure terminates. There are $N$ nonflat labels and
$N_0$ terminal quotient marks, and the rooted-forest count gives
\[
 \#\{\text{retained flat nodes}\}
 \le\max\{0,N+N_0-1\}.
\]
Indeed, each leaf accounts for a distinct label and every retained
flat node has at least two children. Together with
\eqref{eq:modelessentialcount} and $N_0\le\Lambda/\alpha$, this
bounds the number of components and terminal charts in terms of
$A,\Lambda,\alpha$.

Every complete flat limit is a finite Euclidean quotient.
Its orbifold universal cover is Euclidean. For a regular base point
$o$ and a lift $\widetilde o$, let $D$ be a Dirichlet region and
$D_R:=D\cap B_e(\widetilde o,R)$. Then
\begin{align}\label{vole volt}
 \Vol_e(D_R)=\Vol_hB_h(o,R).
\end{align}
For any $m$ distinct deck transformations, their translates of $D_R$
are disjoint up to null sets and lie in $B_e(\widetilde o,R+L_m)$,
where $L_m$ is their maximal displacement of $\widetilde o$. Hence
$m\kappa R^4\le\omega_4(R+L_m)^4$. Letting $R\to\infty$ bounds
$m$, so the deck group is finite. A common fixed point identifies
its action as linear; the parallel complex structure makes it unitary.
Since singularities are isolated, any nontrivial group acts freely
on $S^3$, and the vertex is the unique singular point.

A flat level with just one marked attachment is absorbed into a
connecting edge. For a nontrivial quotient, that attachment must
be the vertex; otherwise the vertex would require a second label.
For a Euclidean level, take the attachment as origin. Thus retained
flat nodes separate labels at distinct points.

For each terminal old point, choose an actual positive scale
$\sigma_j$ with $\sigma_j/s_j\to0$, where $s_j$ is its immediate
parent scale. Choose $\sigma_j$ sufficiently small relative to the
individual uniformizing chart so that, in normal coordinates on
the smooth cover, $\sigma_j^{-2}h_j$ and $J_j$ converge on expanding
balls to the tangent flat quotient. Smoothness of each lifted
metric allows this diagonal choice. The terminal chart therefore
supplies an inner cone collar at a positive scale.

\noindent\textit{Step 6: no further degeneration between consecutive nodes.}
Let $c_j$ be the center of the child cluster. Let the child have
scale \(r_j\) and its parent scale \(s_j\), with
$
r_j/s_j\to0.
$
For an original singular point at a terminal node, use its positive scale $\sigma_j$ for $r_j$. The child's ALE end describes the inner side of the connecting region. A small quotient neighborhood of the parent attachment describes the outer side. Choose the two end cross-sections diagonally, retaining convergence
on collars from half to twice their physical radii, namely
\[
 R_jr_j,\qquad\delta_js_j,\qquad
 R_j\to\infty,\qquad\delta_j\to0,\qquad
 \frac{R_jr_j}{\delta_js_j}\to0.
\]
Thus the child is viewed far out toward infinity, while the parent is viewed very close to its attachment point.

Let $r_j\ll\rho_j\ll s_j$. The complete limit at scale \(\rho_j\) must be flat. Otherwise, a nonflat limit would be equivalent to a selected
label with comparable scale and bounded normalized center
distance. This contradicts the fact that the parent and child are
consecutive.

There is no further curvature concentration on a fixed annulus away
from the central point. For the local estimates, retain the concentric
annulus with half its inner radius and twice its outer radius.
If curvature concentration occurred there, point selection
would give a nonflat limit at a scale smaller than $\rho_j$, with
center at distance comparable to $\rho_j$ from $c_j$. Its label
would lie neither in the child cluster nor at another parent
attachment, contradicting the construction in Step~5. The same
annulus contains no old source quotient mark. Consequently its
normalized curvature is bounded, and Step~3 gives smooth
convergence with coverage of every point of the underlying orbifold on smaller
annuli.

If the complete flat limit has a nontrivial group, its vertex must
be the central point. A vertex elsewhere would lie in one of
these annuli and would require an original singular point or a curvature
concentration there. For a Euclidean limit, take the central
point as origin. Hence every fixed annulus in every intermediate
limit is a centered flat quotient annulus.

\noindent\textit{Step 7: uniform cone charts.}
The preceding statement must hold uniformly over the intervening scales.
If it failed, we could choose offending scales \(\rho_j\), a fixed annular radius ratio, and a derivative order for which the cone approximation fails. By extraction at those scales, we get exactly the flat, smoothly converging annulus established in Step 6, which is a contradiction.

The linear actions agree up to unitary conjugacy:
\[
\Gamma_{\infty,\mathrm{child}}
=
U\Gamma_{\mathrm{parent},q}U^{-1}
\quad\text{for some }U\in U(2).
\]
On overlapping cone annuli, the limiting comparison maps preserve radius and identify the spherical space forms. Their lifts conjugate the orthogonal group actions. Because the complex structures also converge, these lifted identifications are complex linear, hence unitary.

The group-order bound from Section~\ref{Uniform volume growth} leaves only finitely many possible unitary quotient types. Consequently, sufficiently close overlapping charts cannot change type. Passing through the overlaps identifies the child's end action with the parent's local action.

\noindent\textit{Step 8: connecting region is a product.}
Choose overlapping dyadic annuli. In each chart, use the logarithm of the physical radial coordinate and the radial vector field
$
r\partial_r.
$
On a cone,
$
(r\partial_r)(\log r)=1.
$
Write $t_a$ and $V_a$ for these local logarithmic radii and radial
fields. After aligning adjacent charts, on every overlap
$|t_a-t_b|+|V_b(t_a)-1|\le C\epsilon_j$, with $\epsilon_j\to0$
uniformly in the annulus index. Choose a partition of unity
$\{\chi_a\}$ with bounded overlap and uniformly bounded
$V_b(\chi_a)$, and put $T:=\sum_a\chi_at_a$ and
$V:=\sum_a\chi_aV_a$. Since $\sum_aV(\chi_a)=0$, at each point
choose an index $b$ whose chart contains that point and write
\[
 V(T)=\sum_{a,c}\chi_a\chi_cV_c(t_a)
       +\sum_aV(\chi_a)(t_a-t_b)=1+O(\epsilon_j)>0.
\]
The constant is independent of the number of annuli.
Thus \(T\) has no critical points in the connecting region.
Move the two boundary collars slightly to level sets of \(T\).
The region between them is compact, and along the flow of
$V/V(T)$ the value of $T$ increases at unit speed. The flow
therefore identifies all intermediate levels. Each level is diffeomorphic to \(S^3/\Gamma\), so the region is diffeomorphic to
\[
(S^3/\Gamma)\times[a_j,b_j].
\]
Coverage of all points of the underlying orbifolds rules out an additional component hidden between the two boundary levels. The bounded overlap of the dyadic cover keeps these local estimates independent of how many annuli are needed.
This is the annular patching used in the critical bubble
construction; compare \cite[Section~7]{TVmod} and
\cite[Section~4, pp.~901--903]{TV}. Together with the finite
clustering of all nonflat labels and old quotient marks, it
accounts for the degeneration over $\mathcal E$ and proves the
lemma.
\end{proof}

\subsection{Nullness, rigidity and regular points of the limit}

\begin{lemma}[Nullness]\label{lem:modelnullclosed}
The preliminary limit and every ALE orbifold component of the bubble
trees in Lemma~\ref{lem:modelhierarchy} have zero real K\"ahler class.
For each such component $(Y,h,J)$, with K\"ahler form $\omega_h$,
this means $[\omega_h]=0\in H^2_{\mathrm{dR}}(Y;\mathbb R)$.
\end{lemma}
\begin{proof}
Let $\mathcal F\subset Y$ be
its finite set excluded from smooth convergence. On compact subsets of $\mathcal R(Y)\setminus \mathcal F$, the convergence
is described by embeddings $\Phi_j$ for which
\[
 \Phi_j^*(r_j^{-2}\omega_j)\longrightarrow\omega_\infty
 \qquad\text{smoothly},
\]
where $r_j=1$ for the unrescaled limit. For every compact smooth
closed oriented surface $\Sigma\subset\mathcal R(Y)\setminus \mathcal F$,
\[
 \int_\Sigma\omega_\infty
 =\lim_{j\to\infty}
   \int_\Sigma\Phi_j^*(r_j^{-2}\omega_j)
 =\lim_{j\to\infty}
   r_j^{-2}\int_{\Phi_j(\Sigma)}\omega_j
 =0.
\]
Applying Lemma~\ref{nullness condition} to the finite set
$\operatorname{Sing}Y\cup\mathcal F$, we conclude that
$
 [\omega_\infty]=0.
$
The same argument applies to the completed main limit and
every component of the bubble trees.
The terminal flat quotient charts are null because they
are contractible.
\end{proof}

The K\"ahler mass computation of Hein--LeBrun \cite{HL} is local
on a smooth ALE end. We combine it with the orbifold positive-mass
rigidity theorem \cite[Theorem~1.9]{JV}. The following proof
makes explicit how the vanishing K\"ahler class removes the mass
term even when the compact interior contains quotient points.

\begin{lemma}[Null AE rigidity]\label{lem:nullAE}
Let $(X,h,J)$ be a complete scalar-flat K\"ahler ALE
orbifold surface. Suppose that $X$ is null and has an
AE end with differentiable decay order $\tau>1$.
Then $(X,h,J)$ is biholomorphically isometric to
$(\mathbb C^2,e,J_0)$. In particular,
$\operatorname{Sing}X=\emptyset$.
\end{lemma}
\begin{proof}
By Lemma~\ref{lem:oneend}, the given AE end is the only end.
By nullness and orbifold de Rham theory, $$
 \omega=da.
$$ for some smooth orbifold one-form $a$ on $X$.
The Ricci form $\rho$ is a smooth closed orbifold
two-form. 
As in the proof of Lemma~\ref{nullness condition},
the end satisfies $H^2(E;\R)=0$.
Hence there is a smooth one-form $\theta$ on $E$
such that $$d\theta=\rho|_E.$$
Choose an end
cutoff $\chi$ which vanishes near the inner boundary of $E$
and is identically one sufficiently far out. Extend $\chi\theta$
by zero and put
$
 \rho_c:=\rho-d(\chi\theta).
$
Then $\rho_c$ is closed, smooth as an orbifold form, and
compactly supported. Therefore
\[
 \int_X\rho_c\wedge\omega
 =\int_X\rho_c\wedge da
 =\int_Xd(\rho_c\wedge a)=0.
\]

In complex dimension two,
$\rho\wedge\omega=S\,dV_h=0$ because the metric is scalar-flat.
Choose $R_2$ so that $\chi=1$ on the AE end wherever $r\ge R_2$.
For $R>R_2$, let $S_R:=\{r=R\}$ and let $X_R$ be the compact
orbifold obtained by removing $r>R$.
Then $\partial X_R=S_R$, $\operatorname{supp}\rho_c\subset X_R$,
and $\chi=1$ near $S_R$. The last
identity consequently gives
\[
 \begin{aligned}
 0
 =\int_{X_R}\rho_c\wedge\omega
 =\int_{X_R}\rho\wedge\omega
    -\int_{X_R}d(\chi\theta\wedge\omega)
 =-\int_{S_R}\theta\wedge\omega.
 \end{aligned}
\]
At an interior quotient point the forms lift smoothly, so
orbifold Stokes applies and small-sphere terms also vanish if the
calculation is performed by excision.

We apply the boundary formula of \cite[Proposition~2.1]{HL} on
the smooth AE end in complex dimension two to obtain
\[
 \int_{S_R}
 \sum_{i,j=1}^{4}
 (\partial_i h_{ij}-\partial_j h_{ii})\nu^j\,dA_e
 =
 2\int_{S_R}\theta\wedge\omega
 +O(R^{2-2\tau}).
\]
Here $\nu$ and $dA_e$ are the outward unit normal and area element
for the Euclidean background metric.
The preceding vanishing of $\int_{S_R}\theta\wedge\omega$ and
the assumption $\tau>1$ therefore imply
\[
 m(h)
 =
 \frac{1}{12\pi^2}
 \lim_{R\to\infty}
 \int_{S_R}
 \sum_{i,j=1}^{4}
 (\partial_i h_{ij}-\partial_j h_{ii})\nu^j\,dA_e
 =0.
\]

The metric is complete, its real scalar curvature is identically
zero and hence is nonnegative and integrable, its singularities
are finitely many isolated quotient points, and its AE order
satisfies
$\tau>(4-2)/2=1$. These are the hypotheses of
\cite[Theorem~1.9]{JV}. A positive constant change in the mass
normalisation leaves its zero set unchanged. That theorem
therefore implies that $(X,h)$ is isometric to $(\mathbb R^4,e)$,
without orbifold singularities. Since $J$ is
parallel, its image under this isometry is a constant orthogonal
complex structure, which can be identified with the standard
one by a linear Euclidean isometry.
\end{proof}

\begin{lemma}
[No concentration at regular points]
\label{lem:modelregular}
Let \((X,d,p)\) and \(\mathcal E\) be the orbifold and exceptional set from Proposition~\ref{prop:modelprecompact} and Lemma~\ref{lem:modelhierarchy}. For every \(x\in\mathcal R(X)\) and every sequence \(x_j\in X_j\) converging to \(x\) in the pointed GH sense,
\[
\lim_{r\downarrow0}\limsup_{j\to\infty}
\int_{B_{h_j}(x_j,r)}
|\operatorname{Rm}_{h_j}|^2\,dV_{h_j}=0.
\]
Moreover, no sequence of orbifold singular points of \(X_j\) converges to \(x\). Consequently,
\[
\mathcal E=\operatorname{Sing}X.
\]
\end{lemma}
\begin{proof}
By Lemma~\ref{lem:modelhierarchy}, the completed limit is an
orbifold. The local group at $x$ is
$
\Gamma_x=\{1\}.
$
For now, \(x\) might still belong to \(\mathcal E\): the definition of that set records concentration
in the sequence, rather than singularity of the completed metric.
\begin{claim}
Let
$
\nu_j:=|\operatorname{Rm}_{h_j}|^2\,dV_{h_j}$ and
$\nu_j\rightharpoonup\nu
$
on bounded regions. Then $
\nu(\{x\})=0.
$
\end{claim}
\begin{proof}
Suppose
$
\nu(\{x\})>0.
$
The curvature cannot satisfy a uniform bound
$
|\operatorname{Rm}_{h_j}|\leq K
$ in a fixed ball $B_{j,r}:=B_{h_j}(x_j,r)$.
Otherwise,
\[
\int_{B_{j,r}}
|\operatorname{Rm}_{h_j}|^2\,dV_{h_j}
\leq K^2Br^4\rightarrow 0,
\quad r\downarrow0,\] contradicting the positive atom.

Choose a small regular outer annulus around \(x\). Its source
annuli converge smoothly. Apply curvature point selection in the
region enclosed by this annulus, using it as a controlled outer
collar. The positive curvature atom forces the maximum of
$|\Rm|$ times squared distance to that boundary to tend to infinity.
Thus the selected points, rescaled by their curvature, have distance
to the boundary tending to infinity. Step~3 of
Lemma~\ref{lem:modelhierarchy} gives a complete nonflat scalar-flat
K\"ahler ALE limit at a shrinking scale, among the components
retained in that lemma.

Consider the finite chain leading from \(x\) to the selected nonflat component:
\[
x\longrightarrow Y_1\longrightarrow Y_2
\longrightarrow\cdots\longrightarrow Y_m.
\]
Here \(Y_1\) is the first component below \(x\), and \(Y_m\) is the selected nonflat component.
By group matching
$
\Gamma_\infty(Y_1)=\Gamma_x=\{1\}.
$
Therefore \(Y_1\) is AE. By Lemma~\ref{lem:modelnullclosed}, \(Y_1\) is null. By Lemma~\ref{lem:nullAE},
$
Y_1\cong\mathbb R^4.
$
In particular, every attachment point on \(Y_1\) has trivial local group.
Now let \(q_1\in Y_1\) be the attachment point of \(Y_2\). By group matching again,
$
\Gamma_\infty(Y_2)=\Gamma_{q_1}(Y_1)=\{1\}.
$
So \(Y_2\) is also AE and null, hence Euclidean.
Repeating this argument gives
$$
Y_1\cong Y_2\cong\cdots\cong Y_m\cong\mathbb R^4.
$$
But \(Y_m\) was chosen to be nonflat. This is the contradiction.
Thus the curvature measure has no atom at $x$.
\end{proof}
The argument includes flat intermediate components: they carry
the attachment data even though they carry no curvature energy.
For $x_j\to x$ and every sufficiently small continuity
radius $r$, GH compatibility and the upper bound for weak limits
on closed sets give
\[
 \limsup_{j\to\infty}
 \int_{B_{h_j}(x_j,r)}|\Rm_{h_j}|^2\,dV_{h_j}
 \le \nu\bigl(\overline{B_d(x,2r)}\bigr)\longrightarrow0, \quad r\downarrow0,
\]
since $\nu(\{x\})=0$.
If necessary this argument is applied to a subsequence realizing
the original limsup, so it proves the stated assertion for every
sequence of centers converging to $x$.

Finally suppose that singular points in $X_j$ converge to
$x$. Follow the corresponding marked points through the bubble trees. The same
induction makes every intervening complete component Euclidean.
At the final attachment, however, group matching would identify
the trivial group of a smooth Euclidean point with the nontrivial
group of the terminal quotient chart. This is impossible. A mark
attached directly to the original limit gives the same
contradiction without an intervening component. Hence no orbifold singular points can converge to $x$.
Thus the marked measure has no atom at any regular point. Since
$\mathcal E$ consists of its atoms of mass at least $\alpha$, and
$\operatorname{Sing}X\subseteq\mathcal E$, we obtain
$\mathcal E=\operatorname{Sing}X$.
\end{proof}

\subsection{Compatible charts and coverage}

The preceding lemma gives smooth convergence at every regular
limit point. We now construct charts compatible with the GH
identifications, following \cite[proof of Proposition~7.4]{LLW},
and show that their images contain the corresponding points of the underlying orbifolds.

\begin{proposition}
\label{prop:modelhatcharts}
The convergence in Proposition~\ref{prop:modelprecompact} satisfies Definition~\ref{def:modelhat}, after passing to a further subsequence. Moreover, its smooth comparison embeddings satisfy the coverage property on every compact subset of \(\mathcal R(X)\), see Section \ref{Convergence on the regular part}.
\end{proposition}
\begin{proof}
Fix $K\Subset U\Subset\mathcal R(X)$. Lemma~\ref{lem:modelregular},
\eqref{eq:modelEpsRegularity} and noncollapse give uniform curvature
derivative and injectivity-radius bounds on the source region
corresponding to a neighborhood of $\overline U$, with no source
singular points there. These bounds hold for the entire corresponding
region: a violating sequence would converge to a regular point where
the same local estimates apply.

Let $\psi_j^0$ be the original GH maps. Choose a finite coordinate
cover with source centers $\psi_j^0(x_a)$ and frames consistent with
the prescribed GH limit. The metrics and transition maps converge
smoothly. The atlas construction in
\cite[proof of Proposition~7.4]{LLW} combines the coordinate maps
$u_a$ and cutoffs $\eta_a$ into the Euclidean embedding
$x\mapsto(\eta_a(x),\eta_a(x)u_a(x))_a$.
Tubular projection between the resulting nearby submanifolds gives
embeddings $\Phi_j$ near $\overline U$ with
$\Phi_j^*h_j\to h$ smoothly. Lemma~\ref{lem:parallel-complex}
gives simultaneous convergence of $\Phi_j^*J_j$ to the complex
structure already obtained in Proposition~\ref{prop:modelprecompact}.
Because the finite atlas was based at GH-corresponding points,
\begin{equation}\label{eq:modelchartGHerror}
 e_j:=\sup_{\overline U}d_j(\Phi_j(x),\psi_j^0(x))\longrightarrow0.
\end{equation}
Replace $\psi_j^0$ by $\Phi_j$ on $K$. This changes its distortion by
at most $2e_j$, and its density and base-point errors by at most $e_j$.
Increasing the error slightly preserves
\eqref{eq:modelGHbase}--\eqref{eq:modelGHcover} and gives the exact
compatibility \eqref{eq:modelhatcompatible}.

For coverage, first take $U$ smoothly bounded. More generally, if
$\Phi$ is defined near $\overline U$ and $\Phi^*h_j\ge c^2h$, any
source path from $\Phi(x)$ to the complement of $\Phi(U)$ must first
reach $\Phi(\partial U)$. Pulling back to that first exit gives
length at least $c\,d_h(K,\partial U)$. Consequently,
\begin{align}\label{first-exit estimate}
 B_{h_j}(\Phi(x),r)\subset\Phi(U),\qquad
 x\in K,\quad0<r<c\,d_h(K,\partial U).
\end{align}
For $\Phi_j$ one may take $c=1/2$. If $x_j$ converges to $x\in K$
under the original GH identifications, then
$d_j(x_j,\Phi_j(x))\le d_j(x_j,\psi_j^0(x))+e_j\to0$.
Thus $x_j\in\Phi_j(U)$ for large $j$. For general $U$, insert a
smoothly bounded domain between $K$ and $U$. A diagonal argument on
an exhaustion completes the construction.
\end{proof}

\subsection{Proof of Theorem~\ref{thm:stationary}}
Proposition~\ref{prop:stationaryVolume} gives both volume bounds,
and Proposition~\ref{prop:modelprecompact} gives preliminary
pointed compactness. Lemmas~\ref{lem:oneend} and
\ref{lem:modelhierarchy} show that the completed limit is a
connected scalar-flat K\"ahler ALE orbifold with one
end. Its decay has every order below two. The homogeneous Sobolev inequality passes to the completed limit
and its normalized components by
Lemma~\ref{lem:sobolev-passage}. Lemma~\ref{lem:modelnullclosed}
passes nullness to the limit, using the period criterion in
Lemma~\ref{nullness condition}.
Lemma~\ref{lem:modelregular} and
Proposition~\ref{prop:modelhatcharts} give smooth convergence
throughout its regular part, with the required compatibility
and the coverage asserted in \ref{the convergence charts}.

It remains to check \eqref{eq:modelComplexity}.
Let \(\mu\) be the weak limit of the marked measures. Lemma~\ref{lem:modelregular} and Proposition~\ref{prop:modelhatcharts} give
\[
\mu\big|_{\mathcal R(X)}
=
|\operatorname{Rm}_h|^2\,dV_h.
\]
By Proposition~\ref{prop:modelprecompact},
\(\operatorname{Sing}X\subset\mathcal E\), and hence every orbifold singular point has \(\mu\)-mass at least \(\alpha\). Therefore
\[
\int_X|\operatorname{Rm}_h|^2\,dV_h
+\alpha\,\#\operatorname{Sing}X
\le \mu(X)
\le
\liminf_{j\to\infty}\mathscr E_\alpha(X_j,h_j)
\le\Lambda.
\]
Thus the limit belongs to the same model class. Natural volume convergence follows from Proposition~\ref{prop:modelprecompact}, completing the proof.

\subsection{Extension from the end and quantitative gaps}

The extension argument uses the one-convex structure of an ALE
K\"ahler surface \cite[Section~2.2]{HV}, followed by Hartogs
extension on normal Stein spaces \cite[Theorem~1.1]{Ruppenthal}.
This produces a global holomorphic map before the metric
identity is continued.

\begin{lemma}[Extension from the end]\label{lem:endextension}
Let \((X_i,h_i,J_i)\), \(i=1,2\), be null scalar-flat K\"ahler ALE orbifold surfaces, and let \(K_i\subset X_i\) be compact sets such that \(X_i\setminus K_i\) is connected. Then every biholomorphic isometry
\[
F_0:X_1\setminus K_1\longrightarrow X_2\setminus K_2
\]
extends uniquely to a global biholomorphic isometry
$
F:X_1\longrightarrow X_2.
$
Thus
\[
F|_{X_1\setminus K_1}=F_0,
\qquad F^*h_2=h_1.
\]
\end{lemma}
\begin{proof}
We first show that the underlying normal complex surfaces
$X_1,X_2$ are Stein.
In an ALE coordinate,
$J_i$ converges to $J_0$ \eqref{eq:modelALEdefinition}. So, $r^2$ is
strictly plurisubharmonic outside a compact set. Extending this
function over the compact core, we get a proper exhaustion which
is strictly plurisubharmonic outside a compact set. Thus each
$X_i$ is one-convex. Its Remmert reduction is a proper
bimeromorphic holomorphic map
\[
 \pi_i:X_i\longrightarrow Y_i
\]
onto a normal Stein surface $Y_i$, with connected compact fibers; see
\cite[Section~2.2]{HV}.

Nullness, in the cohomological form of Lemma~\ref{nullness condition},
excludes compact analytic curves.
Suppose \(C\subset X_i\) were a nontrivial compact analytic curve. Because \(C\) is complex one-dimensional, its regular part is a real two-dimensional surface, and the K\"ahler form restricts to its positive area form. Hence
$
\int_{C_{\mathrm{reg}}}\omega_i>0.
$
But the same integral computes the cohomological pairing by nullness:
\[
\int_{C_{\mathrm{reg}}}\omega_i
=
\langle[\omega_i],[C]\rangle
=
0.
\]
That is the contradiction
$
0=\int_{C_{\mathrm{reg}}}\omega_i>0.
$
Thus every fiber is
zero-dimensional. The proper map $\pi_i$ is consequently finite.
Since a finite bimeromorphic map onto a normal complex space is an
isomorphism, we know $X_i$ itself is Stein.

Embed $X_1$ and $X_2$ properly as closed complex subspaces of
finite-dimensional complex affine spaces. Every coordinate
function of the embedded map $F_0$ extends holomorphically to
$X_1$ by \cite[Theorem~1.1]{Ruppenthal}. Its hypotheses hold:
$X_1$ is connected, normal, Stein, of pure complex dimension
two, and has only isolated singularities, while
$X_1\setminus K_1$ is connected. Every entire function
vanishing on the embedded target pulls back to zero on this
exterior, and hence to zero on all of $X_1$ by the identity
theorem. The extended coordinate functions therefore define a
holomorphic map
\[
 F:X_1\longrightarrow X_2.
\]
Apply the same argument to $F_0^{-1}$ to obtain
$G:X_2\to X_1$. The identities $G\circ F=\operatorname{id}$
and $F\circ G=\operatorname{id}$ hold on the respective
exteriors, and thus everywhere. Consequently $F$ is a global
biholomorphism.

We next continue the metric identity. In local holomorphic coordinates, a K\"ahler metric has a potential:
$
h_{i\bar j}=\partial_i\partial_{\bar j}\varphi.
$
The scalar-flat equation is fourth order
\[
h^{i\bar j}\partial_i\partial_{\bar j}
\log\det(h_{k\bar\ell})=0.
\]
Its highest order part has the form
\[
h^{i\bar j}h^{k\bar\ell}
\partial_i\partial_{\bar j}\partial_k\partial_{\bar\ell}\varphi.
\]
Accordingly, the principal symbol is a nonzero constant multiple of
$
\left(h^{i\bar j}\xi_i\xi_{\bar j}\right)^2.
$
The positivity of the metric makes this nonzero for every nonzero real covector. Thus the equation is elliptic. Its dependence on the derivatives of \(\varphi\) is analytic on the region where the complex Hessian is positive definite.
Analytic elliptic regularity \cite{Morrey} makes the potential,
and hence the metric coefficients, real analytic.

Now consider the tensor
$
F^*h_2-h_1
$
on \(\mathcal R(X_1)\).
It is real analytic because \(F\) is holomorphic and both metrics are real analytic. It vanishes on the exterior where \(F=F_0\). Since the regular locus is connected, the real-analytic identity theorem yields
\[
F^*h_2=h_1
\qquad\text{on }\mathcal R(X_1).
\]
So the global biholomorphism is also an isometry on the entire regular locus.

Finally, near a quotient point a biholomorphism of the
punctured quotient neighborhoods lifts to their universal
uniformizing covers. The lifts and their inverses extend
holomorphically across the origins by Hartogs extension in
complex dimension two. Their inverse identities continue
across the origins as well, giving the required orbifold
biholomorphism. The lifted metric identities extend by
continuity. Uniqueness follows from the holomorphic identity
theorem.
\end{proof}
\begin{lemma}[Regularity and branching gaps]\label{lem:modelgaps}
Let $(X,h,J,p)\in \mathcal A(A,\Lambda,\alpha)$ be a pointed null
scalar-flat K\"ahler ALE orbifold belonging to this fixed
class in Definition~\ref{def:modelmoduli}.

Fix an integer $k\ge2$ and a number $0<\eta<1/10$.
Then the following assertions hold.

\begin{enumerate}[label=\textup{(\roman*)}]

\item \textbf{Euclidean comparison charts.}
There is 
$\zeta_E=\zeta_E(A,\Lambda,\alpha,k,\eta)>0$ with the following
property. If $(X,h,J,p)$ satisfies
\begin{equation}\label{eq:modelgapsEuclidean}
 d_{\mathrm{pGH}}\left(
  (\overline B_h(p,4),d_h,p),
  (\overline B_{\R^4}(0,4),d_e,0)\right)<\zeta_E,
\end{equation} then
the point $p$ belongs to $\mathcal R(X)$ and there is
a smooth embedding
\[
 \Phi:B_{\R^4}(0,3)\longrightarrow\mathcal R(X),
 \qquad \Phi(0)=p,
\]
such that
\[
 \|\Phi^*h-e\|_{C^k(B_{\R^4}(0,2))}
 +\|\Phi^*J-J_0\|_{C^k(B_{\R^4}(0,2))}<\eta,\qquad
 B_h(p,1)\subset\Phi(B_{\R^4}(0,2)).
\]
The $C^k$ norms are measured in the standard Euclidean
coordinates on $B_{\R^4}(0,2)$. In particular, the chart
is centered at $p$, controls both the metric and the
complex structure, and contains the entire metric ball
$B_h(p,1)$ in its image.

\item \textbf{Lower bounds for the regularity scales.}
Let $a_h$ and $\rho_h$ be the scales defined in
Section~\ref{sec:scr}, and let $N$ be the derivative
order used in their definitions. If $k\ge N+6$, there
is a constant
\[
 a_E=a_E(A,\Lambda,\alpha,N,k,\eta),
 \qquad 0<a_E\le1/8,
\]
such that every $(X,h,J,p)\in\mathcal A$ satisfying
\eqref{eq:modelgapsEuclidean} also satisfies
\[
 a_h(w)\ge a_E
 \quad\text{for every }w\in B_h(p,1/2),
 \qquad
 \rho_h(p)\ge a_E.
\]

\item \textbf{Separation from flat wedges.}
For every $G\ge1$ and every integer $N_B\ge2$, there
is a constant
\[
 \zeta_B=\zeta_B(A,\Lambda,\alpha,G,N_B)>0
\]
with the following property. Choose an integer $b$
and finite subgroups $\Gamma_\ell\subset U(2)$ such that
\[
 2\le b\le N_B,
 \qquad |\Gamma_\ell|\le G,
 \qquad 1\le\ell\le b,
\]
where each $\Gamma_\ell$ acts freely on the unit sphere.
Form the flat wedge
$
 W:=\bigvee_{\ell=1}^{b}(\R^4/\Gamma_\ell)
$
by identifying all cone vertices to a single point $o$,
and let $d_W$ be its induced length distance. Then,
for every $(X,h,J,p)\in\mathcal A$,
\[
 d_{\mathrm{pGH}}\left(
  (\overline B_h(p,4),d_h,p),
  (\overline B_W(o,4),d_W,o)\right)\ge\zeta_B.
\]
The same constant $\zeta_B$ works for every such
$(X,h,J,p)$ and every such wedge $W$.

\end{enumerate}

In all pointed GH distances above, each closed ball
carries the restriction of its ambient length distance.
All constants are uniform over the indicated class
and depend only on the listed parameters.

The assertions also apply at every scale $r>0$ by
replacing $h$ with $r^{-2}h$. In the original metric,
the relevant ball radii are multiplied by $r$, the
pointed GH thresholds are $r\zeta_E$ and $r\zeta_B$,
and the regularity lower bounds are $ra_E$.
The $C^k$ estimates are understood in normalized
coordinates for the metric $r^{-2}h$.
\end{lemma}
\begin{proof}
Suppose first that no $\zeta_E$ exists. Then there is a sequence
of models whose $\overline B_{h_j}(p_j,4)$ converge to the
$\overline B_{\R^4}(0,4)$ but for which the asserted charts
fail. By Theorem~\ref{thm:stationary}, a subsequence converges
in the pointed $\hat C^\infty$ topology to another model
$(X_\infty,h_\infty,J_\infty,p_\infty)$. Its open $B_{h_\infty}(p_\infty,4)$ is Euclidean. In particular every point of its smaller
closed balls has Euclidean tangent cone.

On the Euclidean limiting ball, the K\"ahler condition implies
the parallel complex structure $J_\infty$ is constant. Choose Euclidean coordinates which identify
it with $J_0$.
Once \(p_j\) belongs to the image, we can take its inverse coordinate:
$
v_j:=\Phi_j^{-1}(p_j).
$
Compatibility with the pointed GH convergence implies
$
v_j\rightarrow0.
$
Now precompose the chart with the small translation
$
x\longmapsto x+v_j.
$
The translated chart sends \(0\) exactly to \(p_j\). Because the original chart was defined on a neighbourhood of the closed radius-three ball, there is enough room to make this translation while retaining the required domain.
Thus the coverage in Proposition~\ref{prop:modelhatcharts} gives
\[
p_j\text{ lies in the chart}
\ \Longrightarrow\
\text{the chart can be centered at }p_j.
\]
After recentering the chart, the first-exit estimate from Proposition~\ref{prop:modelhatcharts} gives
\[
B_{h_j}(p_j,1)\subset\Phi_j(B(0,2)).
\]
This contradicts failure of the asserted charts.

The \(C^k\) chart estimate gives uniform constants \(C'_m\) such that
\[
|\nabla^m\operatorname{Rm}_h|\le C'_m
\quad\text{on }B_h(p,3/4),
\qquad 0\le m\le N+4.
\]
The controlled coordinate geometry on $B_h(p,3/4)$, with a margin
of at least $1/4$ around $B_h(p,1/2)$, also gives a uniform
injectivity-radius lower bound
$
\operatorname{inj}_h(w)\ge i_E>0,
\text{for }w\in B_h(p,1/2).
$

Choose \(a_E>0\) sufficiently small that
$
a_E\le\frac18,
2a_E<i_E,
$
and
\[
C'_m<a_E^{-2-m},
\qquad 0\le m\le N+4.
\]
For any \(w\in B_h(p,1/2)\),
$
B_h(w,2a_E)\subset B_h(p,3/4),
$
because \(2a_E\le1/4\).
This ball is contained in the regular coordinate region. It has the required curvature-derivative bounds, and the injectivity condition holds at \(w\).
Thus \(a_E\) is admissible in the definition of \(a_h(w)\), yielding
\[
a_h(w)\ge a_E
\qquad
\text{throughout }B_h(p,1/2).
\]
Now we estimate
$
a_h(w)+d_h(p,w)
$
for every \(w\in X\).
If \(w\in B_h(p,1/2)\), then \(a_h(w)\ge a_E\).
If \(w\notin B_h(p,1/2)\), then \(d_h(p,w)\ge1/2\ge a_E\).
In either case,
$
a_h(w)+d_h(p,w)\ge a_E.
$
Taking the infimum gives
$
\rho_h(p)\ge a_E.
$

For assertion \textup{(iii)}, there are finitely many abstract
groups of order at most $G$, and each has finitely many
complex two-dimensional unitary representations up to
equivalence. Restricting to the faithful free sphere actions
and taking at most $N_B$ branches leaves only finitely many
pointed wedge types.

If one fixed wedge with at least two branches could be approximated arbitrarily
closely by the indicated model balls, compactness of the
models would produce a model whose radius-four ball is that
wedge ball. The vertex of such a wedge has a disconnected
punctured neighborhood. An orbifold point has a
punctured neighborhood covered by a punctured four-ball, and
therefore has a connected punctured neighborhood. Thus the
limiting wedge vertex cannot occur in a model. There is a
positive separation for each wedge type, and the minimum of
these finitely many separations gives $\zeta_B$.
\end{proof}

\begin{remark}\label{rem:modelmotivation}
We follow Chen--Wang \cite{ChenWangI,ChenWangA,ChenWangB} from model moduli
to canonical radius and then to compactness of the flows.
Section~\ref{sec:fixedclass} identifies shrinking-scale Calabi
limits as scalar-flat ALE models. For any fixed real K\"ahler
class, the Hodge and lattice argument in
Lemmas~\ref{lem:fixedHodgeNorm}--\ref{lem:period} gives their
nullness.
\end{remark}

%%%%%%%%%%%%%%%%%%%%%%%%%%%%%%%%%%%%%%%%%%%%%%%%%%%%%%%%%%%%%%%%%%%%%%%%%%%%%%%%%%%%%%%%%%%%%%%%%%%%%%%%%%%%%%%%%%%%%%%%%%%%%%%%%%%%%%%%%%%%%%%%%%%%%%%%%%%%%%%%%%%%%%%%%%%%%%%%%%%%%%%%%%%%%%%%%%%%%%%%%%%%%%%%%%%%%%%%%%%%%%

\section{Canonical radius}
\label{sec:scr}

We use the models and gap estimates of Section~\ref{sec:stationary}
to define a space-time canonical radius. This adapts the
space-time canonical radius $\operatorname{scr}$ of Chen--Wang
\cite[Definition~2.7 and Remark~2.8]{ChenWangRemarks}.
Their definition incorporates two-sided Ricci-flow regularity on the
time scale $r^2$. Here we require model approximation at every smaller
scale and backward regularity on the Calabi-flow time scale $r^4$.
The comparison maps use points of the underlying manifold that are fixed
throughout the time interval.

First, an assumed lower bound for this radius gives stationary
limits at shrinking scales. Their estimates improve the defining
tests. Weighted point selection then gives a positive lower bound
along the original flow. The resulting short-time estimates will be
used in Section~\ref{sec:upgrade}.

\subsection{Definition and elementary properties}

\subsubsection{Choice of constants}
\label{Choice of constants}
Let $E_0$ be the curvature-energy bound in
Lemma~\ref{lem:energy}. Choose $A_*>2A$, larger also than the
Euclidean homogeneous Sobolev constant, and set
$$\Lambda:=2E_0+1,\qquad N:=12.$$ 
The stationary orbifold interior estimate \eqref{eq:modelEpsRegularity}
gives constants $\eps_{\rm stat}>0$ and $H_m$, $0\le m\le N$, such that
\begin{equation}\label{eq:modelinterior}
 \int_{B_h(z,7r/4)}|\Rm_h|^2\le\eps_{\rm stat}
 \quad\Longrightarrow\quad
 \sup_{B_h(z,5r/4)}|\nabla^m\Rm_h|
 \le H_mr^{-2-m}.
\end{equation}
Here the homogeneous Sobolev constant is at most $A_*$, and derivatives
at quotient points are interpreted in orbifold charts. Decrease
$\eps_{\rm stat}$ so that
$\eps_{\rm stat}\le\eps_{\rm reg}(A_*)$, and put
\[
 \eps_R:=\min\{1,\eps_{\rm stat}/4\}<\eps_{\rm stat},\qquad
 \alpha:=\eps_R/4,\qquad
 C_m:=4H_m+1>H_m,\qquad \theta:=1/16.
\]
We use the closed class $\mathcal A:=\mathcal A(A_*,\Lambda,\alpha)$ of
Theorem~\ref{thm:stationary}. Let $V_*$ be its uniform upper volume
ratio, and fix
\[
 B_0:=2\max\{1,\omega_4,V_*\}.
\]
The model comparisons below allow recentering and rescaling, including
recentering at a singular point.

For $z\in \mathcal R(Z)$, let $a_h(z)$ be the supremum of the $a>0$
such that $B_h(z,2a)$ is regular,
\[
 \operatorname{inj}_h(z)>2a,\qquad \sup_{B_h(z,2a)}|\nabla^m\Rm_h|<a^{-2-m},
 \qquad 0\le m\le N+4.
\]
Set $a_h:=0$ at singular points and define the regularity scale
\begin{equation}\label{eq:modelrho}
 \rho_h(z):=\inf_{w\in Z}\{a_h(w)+d_h(z,w)\}.
\end{equation}
For $0<r<\rho_h(z)$, the ball $B_h(z,2r)$ is regular, and the
injectivity radius and curvature satisfy the bounds defining $a_h(z)$
at scale $r$. On non-Euclidean models,
the infimum also gives
\[
\rho_h(z)\le a_h(z),
\qquad
|\rho_h(z)-\rho_h(z')|\le d_h(z,z').
\]
Thus $\rho_h$ is a $1$-Lipschitz regularity scale. It is finite,
positive on the regular locus, and zero at orbifold singular points.
Positivity follows from a uniform lower bound for $a_h$ in a
small regular neighborhood. Finiteness follows by evaluating the
infimum at a point where $a_h$ is finite. Such a point exists at an
orbifold singular point or at a point of nonzero curvature. On Euclidean space we use the convention
$\rho_h\equiv+\infty$, with inequalities interpreted in this extended
sense. The function scales as a length. In particular, a bounded
closed set on which $\rho_h$ has a positive lower bound is a compact
subset of the regular locus.

Apply Lemma~\ref{lem:modelgaps} to the class
$\mathcal A=\mathcal A(A_*,\Lambda,\alpha)$, with
differentiability order $k=N+6$ and with the parameter $\eta$
in that lemma set equal to $1/100$.
This gives constants $0<a_E\le1/8$ and $\zeta_E>0$ such
that, for every $(X,h,J,p)\in\mathcal A$,
\[
 d_{\mathrm{pGH}}\left(
  (\overline B_h(p,4),d_h,p),
  (\overline B_{\R^4}(0,4),d_e,0)\right)<\zeta_E
 \quad\Longrightarrow\quad
 \rho_h(p)\ge a_E.
\]
The value $1/100$ used in applying
Lemma~\ref{lem:modelgaps} is fixed before choosing
the constants $\eta$ and $\delta$ in
\eqref{eq:tolerances}.

Next, set
\[
 G_*:=\lceil G(A_*)\rceil,\qquad
 N_B:=\max\{2,\lceil B_0G_*/\omega_4\rceil\},
\]
where $G(A_*)$ is the group-order bound in
\eqref{eq:modelGroup}, and $\omega_4$ is the volume of the
Euclidean unit ball in $\R^4$.
Each local branch of a provisional flow limit inherits the
homogeneous Sobolev inequality with constant at most $A_*$.
The volume argument giving \eqref{eq:modelGroup} therefore
bounds its local quotient group by
$|\Gamma_\ell|\le G_*$.
At a joining point with $b$ branches, the inherited upper
volume bound gives
\[
 b\,\frac{\omega_4}{G_*}
 \le \sum_{\ell=1}^{b}\frac{\omega_4}{|\Gamma_\ell|}
 \le B_0.
\]
Consequently, $b\le N_B$. We use the upper volume bound
$B_0$ inherited from test \textup{(V)} \eqref{eq:scrV} in this count,
because membership of the provisional flow limit in
$\mathcal A$ has not yet been established.

Apply the branching assertion of Lemma~\ref{lem:modelgaps}
with $G=G_*$ and this value of $N_B$. Choose $\zeta_B>0$
so that, for every $(X,h,J,p)\in\mathcal A$ and every
flat wedge
\[
 W:=\bigvee_{\ell=1}^{b}(\R^4/\Gamma_\ell),
 \qquad 2\le b\le N_B,\qquad |\Gamma_\ell|\le G_*,
\]
one has
\[
 d_{\mathrm{pGH}}\left(
  (\overline B_h(p,4),d_h,p),
  (\overline B_W(o,4),d_W,o)\right)\ge\zeta_B.
\]
Here each $\Gamma_\ell\subset U(2)$ is a finite subgroup
acting freely on the unit sphere, $o$ is the common cone
vertex, and $d_W$ is the induced length distance.
The closed balls in these comparisons carry the
restrictions of their ambient distances.

If necessary, decrease $\zeta_E$ and $\zeta_B$ to account
for the fixed constants involved in passing from the
correspondence estimates on the larger comparison balls
to pointed GH estimates on radius-four balls. Retain the
same notation for the decreased constants, and then fix
\begin{equation}\label{eq:tolerances}
 \eta:=\min\{1/100,a_E/100\},\qquad
 \delta:=\min\{1/100,\eta/10000,\zeta_E/1000,\zeta_B/1000\}.
\end{equation}
The number $\eta$ is a dimensionless regularity threshold
in the normalized comparison model: lower bounds for
$\rho_h$ in terms of $\eta$ determine the regular regions
used in test \textup{(P)}. The number $\delta$ specifies
the allowed distance distortion and the allowed $C^N$
errors in comparing the metrics and complex structures
in that test.

\subsubsection{The three tests}
\label{sec:canonicaltests}
In a fixed manifold $M$, we assume that the flow $g(\cdot)$ exists on $[t-s^4,t]$. We write
\[
 G(\tau):=s^{-2}g(t+s^4\tau),\qquad -1\le\tau\le0,
 \qquad
 W_G(R):=\bigcup_{-1\le\tau\le0}\overline{B_{G(\tau)}(p,R)}.
\]
When comparing different times, points $x,y\in M$ are held fixed;
only the metric varies. For a sequence of flows, these points may
vary with the sequence index, but each is fixed with respect to time.
The following tests at $(p,t,s)$ are dimensionless.

\paragraph{(P): model comparison with fixed points of the
underlying manifold.}
There exist a pointed model $(Z,h,J_Z,o)\in\mathcal A$,
a compact relation $\mathscr C$, and a smooth embedding
$\Phi$ satisfying the conditions below. These objects
may depend on $(p,t,s)$, but are independent of
$\tau\in[-1,0]$.

\begin{enumerate}[label=\textup{(\roman*)}]

\item \textbf{Matching model points with points of $M$.}
The relation $\mathscr C$ is a compact subset
\[
 \mathscr C\subset
 \overline{B_h(o,9)}\times W_G(9),
\]
where compactness is taken in the product topology of
$Z\times M$. An element $(z,x)\in\mathscr C$ is an
ordered pair: its first entry is a point of the model
$Z$, and its second entry is a point of the underlying
manifold $M$. Membership in $\mathscr C$ declares that
these two points are matched.

A point may have more than one match. Once a pair
$(z,x)$ is chosen, both $z$ and $x$ remain fixed as
$\tau$ varies; distances involving $x$ are measured
using the changing metric $G(\tau)$.

We require $(o,p)\in\mathscr C$. We also require that
every $z\in\overline{B_h(o,8)}$ have at least one
partner $x\in W_G(9)$ with $(z,x)\in\mathscr C$,
and that every $x\in W_G(8)$ have at least one
partner $z\in\overline{B_h(o,9)}$ with
$(z,x)\in\mathscr C$.

\item \textbf{Comparison of distances.}
We require the following uniform distortion bound:
\begin{equation}\label{eq:scrPdistance}
 \sup_{\tau\in[-1,0]}
\operatorname{dis}_\tau(\mathscr C)<\delta,\quad 
\operatorname{dis}_\tau(\mathscr C)
:=
\sup_{\substack{(z,x),(z',x')\in\mathscr C}}
\left|
d_h(z,z')-d_{G(\tau)}(x,x')
\right|.
\end{equation}

\item \textbf{Smooth comparison on regular regions.}
Using $\rho_h$ from \eqref{eq:modelrho} and $\eta$
from \eqref{eq:tolerances}, define
\begin{align*}
 K&:=\{z:d_h(o,z)\le6,\ \rho_h(z)\ge\eta\},\\
 K^+&:=\{z:d_h(o,z)\le13/2,\ \rho_h(z)\ge\eta/2\},\\
 U&:=\{z:d_h(o,z)<7,\ \rho_h(z)>\eta/4\}.
\end{align*}
These are nested subsets of the regular locus:
$K\subset K^+\subset U$.

There is an open neighbourhood
$V\subset\mathcal R(Z)$ of $\overline U$ and a smooth
embedding
$
 \Phi:V\longrightarrow M
$
such that
\[
 (z,\Phi(z))\in\mathscr C
 \qquad\text{for every }z\in K^+.
\]
Thus, on $K^+$, the embedding specifies one of the
partners allowed by the relation.

The same embedding is used at every time, and it satisfies
\begin{align}
 \max_{0\le m\le N}
 \sup_{K^+\times[-1,0]}
 |(\nabla^h)^m(\Phi^*G(\tau)-h)|_h
 &<\delta,\label{eq:scrPmetric}\\
 \max_{0\le m\le N}
 \sup_{K^+}
 |(\nabla^h)^m(\Phi^*J_M-J_Z)|_h
 &<\delta.\label{eq:scrPJ}
\end{align}
The domain $V$ may be disconnected.

\item \textbf{Containment of small ambient balls.}
We require the \emph{capture condition}:
\begin{equation}\label{eq:scrcapture}
 B_{G(\tau)}(\Phi(z),\eta/500)
 \subset\Phi(B_h(z,\eta/100)),\quad \forall z\in K,\quad \forall \tau\in[-1,0].
\end{equation}
%It ensures that the entire ambient ball on the left is contained in the comparison chart.
\end{enumerate}
The capture condition follows from the preceding
domain and metric bounds. Indeed, if $z\in K$ and
$w\in\overline{B_h(z,\eta/100)}$, then
\[
 d_h(o,w)\le6+\eta/100<13/2,\qquad
 \rho_h(w)\ge\eta-\eta/100>\eta/2.
\]
Thus this closed ball lies in the interior of $K^+$.
On it, \eqref{eq:scrPmetric} gives
$\Phi^*G(\tau)\ge(1-\delta)h$.
A path starting at $\Phi(z)$ and leaving
$\Phi(B_h(z,\eta/100))$ must therefore have length
at least
\[
 \sqrt{1-\delta}\,\eta/100>\eta/500
\]
before its first exit. This proves
\eqref{eq:scrcapture}, which we record explicitly
as part of the test.

Finally, suppose $(z,x)\in\mathscr C$ with $z\in K$.
Since $(z,\Phi(z))\in\mathscr C$ as well,
\eqref{eq:scrPdistance}, applied to these two pairs,
gives
\[
 \sup_{\tau\in[-1,0]}
 d_{G(\tau)}(x,\Phi(z))
 <\delta<\eta/500.
\]
The capture condition then implies
\[
 x\in\Phi(B_h(z,\eta/100)).
\]
Consequently, $x$ has a unique model coordinate
$\Phi^{-1}(x)$ in this ball, independent of time.
In particular, if $o\in K$, the pair
$(o,p)\in\mathscr C$ ensures that $p$ lies in the
fixed comparison chart.

\paragraph{(R): backward regularity from small curvature energy.}
\begin{equation}\label{eq:scrR}
  \int_{B_{G(0)}(p,2)}|\Rm_{G(0)}|^2<\eps_R
  \quad\Longrightarrow\quad
 \sup_{x\in B_{G(0)}(p,1),\,-\theta\le\tau\le0}
 |\nabla^m\Rm_{G(\tau)}|(x)<C_m,
 \quad 0\le m\le N.
\end{equation}
The spatial set is the terminal ball, regarded as a fixed subset of
the source. In the original scale the backward length is
$\theta s^4$, and the bound is $C_ms^{-2-m}$.

\paragraph{(V): terminal upper volume.}
The terminal ball satisfies
\begin{equation}\label{eq:scrV}
 \Vol_{G(0)}B_{G(0)}(p,1)<B_0.
\end{equation}
A lower radius bound at every time makes this volume bound available
at every time.

\begin{definition}[Space-time canonical radius]
\label{def:scr}
The space-time canonical radius
\[
\scr(p,t)
:=\sup\{r>0:\text{(P), (R), (V) hold for every }0<s<r\},
\]
with the backward testing intervals $[t-s^4,t]$ also required to lie in the flow domain. Its value is zero if no such $r$ exists, and the value
$+\infty$ is allowed. All constants are fixed as above.
\end{definition}

Only scales strictly below the radius are tested. In particular,
$\scr\ge c$ supplies every scale $s<c$, regardless of whether the
supremum is attained. The argument does not require continuity of $\scr$.
It also treats limits in which the energy in (R) approaches $\eps_R$.
The definition also applies to stationary orbifold models, with
curvature derivatives taken in orbifold charts and comparison
embeddings taking values in the regular locus.

\begin{lemma}\label{lem:radiusproperties}
The space-time canonical radius has the following elementary properties:
\begin{itemize}
\item it is positive uniformly on each compact
time slab $M\times[a,b], 0<a<b<\infty$ with \(g(t)\) smooth throughout \([a,b]\);
\item it scales as a length
under Calabi scaling \eqref{Calabi rescaling};

\item for every $(X,h,J,p)\in\mathcal A(A_*,\Lambda,\alpha)$,
regarded as the time-independent Calabi flow
$g(t)\equiv h$ on $(-\infty,0]$, one has
\[
 \scr(p,t)=+\infty
 \qquad\text{for every }t\le0.
\] 
\end{itemize}
\end{lemma}
\begin{proof}
On a compact smooth slab, rescaled fixed normal coordinate charts
converge, in every prescribed finite order as the scale tends to zero,
to static Euclidean K\"ahler geometry. An outer coordinate ball of
normalized radius forty prevents minimizing curves between
radius-nine points from leaving the chart. This gives the ambient
distance comparison, a common relation containing the chart graph,
and ball containment by first exit. Normalized curvature derivatives tend
uniformly to zero, and the volume ratios tend to $\omega_4<B_0$.
The flow is defined on a fixed earlier interval, so all sufficiently
small backward testing intervals are available. Hence all tests hold
uniformly at small scales.

Every test is invariant under the normalization above. On a stationary
model, use the diagonal relation and the identity on the regular
locus for (P). The antecedent of (R) implies
\eqref{eq:modelinterior}, with the strict margin $H_m<C_m$;
(V) follows from $V_*<B_0$. Thus every scale passes.
\end{proof}

The tests impose estimates only through order $N=12$.
In Step~1 of Theorem~\ref{thm:selected}, we construct fixed holomorphic
charts from these bounds and prove the interior smoothing estimate
\eqref{eq:interiorsmoothing}. We restrict to smaller spatial domains
and the latter part of each backward interval when passing to smooth limits.

\subsection{Slice compactness and shrinking scales}
\label{sec:slicecompactness}

We first identify smooth local limits as extremal, then extract a
terminal slice under uniform backward regularity. The extraction
assumes this estimate explicitly; it does not use the global
canonical-radius lower bound proved later. For the selected sequences,
test (R) and interior smoothing will supply the estimate.

Vanishing dissipation identifies the limit on a regular backward
cylinder. We will use this before extending the metric across the
singular points.

\begin{proposition}[Extremality of a local limit]\label{prop:local}
Let $g_j(\tau)$ be Calabi flows on a common interval $$I_L:=[-L,0]$$
with fixed complex structures $J_j$.
Suppose there are time-independent comparison embeddings
$\Phi_j:U\to M_j$ such that
\begin{equation}\label{assumed convergence}
 \begin{aligned}
 \Phi_j^*g_j&\longrightarrow h
 &&\text{in }C^\infty_{\rm loc}(U\times[-\delta,0]),\\
 \Phi_j^*J_j&\longrightarrow J_\infty
 &&\text{in }C^\infty_{\rm loc}(U),
 \end{aligned}
\end{equation}
where $0<\delta\le L$ and $h(\tau)$ is positive definite.
The first convergence includes all mixed space and time derivatives.
Suppose that the spacetime dissipation tends to zero
\begin{equation}\label{spacetime dissipation 0}
 \int_{-L}^{0}D_j(\tau)\,d\tau\longrightarrow0,
\qquad D_j(\tau):=\int_{M_j}|\Dcal_{g_j(\tau)}S_j|^2\,dV_{g_j(\tau)},
\qquad S_j(\tau):=S_{g_j(\tau)}.
\end{equation}
Then $h(\tau)$ is a Calabi flow on $(U,J_\infty)$, and $h(\tau)$ is extremal for all $\tau\in [-\delta,0]$.
\end{proposition}
\begin{proof}
Apply Lemma~\ref{lem:parallel-complex} to
$(\Phi_j^*g_j(\tau),\Phi_j^*J_j)$.
Then \eqref{assumed convergence} identifies
the limiting complex structure with $J_\infty$.
Thus $(h(\tau),J_\infty)$ is K\"ahler for every
$\tau\in[-\delta,0]$.
Since $\Phi_j$ are independent of time,
\eqref{assumed convergence} also passes the Calabi
flow equation to $(h,J_\infty)$.
For $K\Subset U$ and $-\delta<a<b<0$,
\[
 \int_a^b\int_K|\Dcal_{h(\tau)}S_{h(\tau)}|^2\,dV_{h(\tau)}\,d\tau
 =\lim_j\int_a^b\int_{\Phi_j(K)}
 |\Dcal_{g_j(\tau)}S_{g_j(\tau)}|^2\,dV_{g_j(\tau)}\,d\tau=0.
\]
Exhaustion and continuity, including the temporal endpoints,
give extremality on $U\times[-\delta,0]$.
\end{proof}

We now use the elliptic estimates and singularity removal of
Chen--Weber \cite[Theorem~3.1 and Sections~3.3, 4]{CW}
to complete the terminal slice limit.

\begin{proposition}[Single-slice extraction from backward regularity]\label{prop:slice}
Let $(M_j,g_j(\tau),J_j,p_j)$ be connected complete pointed Calabi
flows of real dimension four on a common interval $I_L$,
where $L>0$. Assume the following uniform bounds:
\begin{enumerate}[label=\textup{(\roman*)}]
 \item The full curvature energy is at most $E_0$ on every slice.
 \item There are $A,b_0<\infty$, independent of $j$ and $\tau$,
 such that the uniform Sobolev inequality \eqref{eq:Sob} holds for $(M_j,g_j(\tau))$ with $(A,b_0)$.
 \item There are $r_v,\kappa_-,\kappa_+>0$ such that
 \begin{equation}\label{eq:extractionvolumes}
 \kappa_-r^4\le\Vol_{g_j(\tau)}B(x,r)\le\kappa_+r^4,
 \qquad x\in M_j,\ 0<r\le r_v.
 \end{equation}
 \item There are $r_*,\eps,\theta>0$ and constants $C_k$, $k\ge0$,
 independent of $j$, the center $x$, and the terminal time $t$, such that
 \begin{equation}\label{eq:slicebackward}
 \int_{B_{g_j(t)}(x,2r)}|\Rm|^2\,dV<\eps
 \quad\Longrightarrow\quad
 \sup_{B_{g_j(t)}(x,r)\times[t-\theta r^4,t]}
 r^{2+k}|\nabla^k\Rm|\le C_k
 \end{equation}
 for every $k\ge0$, whenever $0<r<r_*$ and
 $[t-\theta r^4,t]\subset I_L$.
 The spatial set is the fixed terminal ball.
 \item The spacetime dissipation tends to zero \eqref{spacetime dissipation 0}.
\end{enumerate}
Then a subsequence of the terminal slices converges in pointed GH
sense to a complete proper length space $(X,d,p)$, smoothly off a
finite provisional exceptional set $\mathcal E$. The regular
limit is extremal K\"ahler with finite curvature energy. Branchwise
singularity removal completes it to a smooth extremal K\"ahler
multifold. Points of $\mathcal E$ may be smooth after completion.
\end{proposition}
\begin{proof}
\noindent\textit{Step 1:}
The two-sided local volume bounds $(iii)$ give pointed GH precompactness.
After a diagonal extraction, terminal curvature-energy measures
converge weakly on bounded regions to $\mu$, with $\mu(X)\le E_0$.
For the threshold $\eps$ in \textup{(iv)}, set
$\mathcal E:=\{q:\mu(\{q\})\ge\eps\}$; then
$\#\mathcal E\le E_0/\eps$.

\noindent\textit{Step 2:}
For $q\notin\mathcal E$, choose $r>0$ with
$\mu(\overline{B_d(q,4r)})<\eps$, $4r<r_v$, $r<r_*$ and
$\theta r^4<L$. GH ball inclusions and weak convergence give
$\int_{B_{g_j(0)}(q_j,2r)}|\Rm|^2<\eps$ for corresponding $q_j$.
Assumption \textup{(iv)} supplies backward curvature bounds.
At time zero, curvature and noncollapse give injectivity-radius
bounds and time-independent comparison charts $\Phi_j:U\to M_j$
whose smaller images contain definite terminal balls.
Using the flow equation $|\partial_\tau g|_g\le C|\nabla^2\Rm|_g\le Cr^{-4}$,
for $0<\delta<\theta r^4$ we have on these fixed charts
\begin{equation}\label{eq:terminalmetricequivalence}
 e^{-C\delta r^{-4}}g(0)\le g(\tau)
                 \le e^{C\delta r^{-4}}g(0),
 \qquad -\delta\le\tau\le0.
\end{equation}
Connection variation, \textup{(iv)}, and the flow equation give
compactness in every mixed derivative on smaller cylinders.
Lemma~\ref{lem:parallel-complex} gives a time-independent limiting
complex structure, and Proposition~\ref{prop:local} gives extremality.
A diagonal exhaustion patches the terminal metrics and complex
structures, with coverage supplied by the smaller terminal balls.
By lower semicontinuity,
\begin{equation}\label{eq:limitingcurvatureenergy}
 \int_{X\setminus\mathcal E}|\Rm_h|^2\,dV_h\le E_0.
\end{equation}

\noindent\textit{Step 3:}
Coverage and \eqref{eq:extractionvolumes} allow
Lemma~\ref{lem:volume-passage} to pass the volume measures and
local volume bounds to the limit, without atoms at $\mathcal E$.
Lemma~\ref{lem:sobolev-passage} then gives the limiting Sobolev
inequality and its local homogeneous form.

\noindent\textit{Step 4:}
We now apply interior regularity to the limiting
metric. Let $x\in X\setminus\mathcal E$ and $d:=\operatorname{dist}_h(x,\mathcal E)>0$, with $d$ small. The ball
$B_h(x,d/2)$ avoids the deleted set, and its curvature energy
tends to zero as $d\downarrow0$, by
\eqref{eq:limitingcurvatureenergy}. For each fixed $k$, by the
interior estimate \cite[Theorem~3.1]{CW}, and using
\eqref{eq:localhomogeneoussobolev} in Lemma~\ref{lem:sobolev-passage}, we get
\begin{equation}\label{eq:puncturedecay}
 d^{2+k}|\nabla^k\Rm_h|(x)
 \le C_k\left(\int_{B_h(x,d/2)}|\Rm_h|^2\,dV_h\right)^{1/2}
 \longrightarrow0.
\end{equation}
Here the integral is the curvature energy of the regular limiting
metric. An atom of $\mu$ records energy lost from the
approximating sequence and is absent from this integral.

We apply the local cone construction in the proof of
\cite[Section~4]{CW} separately to each branch, using
the two sided volume bounds $(iii)$ and \eqref{eq:puncturedecay}.
On each preliminary punctured quotient cover, we use a radial
coordinate $r$, comparable to distance, and choose radial coordinate
curves with uniformly bounded speed. Set
\[
K:=J_\infty\nabla_hS_h,\qquad A_1:=\nabla K.
\]
We first establish the integrability needed for singularity removal.
This argument applies to the extremal limit.
Since $K$ is Killing, its transport equations along such a curve are
\[
\nabla_rK=A_1(\gamma'),\qquad
\nabla_r A_1=\Rm*\gamma'*K.
\]
For any small $\eps>0$, we shrink the neighborhood so $|\Rm|\le\eps r^{-2}$ and define
\[
H(r):=|K|^2+\eps^{-1}r^2|A_1|^2.
\]
Differentiation and completion of the square give
\[
rH'(r)\ge2\eps^{-1}r^2|A_1|^2-Cr|A_1||K|
\ge-C\eps H(r).
\]
The outer sphere is a compact smooth set. Integrating inward and choosing $\eps$ sufficiently small yields, for every prescribed $\delta>0$,
\[
|K|\le C_\delta r^{-\delta},\qquad
|\nabla K|\le C_\delta r^{-1-\delta}.
\]
The upper $r^4$ volume bound, applied on dyadic annuli, gives $K\in L^p_{\rm loc}$ for every finite $p$ and $\nabla K\in L^2_{\rm loc}$ by choosing $\delta<1$. For annular cutoffs $\chi_\eta$ with $|d\chi_\eta|\le C\eta^{-1}$,
\[
\int|d\chi_\eta|^2|K|^2
+\int|d\chi_\eta||K||\nabla K|
\le C\eta^{2-2\delta}\longrightarrow0.
\]
Thus $K,\nabla K\in L^2$ near the puncture, with the cutoff
control needed for the weak elliptic identities.
Since $\nabla J_\infty=0$, we have
\[
|dS_h|=|K|,
\qquad
|\nabla^2S_h|=|\nabla K|.
\]
Together with \eqref{eq:limitingcurvatureenergy}, this
establishes $S_h\in W^{2,2}$ locally across the puncture.
After shrinking the neighbourhood so that its curvature
energy is sufficiently small, we apply
\cite[Section~3.3]{CW} to obtain
$\Rm_h\in L^p$ for a prescribed finite $p>4$. Applying the interior estimate
\cite[Theorem~3.1]{CW} on balls of radius comparable to
$r=\dist_h(x,\mathcal E)$, H\"older's inequality and the upper
volume bound give $|\Rm_h|(x)\le Cr^{-4/p}=Cr^{-2+\sigma}$,
where $\sigma:=2-4/p>1$. The harmonic coordinate argument in
\cite[Section~4]{CW} then gives a smooth
extension of the metric on each finite cover.

By Lemma~\ref{lem:parallel-complex}, $J_\infty$ extends as a smooth
parallel complex structure over each uniformizing ball. The extended metric is K\"ahler and extremal by
continuity, and the local groups preserve the complex
structure, so they are conjugate into $U(2)$.
The inherited Sobolev inequality bounds each local group order, and
gives a volume lower bound on each branch. This lower bound, combined with the total volume upper bound, bounds the number of branches. Thus
the completion is a smooth extremal K\"ahler multifold. The
case of one branch with trivial group is an ordinary smooth
point. This proves the proposition.
\end{proof}

\begin{remark}
\label{rem:provisional-exceptional-set}
The capacity argument alone does not remove curvature singularities;
Step~4 also uses the integrability of the extremal field and
the singular extremal estimates. The provisional limit may still
have joining points or smooth points carrying lost curvature energy.
The shrinking-scale argument below excludes these possibilities
using the model gaps from Section~\ref{sec:stationary}.
\end{remark}

Before comparing the shrinking-scale limits at different times, we
explain how to retain the original comparison embeddings when new
charts are used to obtain estimates.

\begin{figure}[htbp]
\centering
\begin{tikzpicture}[x=1cm,y=1cm,font=\small]
 \node at (2.5,4.2) {Reference domains};
 \filldraw[fill=black!4,draw=black!55] (2.5,2) circle (1.8);
 \filldraw[fill=blue!9,draw=black!65] (2.5,2) circle (0.95);
 \node at (2.5,2) {$U$};
 \node at (3.45,3.25) {$U^+$};
 \draw[<->,black!65] (3.45,2) -- (4.3,2);
 \node[below,inner sep=4pt] at (3.87,2) {$>0$};
 \draw[->] (4.8,2) -- (6.7,2) node[midway,above] {$\Phi_j$};
 \node at (9.65,4.2) {Images in $M_j$};
 \path[fill=black!4,draw=black!55]
 (7.65,1.1) .. controls (7.2,2.1) and (7.8,3.7) .. (9.0,3.8)
 .. controls (10.0,4.0) and (11.4,3.6) .. (11.7,2.8)
 .. controls (12.25,1.7) and (11.8,0.7) .. (10.7,0.35)
 .. controls (9.6,0.0) and (8.1,0.3) .. (7.65,1.1) -- cycle;
 \path[fill=blue!9,draw=black!65]
 (8.7,1.4) .. controls (8.3,2.15) and (8.85,3.0) .. (9.6,3.0)
 .. controls (10.5,3.0) and (10.9,2.65) .. (10.9,2.05)
 .. controls (10.85,1.15) and (9.4,0.85) .. (8.7,1.4) -- cycle;
 \node at (9.6,2.05) {$\Phi_j(U)$};
 \node at (10.35,3.43) {$\Phi_j(U^+)$};
 \draw[<->,black!65] (10.91,2.05) -- (11.93,2.05);
 \node[below,inner sep=4pt] at (11.42,2.05) {$\ge cr$};
\end{tikzpicture}
\caption{The chart is controlled on $U^+$ and used on $U\Subset U^+$.}
\label{fig:nesteddomains}
\end{figure}
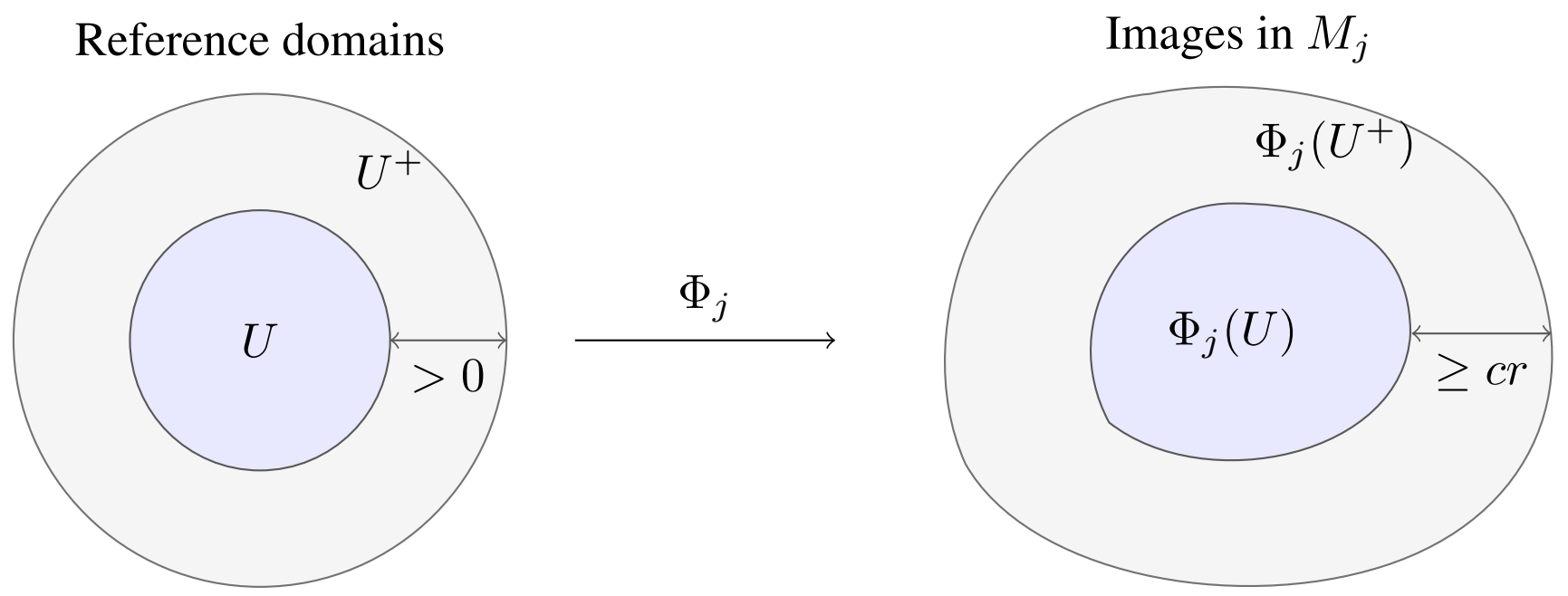

For local estimates we choose nested regular domains
\[
 U\Subset U^+\Subset \mathcal R(X).
\]
Charts and estimates are available on $U^+$, and the conclusion is
used on $U$; see Figure~\ref{fig:nesteddomains}. For a fixed smooth
reference metric $h$, compact containment gives
$d_h(\overline U,X\setminus U^+)>0$.
For comparison embeddings $\Phi_j:U^+\to M_j$, a uniform spatial
margin at scale $r$ means
\[
 d_{g_j(t)}\bigl(\Phi_j(\overline U),
                  M_j\setminus\Phi_j(U^+)\bigr)\ge cr,
\]
where $c>0$ is independent of $j$ and of $t$ in the stated interval.
A fixed separation in $h$ and uniform ellipticity of
$r^{-2}\Phi_j^*g_j(t)$ relative to $h$ imply this bound by
the first-exit argument.
Compact containment for each $j$ alone does not give a uniform margin.
For an interior parabolic estimate on $U\times[a,b]$, we use
estimates on $U^+\times[a-\sigma,b]$. A uniform past-time margin
at scale $r$ means $\sigma\ge cr^4$. These constants may depend
on the chosen regular domains and the data of the estimate.

The next lemma lets us keep the same comparison embeddings when
estimates are proved in another family of charts. At the terminal
time we use the compatible GH and smooth identifications of
\cite[Proposition~7.4]{LLW}. The proof explains why the change of
charts is independent of time.

\begin{lemma}[Retaining the comparison embeddings]\label{lem:retainmaps}
Let $I$ be a compact time interval, $b\in I$, and
$U\Subset V$ be regular domains in a fixed limit manifold.
Let $M_j$ have the same dimension as $V$. Suppose there are two time-independent embeddings
$\Phi_j,\Psi_j:V\to M_j$ such that their
pullback metrics  \[
a_j:=\Phi_j^*g_j(b)\to h_b, \qquad b_j:=\Psi_j^*g_j(b) \to h_b \qquad \text{ smoothly}.
\]
Here, \(h_b\) is a positive definite metric.
Suppose also that, for
every $K\Subset V$,
\begin{equation}\label{eq:retainmatch}
 \sup_{z\in K}d_{g_j(b)}(\Phi_j(z),\Psi_j(z))\longrightarrow0.
\end{equation}
If $\Psi_j^*g_j(t)\to h(t)$ smoothly on $V\times I$, locally in
space, then
\[
\Phi_j^*g_j(t)\to h(t)\text{ smoothly on }U\times I.
\]
The same assertion holds for additional smoothly convergent tensors.

If time-independent GH relations to the same limit family contain
the pairs $(\Psi_j(z),z)$ and have distortion tending uniformly
to zero on $I$, they may be enlarged to contain $(\Phi_j(z),z)$
for $z\in U$, with the same limiting distortion. Compatibility with the same GH identification
at $b$ implies \eqref{eq:retainmatch}.
\end{lemma}
\begin{proof}
Choose $U\Subset V'\Subset V$. The two boundary separations in
$h_b$ are positive; terminal smooth convergence makes the
corresponding margins uniform in $j$. Terminal
metric equivalence and first exit give a uniform ball about each
$\Psi_j(z)$, $z\in\overline U$, contained in $\Psi_j(V')$.
Equation~\eqref{eq:retainmatch} therefore defines
\[
 T_j:=\Psi_j^{-1}\circ\Phi_j\quad\hbox{on }U,
 \qquad T_j\longrightarrow\mathrm{id}\quad\hbox{in }C^0(\overline U).
\]
For the $C^0$ assertion, a sufficiently short source path cannot
leave $\Psi_j(V')$ when it starts in $\Phi_j(\overline U)$;
its coordinate length is controlled by
its metric length. The same argument works on a neighbourhood of
$\overline U$.

The identity $a_j=T_j^*b_j$ bounds $DT_j$ by ellipticity, and
\[
 \partial_r\partial_sT_j^a
 =\Gamma(a_j)^\ell_{rs}\partial_\ell T_j^a
  -\Gamma(b_j)^a_{cd}(T_j)
       \partial_rT_j^c\partial_sT_j^d
\]
bounds all higher derivatives by induction. Thus
$T_j\to\mathrm{id}$ smoothly on $\overline U$, using estimates on
a fixed larger domain compactly contained in $V'$.

Both embeddings are fixed, so $T_j$ and its domain are independent
of time. The identity
\begin{equation}\label{eq:retainedpullback}
 \Phi_j^*g_j(t)=T_j^*(\Psi_j^*g_j(t))
\end{equation}
proves the assertion, including every mixed derivative.

For the relations, the old pair associated to $\Phi_j(z)$ is
$(\Phi_j(z),T_j(z))$. Adding $(\Phi_j(z),z)$ increases the
distortion by at most
$2\sup_{z\in\overline U,\,t\in I}d_{h(t)}(T_j(z),z)=o(1)$.
Existing pairs, and hence the points covered by either projection, are retained. Finally, if
both graphs correspond to the same limit point $z$ at time $b$, comparison
of these two pairs gives \eqref{eq:retainmatch}.
\end{proof}

We now consider shrinking rescalings with a lower canonical-radius
bound. Proposition~\ref{prop:slice} gives limits at each time and
removes their singularities branchwise. The gap estimates place the
completed slices in the model class. To compare different times using
the same points of the underlying manifold, we prove distance continuity and extend the
isometries from a common stationary ALE end.
\begin{theorem}[Spacetime compactness at shrinking scales]
\label{thm:selected}
Let
$
(M_j,g_j(\tau),J_j,p_j), -T_j\le\tau\le0,
$
be smooth compact pointed Calabi flows, where $T_j\to\infty$.
Assume the following.
\begin{enumerate}
\item
There exist a fixed compact complex surface $(M,J)$, a fixed real
K\"ahler class
$
 \Omega\in H^2(M;\mathbb R),
$
constants $\lambda_j\to\infty$, and diffeomorphisms
$
 F_j:M\longrightarrow M_j
$
independent of $\tau$, such that
\[
 F_j^*J_j=J,\qquad
 \left[\lambda_j^{-1}F_j^*\omega_{g_j(\tau)}\right]=\Omega\quad
 \text{ for every }\tau\in[-T_j,0].
\]

\item
There is a constant $c>0$ such that, for every $L>0$ and all
sufficiently large $j$,
\[
 \scr_j(x,\tau)\ge c
 \qquad
 \text{for all }x\in M_j,\quad \tau\in I_L.
\]

\item
There are non-negative numbers $\beta_j\to0$, such that
for every $L>0$, all sufficiently large $j$, and every
$\tau\in I_L$, uniform energy and Sobolev bounds hold
\[
 \int_{M_j}|\Rm_{g_j(\tau)}|^2\,dV_{g_j(\tau)}
 \le E_0,
\qquad
 \|u\|_{L^4(g_j(\tau))}^2
 \le A\left(
 \int_{M_j}|\nabla u|_{g_j(\tau)}^2\,dV_{g_j(\tau)}
 +\beta_j\int_{M_j}u^2\,dV_{g_j(\tau)}
 \right)
\]
for every smooth function $u$ on $M_j$.
\item
For every $L>0$, the spacetime dissipation tends to zero
\eqref{spacetime dissipation 0}.

\end{enumerate}
Then, after passing to a single subsequence, the flows converge to
a pointed stationary model
\[
 (X,h,J_X,p)\in\mathcal A(A_*,\Lambda,\alpha)
\]
in the spacetime pointed-$\hat C^\infty$-Cheeger--Gromov topology of Definition~\ref{def:topology}
on every bounded backward time interval. The limiting flow is
$
 h(\tau)\equiv h, \forall \tau\in(-\infty,0],
$
and the basepoint $p$ may be singular.
More precisely, the following conclusions hold.
\begin{enumerate}[label=\textup{(\roman*)},itemsep=8pt]
\item
Pointed GH convergence holds on bounded pointed
regions, uniformly on every interval $I_L$, using comparison
relations independent of time. In particular, for every $L,R>0$,
\begin{equation}\label{eq:sourceconstant}
 \sup_{\tau\in I_L}
 \sup_{x,y\in B_{g_j(0)}(p_j,R)}
 \left|
 d_{g_j(\tau)}(x,y)-d_{g_j(0)}(x,y)
 \right|
 \longrightarrow0.
\end{equation}
Moreover, there are numbers
$\eps_j=\eps_j(L,R)\to0$ such that, uniformly for
$\tau\in I_L$,
\[
 B_{g_j(\tau)}(p_j,R)
 \subset B_{g_j(0)}(p_j,R+\eps_j),\qquad
 B_{g_j(0)}(p_j,R)
 \subset B_{g_j(\tau)}(p_j,R+\eps_j).
\]

\item
For every $K\Subset U\Subset \mathcal R(X)$ and $L>0$,
there are time-independent convergence embeddings $\Phi_j$,
defined on a regular neighborhood of $\overline U$, such that
\[
 \Phi_j^*g_j(\tau)\longrightarrow h
 \quad\text{in }C^\infty(K\times I_L),
\qquad
 \Phi_j^*J_j\longrightarrow J_X
 \quad\text{in }C^\infty(K).
\]
The first convergence includes all mixed spatial and time
derivatives. These embeddings are compatible with the pointed
GH approximation maps $f_{j,\tau}$ from the corresponding bounded
regions of $M_j$ to $X$. The maps may be chosen Borel, and the
embeddings satisfy uniform coverage on compact time intervals:
\[
 f_{j,\tau}^{-1}(K)\subset\Phi_j(U),
 \qquad \tau\in I_L,\quad j\ge j_0.
\]

\item
The volume measures converge to the natural orbifold volume
measure $dV_h$ on bounded pointed regions, uniformly on bounded
backward time intervals.

Ball volumes also converge at moving centers, times, and radii.
More precisely, if
\[
 \tau_j\in I_L,\qquad x_j\longrightarrow x,\qquad
 \rho_j\longrightarrow\rho>0
\]
under the comparison identifications, and the balls remain in a
bounded comparison region, then
\[
 \Vol_{g_j(\tau_j)}
 B_{g_j(\tau_j)}(x_j,\rho_j)
 \longrightarrow
 \Vol_h B_h(x,\rho).
\]

\item
Under the time-independent comparison relations, the limiting
singular spacetime set is
\[
 \Sing X\times(-\infty,0].
\]
\end{enumerate}
All these conclusions hold along the same subsequence on every
bounded pointed region and every bounded backward time interval.
\end{theorem}
\begin{proof}
Write $d_j(\tau):=d_{g_j(\tau)}$.
The assumed bound $\scr_j\ge c$ makes tests \textup{(P)},
\textup{(R)}, and \textup{(V)} of
Section~\ref{sec:canonicaltests} available at every center and time
under consideration and at every scale $s<c$.

\medskip
\noindent\textbf{Step 1: Extract a provisional limit at an arbitrary moving time.}
Choose any sequence $\tau_j\in I_L$ and set
\[
 \widetilde g_j(t):=g_j(\tau_j+t),\qquad -2\le t\le0.
\]
We apply Proposition~\ref{prop:slice} to their restrictions to
$[-1,0]$, verifying its hypotheses below. The larger interval
$[-2,0]$ supplies a uniform past-time margin for the backward
tests and local smoothing. These shifted intervals lie in
$[-L-2,0]$ for large $j$. Test \textup{(V)}, namely \eqref{eq:scrV}, gives
\begin{equation}\label{eq:selectedsmallvolume}
 \Vol_{g_j(\tau)}B(q,r)\le B_0r^4,
 \qquad 0<r<c.
\end{equation}
For each fixed $R<\infty$, take $j$ large enough that
$\beta_jR^2\le1$. The cutoff proof of \lemref{lem:noncollapse}
applied to the given Sobolev inequality gives
\[
 \Vol_{g_j(\tau)}B(q,r)\ge\kappa r^4,\qquad
 0<r\le R,\qquad \kappa=(20A)^{-2}.
\]
Thus we have required two-sided volume estimates $(iii)$
on a uniform small radius range
$r_v:=\min\{c/4,1\}$.

We prove the backward estimates $(iv)$ at every finite derivative
order by combining \textup{(R)} with the local smoothing argument
below. Test \textup{(R)} initially gives them only through order $N=12$.
Restrict $0<r<\min\{c,1\}/4$, so testing cylinders based in $[-1,0]$ remain inside $[-2,0]$. Suppose that at the terminal time
\[
 \int_{B(q,2r)}|\Rm|^2<\eps_R.
\]
For each $x\in B(q,r)$, apply \textup{(R)} at scale $r/4$.
Its energy hypothesis holds because
\[
 B(x,r/2)\subset B(q,2r).
\]
We obtain curvature derivatives through order twelve on a fixed
terminal neighborhood of $x$, backward for time $\theta(r/4)^4$.

Curvature bounds and noncollapse give terminal injectivity-radius
bounds. The metric and connection variation equations give finite
bounds in fixed charts over a smaller backward cylinder.
The parallel integrable complex structures and quantitative
Newlander--Nirenberg \cite[Theorem~4.5]{StreetNN} provide fixed local
holomorphic charts. Their normalized sizes, ellipticity bounds, and
spatial margins are uniform in $j$. The required starting norm is
below the defining order twelve.

\smallskip
\noindent\emph{Interior smoothing.}
We give the local argument for arbitrary flows satisfying finite
chart bounds, independently of the other hypotheses of this theorem.
Let $k_j(t)$ be smooth Calabi flows on a fixed holomorphic coordinate
ball $B_2\subset\C^n$, $-2\le t\le0$, with scalar curvatures $S_{k_j}$.
Suppose their matrices are
uniformly positive definite with uniformly bounded spatial $C^8$
norms throughout this cylinder. We claim that, for all integers
$a,b\ge0$,
\begin{equation}\label{eq:interiorsmoothing}
 \sup_{B_1\times[-1,0]}|D^a\partial_t^b k_j|\le C_{a,b}.
\end{equation}
Here $D$ denotes differentiation in the fixed coordinates.
The constants depend on the derivative orders, the original finite
bounds, ellipticity, and the fixed separations
$B_1\Subset B_2$ and $[-1,0]\subset(-2,0]$.
We use the local potential construction as in
\cite[Lemma~2.2 and the proof of Theorem~2.3]{LWZ} and interior
fourth-order parabolic Schauder estimates \cite[Chapter~VI]{EZ}.

On a smaller ball choose a real potential $\varphi_j(-2)$ with
controlled norm by fixing its pluriharmonic freedom, and put
$\varphi_j(t):=\varphi_j(-2)+\int_{-2}^tS_{k_j}(s)\,ds$.
Then $(\varphi_j)_{a\bar b}=(k_j)_{a\bar b}$ and
\begin{equation}\label{eq:potentialbootstrap}
\partial_t\varphi_j
=-k_j^{a\bar b}k_j^{c\bar d}(\varphi_j)_{a\bar b c\bar d}
+Q(k_j^{-1},D^3\varphi_j),
\end{equation}
where $Q$ is analytic on the positive-Hermitian-Hessian cone.
The scalar bound controls $\varphi_j$ in $C^0$, and the Euclidean
Laplacian identity $\Delta_E\varphi_j=\mathrm{const}\cdot
\operatorname{tr}_Ek_j$ gives interior spatial $C^{9,\gamma}$
bounds for each fixed $0<\gamma<1$. The flow equation bounds $\partial_t\varphi_j$ through
spatial order six and $\partial_tk_j$ through order four, hence
gives the initial parabolic H\"older coefficient bounds.

The principal symbol of \eqref{eq:potentialbootstrap}, with the
principal term moved left, is a positive multiple of $|\xi|_{k_j}^4$.
Differentiate spatially and apply interior fourth-order Schauder
estimates \cite[Chapter~VI]{EZ} on nested cylinders. Each improvement
of the potential improves the coefficients. Iteration gives every
spatial derivative, and the equation gives the time derivatives.
This proves \eqref{eq:interiorsmoothing}. The argument is local and
does not use the global curvature theorem of \cite{LWZ}.
The same estimates improve the holomorphic charts and their
transitions to every prescribed order.

Returning to $g_j$, apply \eqref{eq:interiorsmoothing} after rescaling
on smaller concentric terminal balls and the latter half of each
backward interval. It gives every prescribed derivative order,
with constants depending on the fixed spatial and temporal separations.

Repeating the construction at every $x\in B(q,r)$ controls the
whole target ball. Only the backward time coefficient is reduced,
for example to $\theta/(2\cdot4^4)$; the energy threshold remains
$\eps_R$.

The energy $(i)$ and Sobolev $(ii)$ bounds are hypotheses of \thmref{thm:selected}.
The dissipation on a shifted bounded interval $(v)$ is bounded by the
dissipation on the larger fixed interval $[-L-2,0]$ and hence
tends to zero. All hypotheses of \propref{prop:slice} are therefore
verified.
After taking a subsequence, we have a complete
proper pointed metric limit $(X,h,p)$, a weak curvature-energy measure
\[
 |\Rm_{g_j(\tau_j)}|^2\,dV_{g_j(\tau_j)}
 \rightharpoonup\mu,
\]
and a finite provisional exceptional set
\[
 \mathcal E:=\{q:\mu(\{q\})\ge\eps_R\},
 \qquad
 \#\mathcal E\le E_0/\eps_R.
\]
Away from $\mathcal E$, fixed local backward charts converge
smoothly and have coverage.
The branchwise completion of the limit is an extremal K\"ahler multifold.

Two possibilities remain open at this stage: several completed
branches might meet at a point, and a point of $\mathcal E$ might
become geometrically smooth after completion while still recording
source curvature concentration.

\lemref{lem:volume-passage} applies because the regular charts have coverage and
\eqref{eq:selectedsmallvolume} holds. It gives convergence to
natural volume. The exact upper constant $B_0$ passes to the limit:
verify ball-volume convergence at continuity radii and then
approach an arbitrary $r<c$ from above to obtain
\[
 \Vol_hB(q,r)\le B_0r^4.
\]

Since $\beta_j\to0$, smooth passage on compact regular supports
gives the homogeneous Sobolev inequality with coefficient $A$.
\lemref{lem:sobolev-passage} extends it across $\mathcal E$.
The same inequality, with the same coefficient, holds on each
normalized branch by zero extension and the capacity cutoff
argument.

Everything in this step applies at arbitrary moving times and
after arbitrary recentering.

\medskip
\noindent\textbf{Step 2: Show that every completed slice belongs to $\mathcal A$.}

We first eliminate the two provisional defects using
\textup{(P)} and \lemref{lem:modelgaps}.

\begin{claim}
    $X$ is an orbifold.
\end{claim}
\begin{proof}
The branchwise homogeneous Sobolev inequality bounds each local
group order by the Sobolev cutoff argument used in the proof of \propref{prop:stationaryVolume}:
$
 |\Gamma_\ell|\le G_*.
$
At a joining point the tangent volume ratio is
$
 \sum_{\ell=1}^{b}\frac{\omega_4}{|\Gamma_\ell|}.
$
The exact inherited upper bound gives
$
 b\,\frac{\omega_4}{G_*}
 \le\sum_{\ell=1}^{b}\frac{\omega_4}{|\Gamma_\ell|}
 \le B_0.
$
Hence $b\le N_B$.
This uses the provisional bound $B_0$, before membership in
$\mathcal A$ is known.
If $b\ge2$, choose a sufficiently small fixed scale $s<c$ at which
the completed neighborhood is arbitrarily close to the wedge
\[
 W:=\bigvee_{\ell=1}^{b}\R^4/\Gamma_\ell.
\]
At this scale, test \textup{(P)} compares the corresponding ball
in the underlying manifold with an allowed model. Restrict the
radius-eight comparison regions to radius-four balls.
Slice convergence and \eqref{eq:scrPdistance} then give a model
ball close to $W$, with error
\[
 O(\delta)+o_j(1)+o_s(1).
\]
The fixed conversion factor in $O(\delta)$ is absorbed in the
choice of $\delta$ in Section~\ref{Choice of constants}. Choose $s$ small and then $j$
large. The resulting error is below the branching gap $\zeta_B$,
contradicting \lemref{lem:modelgaps}. Thus each point has only one branch, so $X$ is an orbifold.
\end{proof}

\begin{claim}
The convergence is smooth with coverage on $\mathcal R(X)$.
\end{claim}
\begin{proof}
Suppose $q\in\mathcal E$ is smooth with energy concentration.
A sufficiently small fixed neighborhood of $q$, after rescaling,
is close to a Euclidean ball. Test \textup{(P)} makes its comparison
model sufficiently Euclidean as well. The Euclidean gap in
\lemref{lem:modelgaps} gives
\[
 \rho_h(o)\ge a_E>\eta
\]
at the comparison-model center, together with uniform regular
geometry nearby.

Let $q_j$ be the source center corresponding to $q$. Both
$(o,q_j)$ and $(o,\Phi_j(o))$ belong to the comparison relation, so
\[
 d_{G_j(0)}(q_j,\Phi_j(o))<\delta.
\]
The capture condition \eqref{eq:scrcapture} therefore places
a smaller ball around $q_j$ of normalized
radius $\eta/1000$, inside the controlled chart.
The whole smaller ball is controlled; a bound only at its center
would not exclude concentration. Model bounds
and the $C^N$ metric comparison \eqref{eq:scrPmetric} give, in
the scale before this normalization,
\[
 |\Rm_{g_j}|\le Cs^{-2}
\]
on a fixed-radius neighborhood of $q_j$. For sufficiently small
$u$,
\[
 \int_{B(q_j,u)}|\Rm|^2
 \le C^2s^{-4}\Vol B(q_j,u)
 \le C^2s^{-4}B_0u^4.
\]
Passing to the measure limit and then letting $u\downarrow0$
gives $\mu(\{q\})=0$, contradicting $q\in\mathcal E$.
Thus we prove the claim.
\end{proof}

The stationary compactness work of Section~\ref{sec:stationary} enters through
\lemref{lem:modelgaps}. Its Euclidean gap incorporates the earlier exclusion
of bubbles above regular points.

By the homogeneous Sobolev inequality, we now have a complete noncompact
extremal K\"ahler orbifold with finite curvature energy. \lemref{lem:scalarflat} applies and gives
\[
 S_h\equiv0.
\]
\propref{prop:ALE} gives ALE structure with positive differentiable
decay, and \lemref{lem:oneend} gives one end.
The fixed-class hypothesis supplies nullness through
\lemref{lem:fixedHodgeNorm}, \lemref{lem:transportedLattice},
and \lemref{lem:period}. On the fixed complex surface $(M,J)$, set
\[
 \omega_j(\tau):=\lambda_j^{-1}F_j^*\omega_{g_j(\tau)}.
\]
The hypothesis of \thmref{thm:selected} gives
$[\omega_j(\tau)]=\Omega$ for every $\tau\in[-T_j,0]$.
For a compact regular limiting two-cycle $\Sigma_0$, transport it
to $M_j$ by a convergence embedding and then to $M$ by $F_j^{-1}$.
Let $c_j$ be the resulting free integral homology class.
Smooth rescaled convergence gives
\[
 \lambda_j\langle\Omega,c_j\rangle
 \longrightarrow\int_{\Sigma_0}\omega_h.
\]
In particular, $\langle\Omega,c_j\rangle\to0$.

The essential strengthening is that $c_j$ ranges over a finite
set. \lemref{lem:fixedHodgeNorm} fixes the harmonic $L^2$ norms of cohomology classes
when $J$ and $\Omega$ are fixed. In real dimension four, these
two-form norms are invariant under constant metric scaling.
\lemref{lem:transportedLattice} pairs the pulled-back forms with
a fixed compactly supported Poincar\'e dual form for $\Sigma_0$.
Uniform metric equivalence on its support and Cauchy--Schwarz
bound all integer coordinates of $c_j$.
Its $\Omega$-period therefore lies in a finite set and is
eventually exactly zero.
\lemref{lem:period} yields
\[
 \int_{\Sigma_0}\omega_h=0.
\]
\lemref{nullness condition} identifies
regular-locus periods with the class on the completed orbifold.

Every geometric singular point belongs to the provisional
exceptional set and carries a $\mu$-atom of mass at least
$\eps_R$. Regular limiting curvature energy and these atoms occupy
disjoint parts of $\mu$. Since $\alpha=\eps_R/4$,
\begin{equation}\label{eq:flowcomplexity}
 \int_X|\Rm_h|^2\,dV_h+\alpha\,\#\Sing X
 \le\mu(X)\le E_0<\Lambda.
\end{equation}
This permits curvature loss into atoms and escape to infinity.
Every completed moving slice consequently belongs to
$\mathcal A(A_*,\Lambda,\alpha)$.
These are still separate slice limits; their points must next be
identified consistently through time.

\medskip
\noindent\textbf{Step 3: Establish short-time distance continuity for fixed points of the underlying manifold.}

Let $a_j\le b_j$ lie in a fixed bounded interval and suppose
$b_j-a_j\to0$. Extract the later slice at $b_j$ by Steps~1--2,
and denote its finite singular set by $\mathcal S$.

Choose small disjoint neighborhoods of radius comparable to
$\rho$ around the points of $\mathcal S$, with regular annular
collars. Also choose a large outer boundary in the ALE end.
Regular convergence with coverage transfers these collars to
fixed subsets of the source.

Full coverage of a thick annulus is needed to identify the
enclosed domains. Distance comparison places the inner and outer
source shells on opposite sides of the collar's middle link.
The radial coordinate on the smooth collar is continuous, so
every path joining the inner shell to the outer shell crosses
that link. Hence the inner component is a genuine small core,
contained in a terminal ball of radius $C\rho$, with the intended
single connected boundary. It cannot leave through an
unrepresented part of the source or connect to another core.
Call these fixed cores $U_{j,\alpha}$ and denote the large
enclosing domain by $W_j$.

On the regular exterior and its boundary collars, \textup{(R)}
and the interior estimate~\eqref{eq:interiorsmoothing} give
\[
 |\partial_\tau g_j|_{g_j}\le K_\rho
\]
over a common short backward interval. Thus
\[
 e^{-K_\rho|\tau-b_j|}g_j(b_j)
 \le g_j(\tau)
 \le e^{K_\rho|\tau-b_j|}g_j(b_j).
\]

Since $dV_{g_j(\tau)}=\omega_{g_j(\tau)}^2/2$ and
$\partial_\tau\omega_{g_j(\tau)}=dd^cS_j$, we have
\[
 \partial_\tau dV_j
 =dd^cS_j\wedge\omega_{g_j(\tau)}
 =d(d^cS_j\wedge\omega_{g_j(\tau)}).
\]
Stokes' theorem on a fixed core gives
\begin{equation}\label{eq:selectedflux}
 \frac{d}{d\tau}\Vol_{g_j(\tau)}U_{j,\alpha}
 =
 \int_{\partial U_{j,\alpha}}d^cS_j\wedge\omega_{g_j(\tau)}.
\end{equation}
The controlled boundary charts bound this flux by a constant
$F_\rho$. The terminal volume is $O(\rho^4)$ by
\eqref{eq:selectedsmallvolume}, so
\[
 \Vol_{g_j(\tau)}U_{j,\alpha}
 \le C\rho^4+F_\rho(b_j-a_j).
\]
For fixed $\rho$, the second term is at most another $C\rho^4$
for large $j$. Hence the core volumes remain $O(\rho^4)$ on
$[a_j,b_j]$. The boundary metric comparison also preserves their
$O(\rho)$ intrinsic boundary diameters.

Let $x$ lie in a core $U$. If its distance to the boundary
exceeded $r$, then $B(x,r)\subset U$, and noncollapse would give
$\Vol U\ge\kappa r^4$. Therefore every point is within
$(\Vol U/\kappa)^{1/4}$ of the boundary. Join two core points to
the connected boundary and then join along that boundary to get
\begin{equation}\label{eq:selectedcorediam}
 \operatorname{diam}^{\mathrm{int}}\overline U
 \le
 2(\Vol U/\kappa)^{1/4}
 +\operatorname{diam}^{\mathrm{int}}\partial U
 \le C\rho.
\end{equation}
For small $\rho$, only the already established small-scale
noncollapse is needed in this argument.

Collapse each core separately to a point. In the quotient length
spaces, all nonzero path lengths lie in the controlled exterior
and compare by the metric-velocity bound. Restoring a core
changes distances by at most its intrinsic diameter. A shortest
quotient path can be chosen to visit each collapsed core at most
once, so the total restoration error is at most
$C\,\#\mathcal S\,\rho$.

Choose the outer boundary far beyond the pointed region of
interest. The same quotient comparison bounds its earlier
distance from the basepoint from below. Relevant minimizing
paths therefore cannot leave the enclosing domain. This
first-exit argument also controls points in an earlier pointed
ball before a terminal distance bound for them is known.

For the same points of the underlying manifold in the union of the earlier and later
pointed radius-$R$ balls, we obtain
\begin{equation}\label{eq:selectedIDC}
 \begin{split}
 |d_{g_j(a_j)}(x,y)-d_{g_j(b_j)}(x,y)|
\le
 C_R\bigl(e^{K_\rho(b_j-a_j)/2}-1\bigr)
 +C\,\#\mathcal S\,\rho.
 \end{split}
\end{equation}
First let $j\to\infty$ with $\rho$ fixed, and then let
$\rho\downarrow0$. The distance difference tends to zero.
The order is important because $K_\rho$ and $F_\rho$ may grow
as $\rho$ decreases.

\medskip
\noindent\textbf{Step 4: Construct a stationary common end and extend it across each slice.}

Work first in fixed regular charts on
$D^+\times[-2\sigma,0]$, where $D\Subset D^+$ and $\sigma>0$
are independent of $j$. The finite defining bounds, followed by
the interior estimate~\eqref{eq:interiorsmoothing}, give the required higher estimates on
$D\times[-\sigma,0]$. We show that scalar-curvature derivatives
tend uniformly to zero on these smaller cylinders.
Otherwise one could
choose moving points and times where a prescribed derivative
stays bounded away from zero. Recenter and shift time there.
Steps~1--2 produce a scalar-flat completed limit, while the
existing controlled chart keeps the selected point regular.
Smooth convergence preserves the nonzero derivative, a
contradiction. Thus, on $D\times[-\sigma,0]$,
\[
 |\nabla^mS_j|\longrightarrow0
\]
for every prescribed $m$. The Calabi equation gives vanishing
metric velocity and its covariant spatial derivatives there.

Choose a fixed $0<s<\min\{c/100,1\}$ and set
\[
 \sigma:=\theta s^4/2.
\]
Far enough out in a terminal ALE model, all relevant terminal
$2s$-balls lie entirely in its smooth ALE end and have energy
less than $\eps_R/4$. On every fixed compact annulus, terminal
smooth convergence with coverage transfers a strict bound below
$\eps_R$ to the underlying manifold.
Test \textup{(R)} supplies backward cylinders of length
$\theta s^4$. Apply the interior estimate~\eqref{eq:interiorsmoothing} on concentric
$s/2$-balls and the latter half of each interval.
This gives smoothing over the common lifetime $\sigma$.

Exhaust the tail by compact annuli. The preceding vanishing
velocity observation, integrated in the retained terminal charts,
gives one stationary end in time-independent source charts.
Include a collar of the compact inner link. The starting radius
of the tail may depend on the terminal model, but the lifetime
$\sigma$ depends only on the fixed scale and testing constants.

The common end must stay in bounded pointed regions at earlier
times. Take a source domain enclosed by one of its regular
links. Its terminal volume is bounded by single-slice volume
convergence. Formula~\eqref{eq:selectedflux} bounds its volume on
the interval, and its boundary geometry is controlled.
The argument of \eqref{eq:selectedcorediam} bounds its
interior diameter. Here noncollapse is available on every fixed
bounded radius range for large $j$, because $\beta_j\to0$.
Thus the basepoint and each fixed end annulus remain at bounded
pointed distance throughout the interval.

Whenever an earlier or moving-time slice is extracted, the fixed
end charts give a biholomorphic local isometry from the common
end into that completed slice. We verify the hypotheses needed
to apply \lemref{lem:endextension}.

First, the end map is injective. Given distinct end points
$z,z'$, choose $a>0$ with $2a<d_h(z,z')$ and
$B_h(z,2a)$ contained in the regular end.
Metric control on $B_h(z,2a)$ and first exit from the image of
$B_h(z,a)$ give a fixed positive lower bound for the ambient
distance between their images in the underlying manifold.
Thus distinct end points cannot become identified in the limit.

Second, the map is proper at infinity. Otherwise a divergent
sequence of end points would have images in one bounded target
region. Choose disjoint fixed-radius regular balls about a
subsequence of these points. Their images are disjoint, have
uniformly positive volume, and lie in one bounded target ball,
contradicting its finite volume.

Third, the image contains an exterior region. Properness, local
openness, and the collar at the compact inner link show that
there is no additional finite boundary. By \lemref{lem:oneend}, the target
has one end. Outside a sufficiently large compact set, the image
is therefore a nonempty open and closed subset of that connected
end and contains it entirely.

After restricting to suitable connected exteriors, \lemref{lem:endextension}
extends the common-end isometry uniquely to a global
biholomorphic isometry of the completed slices.

Its nullness hypothesis was established in Step~2. In
\lemref{lem:endextension}, nullness excludes compact analytic curves and makes
the one-convex ALE complex surface Stein. Holomorphic extension
across the compact core, followed by analytic continuation of
the scalar-flat metric identity, gives the global isometry.
Uniqueness will align successive time intervals.

\medskip
\noindent\textbf{Step 5: Identify fixed points of the underlying manifold through time.}

The global isometries identify the slice geometries with a
reference model $X$. We must determine what happens to a fixed
source sequence $x_j$ in a bounded terminal region.

Choose a countable dense set of reference points $\{a_\ell\}$,
called \emph{anchors}, in an open regular part of the common end. Let $\Psi_j$ be the time-independent end charts
and put
\[
 f_{j,\ell}(\tau)
 :=d_{g_j(\tau)}(x_j,\Psi_j(a_\ell)).
\]
The enclosing-domain estimates bound these functions.
Equation~\eqref{eq:selectedIDC} gives asymptotic
equicontinuity. After diagonal extraction over $\ell$, they
converge uniformly to continuous functions $f_\ell$.

At any time, moving-time slice compactness and the global end
extension represent these values as
\[
 f_\ell(\tau)=d_h(\gamma(\tau),a_\ell)
\]
for some point $\gamma(\tau)\in X$.

Distances to an open regular set separate points. Indeed, if
two points have the same distance functions on that set, choose
an anchor where both functions are differentiable. Their equal
values and gradients specify the same initial direction and
length of a minimizing geodesic, hence the same endpoint.
A minimizing geodesic cannot pass through a nontrivial isolated
orbifold vertex in its interior, because its spherical quotient
link has diameter strictly less than $\pi$. Equality at the
dense anchors implies equality on the open set by continuity.

Distance to one anchor bounds all candidate points in a compact
subset of the proper space $X$. On this compact set, the map
to the countable collection of anchor distances is a continuous
injection with continuous inverse onto its image. Thus
$\gamma$ is uniquely determined and continuous.

Suppose $\gamma(\tau_0)$ is regular. Join its neighborhood to an
open subset of the common end by a \emph{regular tube}, meaning a
connected open set with compact closure in $\mathcal R(X)$.
Terminal regular convergence, \textup{(R)}, and smoothing give
fixed backward charts on this tube. Their limits are stationary
by Step~4. Retain a larger regular tube containing the closure
of the first. Within it, paths to the
fixed end anchors have uniformly bounded length in these charts.
Their images therefore remain in a bounded pointed region at
every time in this backward interval. At each earlier time,
first exit from the larger tube bounds the ambient distance
from the smaller tube to the chart boundary below. The smooth
charts consequently cover a whole ambient neighborhood of each
image point. Steps~1--2 identify these points as regular in the
earlier completed slice and give a local isometry from the tube
into that slice.

On the overlap with the common end, \lemref{lem:retainmaps} aligns the new
charts with the retained end charts through transition maps
converging smoothly to the identity. Those transition maps are
independent of time. The \emph{end-normalized identification} is the global isometry
extending the prescribed identification on the common end.
The tube identification agrees with it on an open set and hence
throughout the connected tube.

Coverage puts $x_j$ inside this fixed chart. Its inverse chart
coordinate is independent of time and converges at $\tau_0$.
Stationary convergence in the same chart therefore gives
\[
 \gamma(\tau)=\gamma(\tau_0)
\]
immediately backward from $\tau_0$.

This local backward constancy forces constancy on each
connected component of the \emph{regular-time set}
$\{\tau:\gamma(\tau)\in\mathcal R(X)\}$. If two values
differed, take the first time at which the distance from the
earlier value reaches a fixed positive intermediate value.
Backward constancy at that time contradicts its definition.
Such a constant regular component cannot end at an interior
time: continuity would give the same regular value at the
endpoint, contradicting maximality of the component.

If $\gamma$ is regular anywhere, it is consequently constant
on the whole interval. If it is regular nowhere, its continuous
image lies in the finite set $\Sing X$ and is again constant.
Sequences of points of the underlying manifold held fixed in time, including those approaching
singular points, therefore have time-independent limits.

\medskip
\noindent\textbf{Step 6: Iterate the identification and construct the convergence data.}

The lifetime $\sigma$ is uniform, so finitely many overlapping
intervals cover any prescribed $I_L$. At each restart,
normalize the new end identification to agree with the preceding
one at an interior overlap time. Choose fixed regular domains
$A\Subset A^+$ in the common end, with both embeddings defined on
$A^+$. Applying \lemref{lem:retainmaps} retains the earlier end
embedding on $A$ through a time-independent transition.
Uniqueness in \lemref{lem:endextension} makes the global
normalizations agree.
By Step~5, the limits of fixed points of the underlying manifold
are constant on both intervals and agree throughout their overlap.

This proves the uniform distance assertion \eqref{eq:sourceconstant}, namely
\[
 \sup_{\tau\in I_L}
 \sup_{x,y\in B_{g_j(0)}(p_j,R)}
 |d_{g_j(\tau)}(x,y)-d_{g_j(0)}(x,y)|
 \longrightarrow0.
\]
For completeness, if uniformity failed, we could choose moving
times $\tau_j$ and terminally bounded points $x_j,y_j$ with a
positive distance discrepancy. Extract their terminal limits.
The preceding construction gives the same limiting points
throughout time, and Step~3 handles the moving times.
Both distances then tend to the same number, a contradiction.

The outer first-exit comparison also supplies mutual
containment of moving and terminal pointed balls, with radius
errors tending to zero.

\smallskip
\noindent\emph{Time-independent comparison relations.}
Choose a terminal relation $\mathscr C_j^0$ compatible
with the terminal regular embeddings and containing $(p_j,p)$.
For fixed $R,L$, arrange that its projections cover both terminal
radius-$(R+1)$ balls and its source projection lies in
$B_{g_j(0)}(p_j,2R+4)$.

For any two pairs $(x,z),(y,w)\in\mathscr C_j^0$,
\[
 |d_{g_j(\tau)}(x,y)-d_h(z,w)|
 \le |d_{g_j(\tau)}(x,y)-d_{g_j(0)}(x,y)|
+|d_{g_j(0)}(x,y)-d_h(z,w)|.
\]
The first term tends uniformly to zero by
\eqref{eq:sourceconstant}; the second tends to zero
by terminal GH convergence. Thus the same terminal relation
works throughout the interval.

Mutual ball containment ensures that its source projection
covers the required moving balls. The target metric is
stationary, and the regular chart graphs remain included.
Exhausting bounded regions and intervals gives the relations
required by \defnref{def:topology}.

\smallskip
\noindent\emph{Smooth convergence in the original terminal embeddings.}
Fix regular domains
\[
 K\Subset U\Subset U^+\Subset\mathcal R(X)
\]
and retain a terminal smooth embedding $\Phi_j$ on a neighborhood
of $\overline{U^+}$. The domains are fixed, and their positive
separations are measured in $h$. Join them to the chosen regular
part of the common end by a regular tube if needed.
Step~5 ensures that the points $\Phi_j(z)$ of the underlying
manifold correspond to the same regular points $z$ at every time.
For smoothing on $U\times I_L$, use bounds on the larger cylinder
$U^+\times[-L-\sigma,0]$, keeping a fixed positive past-time
separation.

A moving-point contradiction using Steps~1--2 gives uniform
intrinsic curvature-derivative bounds there and
$|\nabla^mS_j|\to0$. Put
\[
 q_j(\tau):=\Phi_j^*g_j(\tau),\qquad
 T_j(\tau):=\partial_\tau q_j(\tau).
\]
Then, for every $m$,
\[
 |\nabla^mT_j|_{q_j}\longrightarrow0
\]
uniformly on smaller regular sets and bounded time intervals.

Integrating the zeroth-order estimate gives uniform metric
equivalence to $q_j(0)$ and $C^0$ convergence to $h$.
For higher coordinate derivatives, use
\[
 \partial_\tau\Gamma(q_j)^k_{ab}
 =\frac12q_j^{k\ell}
 \bigl(
 \nabla_a(T_j)_{b\ell}
 +\nabla_b(T_j)_{a\ell}
 -\nabla_\ell(T_j)_{ab}
 \bigr).
\]
Fix a smooth background connection $D$ and write
$A_j:=\nabla^{q_j}-D$. Schematically,
\[
 Dq_j=A_j*q_j,\qquad
 \partial_\tau A_j=q_j^{-1}*\nabla^{q_j}T_j.
\]
Integration first controls $A_j$. Differentiation and induction
control all its background derivatives: the estimate for
$D^rA_j$ uses intrinsic derivatives of $T_j$ through order
$r+1$ and already controlled derivatives of $A_j$ of lower
order. Since all the intrinsic velocity derivatives tend to
zero,
\[
 q_j(\tau)-q_j(0)\longrightarrow0
\]
in every spatial derivative order, uniformly in time.

The pulled-back complex structures are fixed in time and
already converge at time zero. The Calabi equation and the
smoothing bounds then give convergence of all mixed spatial
and time derivatives.

At each restart, take fresh regular domains $V\Subset V^+$
covering the same compact regular set and derive estimates on
$V$ from those on $V^+$. Their positive terminal separations are
independent of $j$. The time-independent terminal transitions
from \lemref{lem:retainmaps} transfer these estimates back to
$\Phi_j$, so the retained domain does not shrink with the
number of restarts.
First exit on the larger retained domain gives coverage of
all points corresponding to $K$.
For GH approximation maps $f_{j,\tau}$ chosen compatibly with
the comparison relations, this is precisely
\[
 f_{j,\tau}^{-1}(K)\subset\Phi_j(U)
 \qquad\text{for every }\tau\in I_L,\quad j\ge j_0.
\]
This proves the smooth convergence and uniform coverage required
by \thmref{thm:selected} and \defnref{def:topology}.

The constancy of limits of points of the underlying manifold held fixed in time also gives the
fixed singular spacetime set
\[
 \Sing X\times(-\infty,0].
\]

\medskip
\noindent\textbf{Step 7: Prove uniform natural-volume convergence.}

Apply the cutoff argument of \lemref{lem:volume-passage} uniformly in time.
Outside small neighborhoods of the finite singular set, the
retained embeddings give smooth uniform convergence of volume
forms. Inside these neighborhoods,
\eqref{eq:selectedsmallvolume} bounds the source volume by
\[
 C\,\#\mathcal S\,\rho^4
\]
uniformly in time. The limiting orbifold volume has the same
small-neighborhood bound.

For a continuous compactly supported test function, the discrepancy of
integrals is therefore bounded by a regular-region error
tending to zero with $j$, plus a constant times $\rho^4$.
First send $j\to\infty$ and then $\rho\downarrow0$.
This proves weak convergence to natural orbifold volume,
uniformly on every bounded backward interval. Curvature-energy
atoms may remain at singular points, while natural volume
has no atoms there.

Finally, uniform distance convergence gives inner and outer
ball inclusions for moving centers and radii. Distance spheres
in a Riemannian orbifold have zero natural volume, so these
inclusions imply ball-volume convergence. In particular, if
$\tau_j\in I_L$, $x_j$ corresponds to $x$, and
$\rho_j\to\rho>0$, with the balls in a bounded comparison region,
then
\[
 \Vol_{g_j(\tau_j)}
 B_{g_j(\tau_j)}(x_j,\rho_j)
 \longrightarrow
 \Vol_hB_h(x,\rho).
\]

A diagonal selection over bounded backward intervals, pointed
regions, regular compact sets, and derivative orders gives one
subsequence with all four conclusions of \thmref{thm:selected}.
\end{proof}

\subsection{Improvement of the defining estimates}

The selected limit is a model for which the defining estimates have
strict margins. We now transfer these margins back to the flows.
For (R), the terminal energy may approach the threshold in the
limit. At a testing scale $s$, applying the inherited estimate
at every center of $B(p,15s/8)$, inside the terminal
energy-controlled ball $B(p,2s)$, resolves this issue.

\begin{lemma}[Improvement of the defining estimates]
\label{lem:radiusimprovement}
Let $(M_j,g_j(\tau),J_j,p_j)$ satisfy the hypotheses
of Theorem~\ref{thm:selected}. Suppose that, for every
$L>0$ and all sufficiently large $j$,
\[
 \scr_j(q,\tau)\ge\frac12
 \qquad
 \text{for every }q\in M_j,\quad -L\le\tau\le0.
\]
Then
$\scr_j(p_j,0)\ge2$ for all sufficiently large $j$.
In particular, such a sequence cannot also satisfy
$\scr_j(p_j,0)=1$ for every $j$.
\end{lemma}

\begin{proof}
Suppose that $\scr_j(p_j,0)<2$ for infinitely many $j$,
and restrict to those indices.
For large $j$, all testing intervals associated with
scales below $2$ lie in the flow domain, since
$T_j\to\infty$.
The assumed radius bound makes all three tests valid
at every scale below $1/4$. Since
$\scr_j(p_j,0)<2$, there must therefore be a scale
$
 s_j\in[1/4,2)
$
at which at least one of \textup{(P)}, \textup{(R)},
or \textup{(V)} fails at $(p_j,0)$.
After passing to a subsequence, we may assume that
\[
 s_j\longrightarrow s\in[1/4,2]:=I
\]
and that the same test fails for every $j$.
By Theorem~\ref{thm:selected}, a further subsequence
converges to a stationary model
$
 (X,h,J_X,p)\in\mathcal A(A_*,\Lambda,\alpha).
$
We use the convergence on $[-32,0]$ and on sufficiently
large bounded pointed regions. This covers all testing
intervals $[-s_j^4,0]\subset[-16,0]$ and all spatial
comparison regions, whose outer radii are at most $18$.
We show that none of the three tests can fail.

\smallskip
\noindent\emph{Improvement of \textup{(P)}.}
Fix $\sigma\in I$ and define
$G(\tau):=\sigma^{-2}g_j(\sigma^4\tau)$, $-1\le\tau\le0.$
For the model required by test \textup{(P)}, take
$
 (Z,h_\sigma,J_Z,o)
 :=(X,\sigma^{-2}h,J_X,p),
$
which belongs to $\mathcal A$ by its constant rescaling invariance. We suppress the dependence on
$j$ and $\sigma$ in the notation below.

As in test \textup{(P)}, we define $W_G(R)$ and $U,K^+,K$ by using $p_j$ and $h_\sigma$ instead.
Since $\rho_{h_\sigma}=\sigma^{-1}\rho_h$, every
$\overline U$ arising for $\sigma\in I$ lies in
the common compact regular set
$
 Q:=\{z:d_h(p,z)\le14,\
          \rho_h(z)\ge\eta/16\}.
$
The spacetime convergence therefore provides a single
time-independent embedding $\Phi_j$, defined on a
regular neighborhood of $Q$, which can be used for
every $\sigma$. We take $\Phi:=\Phi_j$.

We next construct the compact relation required in
\textup{(P)}. Choose the compatible convergence
relations
$
 \mathscr R_j\subset X\times M_j
$
on $[-32,0]$ and sufficiently large pointed regions.
They contain $(o,p_j)$ and the graph pairs
$(z,\Phi(z))$ on the common regular neighborhood.
Their unscaled distance errors are bounded by
$\epsilon_j\to0$. After normalization, every two
pairs $(z,x),(z',x')\in\mathscr R_j$ satisfy
\begin{align}\label{verify P dt}
 \left|
 d_{h_\sigma}(z,z')
 -d_{G(\tau)}(x,x')
 \right|
 \le\frac{\epsilon_j}{\sigma},
 \qquad -1\le\tau\le0.
\end{align}
Because $\sigma\ge1/4$, this error is at most
$4\epsilon_j$, uniformly in $\sigma$.

Define
$
 \mathscr C:=
 \overline{\mathscr R_j\cap
 \bigl(\overline{B_{h_\sigma}(o,9)}
       \times W_G(9)\bigr)},
$
where the closure is taken in $Z\times M_j$.
The containing product is compact, so $\mathscr C$
is compact. It contains $(o,p_j)$.
We verify the two projection requirements.
For large $j$, $\epsilon_j/\sigma<1$ uniformly
in $\sigma$.
If $z\in\overline B_{h_\sigma}(o,8)$, choose a
partner $x$ in $\mathscr R_j$. We compare
$(z,x)$ with $(o,p_j)$ at $\tau=0$ by \eqref{verify P dt},
$
 d_{G(0)}(p_j,x)
 \le d_{h_\sigma}(o,z)+\frac{\epsilon_j}{\sigma}<9.
$
Thus $x\in W_G(9)$, and this pair belongs to
$\mathscr C$.
Conversely, if $x\in W_G(8)$, then
$d_{G(\tau)}(p_j,x)\le8$ at some
$\tau\in[-1,0]$. Choose a partner $z$ in
$\mathscr R_j$. At that time,
$
 d_{h_\sigma}(o,z)
 \le d_{G(\tau)}(p_j,x)
       +\frac{\epsilon_j}{\sigma}
<9.
$
Hence this pair also belongs to $\mathscr C$.

For $z\in K^+$, compatibility supplies
$(z,\Phi(z))\in\mathscr R_j$ and
$
 d_{G(0)}(p_j,\Phi(z))
 \le13/2+\frac{\epsilon_j}{\sigma}<9.
$
Thus
$
 (z,\Phi(z))\in\mathscr C$
for every $z\in K^+.$ Since $\frac{\epsilon_j}{\sigma}
 \le4\epsilon_j<\delta$, \eqref{verify P dt} gives \eqref{eq:scrPdistance} for $h_\sigma$
for all sufficiently large $j$.

The smooth comparison requirements also follow
uniformly in $\sigma$. Indeed, constant rescaling
leaves the Levi-Civita connection unchanged, and
the order-$m$ norms acquire only the factor
$\sigma^m$. These factors are uniformly bounded
for $\sigma\in I$ and $0\le m\le N$.
Since $K^+\subset Q$ and
$\sigma^4\tau\in[-16,0]$, smooth spacetime
convergence gives \eqref{eq:scrPmetric} and \eqref{eq:scrPJ} for $h_\sigma$.

Finally, for $z\in K$, the distance bounds and the
Lipschitz property of $\rho_{h_\sigma}$ imply
$
 \overline B_{h_\sigma}(z,\eta/100)
 \subset\operatorname{int}K^+.
$
On this ball, by the metric comparison
$
 \Phi^*G(\tau)\ge(1-\delta)h_\sigma.
$
A path starting at $\Phi(z)$ and leaving the image
of this ball must have length at least
$\sqrt{1-\delta}\,\eta/100>\eta/500$
before its first exit. Therefore
we obtain the capture condition
\eqref{eq:scrcapture} for $h_\sigma$.

Every requirement of test \textup{(P)} thus holds
for all sufficiently large $j$, uniformly for
$\sigma\in I$. 

\smallskip
\noindent\emph{Improvement of \textup{(R)}.}
At the selected scales $s_j\to s\in I$, write
$G_j(\tau):=s_j^{-2}g_j(s_j^4\tau)$
and $h_s:=s^{-2}h$.
We regard $G_j$ on its full available backward
interval $[-T_j/s_j^4,0]$. By the previously obtained spacetime convergence 
$
 (M_j,G_j(\tau),J_j,p_j)
 \longrightarrow (X,h_s,J_X,p),
$
where the limiting metric is independent of $\tau$.
By the scaling property of the canonical radius, for every fixed
$L>0$ and all sufficiently large $j$, we have on
$M_j\times[-L,0]$
\begin{align}\label{text R 14}
 \scr_{G_j}(q,\tau)
 =s_j^{-1}\scr_j(q,s_j^4\tau)
 \ge\frac1{2s_j}\ge\frac14.
\end{align}
If \textup{(R)} fails at the selected scales,
its normalized energy hypothesis holds:
\[
 E_j:=
 \int_{B_{G_j(0)}(p_j,2)}
 |\Rm_{G_j(0)}|^2\,dV_{G_j(0)}
 <\eps_R.
\]
Fix $b:=1/1024$.
For every $q\in B_{G_j(0)}(p_j,15/8)$, we get
$
 B_{G_j(0)}(q,2b)
 \subset B_{G_j(0)}(p_j,2).
$ Thus, by \eqref{text R 14}, we apply \textup{(R)} at scale $b$ to obtain
\begin{align}\label{text R Rm m}
 \sup_{\substack{B_{G_j(0)}(q,b)
                \times [-\theta b^4,0]}}
 |\nabla^m\Rm_{G_j(\tau)}|(y)
 \le C_mb^{-2-m},
 \qquad 0\le m\le N.
\end{align}

\begin{claim}
$
 B_{h_s}(p,15/8)\subset\mathcal R(X).
$
\end{claim}
\begin{proof}
Fix $z\in B_{h_s}(p,15/8)$ and choose points $q_j\to z$. Then
$
 d_{G_j(0)}(p_j,q_j)
 \longrightarrow d_{h_s}(p,z)<15/8.
$
Thus $q_j\in B_{G_j(0)}(p_j,15/8)$ for all
sufficiently large $j$, and the curvature bounds
\eqref{text R Rm m} apply around $q_j$.

The lower-order Sobolev coefficients are
$s_j^2\beta_j\to0$. The noncollapse argument in
Step~1 of the proof of Theorem~\ref{thm:selected}
therefore gives uniform volume lower bounds for
$B_{G_j(0)}(q_j,r)$, $0<r\le b$.
Together with \eqref{text R Rm m}, these give
uniform local injectivity-radius bounds at $\tau=0$.
The chart construction provides
time-independent holomorphic embeddings
$\Psi_j:B_2\subset\C^2\to M_j$, with $\Psi_j(0)=q_j$,
such that the metrics $\Psi_j^*G_j(\tau)$ are
uniformly elliptic and have uniform spatial
$C^8$ bounds for every $\tau\in[-\theta b^4,0]$.
The first-exit argument
used for \eqref{eq:scrcapture} consequently gives
$
 B_{G_j(0)}(q_j,\rho)\subset\Psi_j(B_1)
$
for a fixed $\rho>0$. The ellipticity bounds, spatial $C^8$ bounds,
and $\rho$ are independent of $j$ and $q_j$.

Applying \eqref{eq:interiorsmoothing} gives uniform
bounds for every derivative of $\Psi_j^*G_j(\tau)$ on $B_1\times[-\theta b^4/2,0]$.
After passing to a local subsequence,
$
 \Psi_j^*G_j(0)\rightarrow k_\infty$
in $C^\infty_{\mathrm{loc}}(B_1),
$
where $k_\infty$ is smooth and positive definite.

On $B_{G_j(0)}(q_j,\rho/4)$, minimizing
curves remain inside the chart by the same
first-exit argument. Hence the smooth local limit
also describes the ambient metric-space limit
near the centers. By uniqueness of the GH limit,
it identifies a neighborhood of $z$ in $X$ as
smooth. Thus $z\in\mathcal R(X)$ and the claim follows.
\end{proof}
Choose a smaller ball $\overline B_{h_s}(p,7/4)\subset\mathcal R(X)$.
Smooth convergence therefore yields
\[
 \int_{B_{h_s}(p,7/4)}
 |\Rm_{h_s}|^2\,dV_{h_s}
 \le\liminf_{j\to\infty}E_j
 \le\eps_R<\eps_{\rm stat}.
\]
Apply the stationary estimate
\eqref{eq:modelinterior} to $h_s$ with $r=1$:
\[
 \sup_{B_{h_s}(p,5/4)}
 |\nabla^m\Rm_{h_s}|
 \le H_m,
 \qquad 0\le m\le N.
\]
Choose a fixed regular domain $U$ satisfying
$
 \overline B_{h_s}(p,9/8)
 \subset U\Subset B_{h_s}(p,5/4).
$
The coverage conclusion of
Theorem~\ref{thm:selected} ensures that, for large $j$,
$
 B_{G_j(0)}(p_j,1)\subset\Phi_j(U).
$
Since $\Phi_j$ is independent of time, smooth spacetime convergence implies
\[
 \limsup_{j\to\infty}
 \sup_{\substack{ B_{G_j(0)}(p_j,1)\times
                 [-\theta,0]}}
 |\nabla^m\Rm_{G_j(\tau)}|(q)
 \le H_m<C_m.
\]
Because $C_m=4H_m+1$ and there are only finitely
many orders $0\le m\le N$, all the strict bounds
required by \textup{(R)} hold simultaneously
for large $j$. This contradicts its failure.

\smallskip
\noindent\emph{Improvement of \textup{(V)}.}
The ball-volume convergence in
Theorem~\ref{thm:selected} passes to the same
normalized flows. Thus
$
 \Vol_{G_j(0)}B_{G_j(0)}(p_j,1)
 \longrightarrow
 \Vol_{h_s}B_{h_s}(p,1).
$
Since $(X,h_s,J_X,p)\in\mathcal A$, the model
volume bound gives
$
 \Vol_{h_s}B_{h_s}(p,1)\le V_*<B_0.
$
Hence
\[
 \Vol_{G_j(0)}B_{G_j(0)}(p_j,1)<B_0
\]
for all sufficiently large $j$. This is precisely
test \textup{(V)}, contradicting its failure.

Every possible failing test at the selected scales has been excluded.
Thus there cannot be infinitely many indices with
$\scr_j(p_j,0)<2$. It follows that
$\scr_j(p_j,0)\ge2$ for all sufficiently large $j$,
which also rules out the normalization
$\scr_j(p_j,0)=1$.
\end{proof}

\subsection{A positive lower bound for the canonical radius}

We use weighted point selection as in
\cite[Section~3.4, Proposition~3.43]{ChenWangB}. An approximate infimum suffices.

\begin{theorem}[Lower bound for the canonical radius]
\label{thm:radiuslower}
Let $g(t)$, $t\in[0,\infty)$, be a global Calabi flow on a compact K\"ahler surface in a fixed
real K\"ahler class, and suppose it satisfies the uniform Sobolev bound. For each
$t_0>0$ there is $r_0>0$ such that
\[
 \inf_{M\times[t_0,\infty)}\scr\ge r_0.
\]
The constant $r_0$ may depend on the fixed flow.
\end{theorem}
\begin{proof}
Choose $a\in(0,t_0)$.
Suppose that the conclusion fails. Then there are
$x_j\in M$ and $b_j\ge t_0$ such that
\[
 \varepsilon_j:=\scr(x_j,b_j)\longrightarrow0.
\]
By Lemma~\ref{lem:radiusproperties}, $\scr$ has a
positive lower bound on every compact time slab
contained in $M\times(0,\infty)$. Hence $b_j\to\infty$:
otherwise a subsequence of the points $(x_j,b_j)$
would lie in one such slab, contradicting
$\varepsilon_j\to0$.

For $(x,t)\in M\times(a,b_j]$, define
\[
 m_j:=\inf_{M\times(a,b_j]}F_j,\qquad
 F_j(x,t):=\frac{\scr(x,t)}{(t-a)^{1/4}}.
\] For each fixed $j$,
Lemma~\ref{lem:radiusproperties} gives $c_j>0$ such that
$\scr\ge c_j$ on $M\times[a,b_j]$.
Since $t-a\le b_j-a$ throughout the domain of $F_j$,
\[
 m_j\ge\frac{c_j}{(b_j-a)^{1/4}}>0.
\]
On the other hand,
\[
 m_j\le F_j(x_j,b_j)
 =\frac{\varepsilon_j}{(b_j-a)^{1/4}}
 \le\frac{\varepsilon_j}{(t_0-a)^{1/4}}
 \longrightarrow0.
\]
In particular, $0<m_j<\infty$.
That enables us to
choose $(p_j,t_j)\in M\times(a,b_j]$ such that
$
 F_j(p_j,t_j)\le\frac32m_j.
$
This is an approximate minimizer of $F_j$; its existence
requires neither continuity of $\scr$ nor attainment
of the infimum.
Since $t_j\le b_j$, setting 
$
 r_j:=\scr(p_j,t_j)
$ we have
\[
 \begin{aligned}
 r_j
 &=F_j(p_j,t_j)(t_j-a)^{1/4}
 \le\frac32F_j(x_j,b_j)(t_j-a)^{1/4}
 =\frac32\varepsilon_j
   \left(\frac{t_j-a}{b_j-a}\right)^{1/4}
 \le\frac32\varepsilon_j\longrightarrow0.
 \end{aligned}
\]
We further have $t_j\to\infty$.
Otherwise, if $t_j$ remained bounded along a subsequence,
the selected points would lie in a fixed slab
$M\times[a,T]$, on which $\scr$ has a positive lower
bound, contradicting $r_j\to0$.

For every $q\in M$ and
$u\in[(t_j+a)/2,t_j]$, we use $m_j\ge\frac23F_j(p_j,t_j)$ to get
\begin{equation}\label{eq:weightedselection}
 \begin{aligned}
 \scr(q,u)
 &\ge m_j(u-a)^{1/4}
 \ge 2^{-1/4}m_j(t_j-a)^{1/4}\\
 &\ge\frac23\,2^{-1/4}
       F_j(p_j,t_j)(t_j-a)^{1/4}
 =\frac23\,2^{-1/4}r_j
 >\frac12r_j.
 \end{aligned}
\end{equation}
Define
\[
 L_j:=\frac{t_j-a}{2r_j^4}
     =\frac{1}{2F_j(p_j,t_j)^4}.
\]
Since $F_j(p_j,t_j)\le3m_j/2\to0$, it follows that
$L_j\to\infty$.
Now consider the Calabi rescalings
\[
 \widehat g_j(\tau):=
 r_j^{-2}g(t_j+r_j^4\tau),
 \qquad -T_j\le\tau\le0,
 \qquad T_j:=\frac{t_j}{r_j^4}.
\]
The scaling property in
Lemma~\ref{lem:radiusproperties} gives
\begin{align}\label{tilde scr 0}
 \widehat{\scr}_j(q,\tau)
 =r_j^{-1}\scr(q,t_j+r_j^4\tau),
 \qquad
 \widehat{\scr}_j(p_j,0)=1.
\end{align}
If $-L_j\le\tau\le0$, then
$
 \frac{t_j+a}{2}
 \le t_j+r_j^4\tau\le t_j.
$
Therefore \eqref{eq:weightedselection} implies
\begin{align}\label{tilde scr 12}
 \widehat{\scr}_j(q,\tau)\ge\frac12
 \qquad
 \text{for every }q\in M,\quad -L_j\le\tau\le0.
\end{align}
Moreover,
$
 T_j-L_j=\frac{t_j+a}{2r_j^4}\longrightarrow\infty,
$
so for every fixed $L>0$, the radius
lower bound holds on all of $M\times[-L,0]$ for
all sufficiently large $j$.

We verify the remaining hypotheses of
Theorem~\ref{thm:selected}. By Lemma~\ref{lem:energy}
and the scale invariance of curvature energy in real
dimension four, we have the same energy bound for $\widehat g_j(\tau)$.
Rescaling \eqref{eq:Sob} gives
\[
 \begin{gathered}
 \|v\|_{L^4(\widehat g_j(\tau))}^2
 \le A\left(
 \int_M|\nabla v|_{\widehat g_j(\tau)}^2\,
 dV_{\widehat g_j(\tau)}
 +\beta_j\int_Mv^2\,dV_{\widehat g_j(\tau)}
 \right),\qquad
 \beta_j:=r_j^2V^{-1/2}\longrightarrow0,
 \end{gathered}
\]
for every smooth function $v$.
For large $j$, $L\le L_j$, and consequently
$
 t_j-Lr_j^4
 \ge\frac{t_j+a}{2}\longrightarrow\infty.
$
Let $\widehat D_j$ denote the dissipation rate of
$\widehat g_j$. By \eqref{eq:scaleddissipation} and
\eqref{eq:totaldissipation}, for every fixed $L>0$,
\[
 \begin{aligned}
 \int_{-L}^0\widehat D_j(\tau)\,d\tau
 &=\int_{t_j-Lr_j^4}^{t_j}D(u)\,du
 =\frac12\left(
 \Ca(t_j-Lr_j^4)-\Ca(t_j)\right)\longrightarrow0.
 \end{aligned}
\]
Since $\Ca$ is nonnegative and nonincreasing, it has
a finite limit.
Finally, the complex structure remains $J$. With
$\lambda_j:=r_j^{-2}\to\infty$, the rescaled
K\"ahler forms satisfy
\[
 \left[\lambda_j^{-1}
       \omega_{\widehat g_j(\tau)}\right]
 =[\omega(t_j+r_j^4\tau)]
 =\Omega.
\]
Thus the fixed-class hypothesis of
Theorem~\ref{thm:selected} holds with the identity
diffeomorphisms of $M$.

The rescaled pointed flows therefore satisfy all
hypotheses of Theorem~\ref{thm:selected}. They also
satisfy
\eqref{tilde scr 0} and \eqref{tilde scr 12}.
Lemma~\ref{lem:radiusimprovement} rules out precisely
such a sequence. This contradiction proves the theorem.
\end{proof}

\subsection{Consequences for the original flow}

The lower bound just proved makes the defining backward
estimate available uniformly along the fixed global flow.
Interior smoothing gives every prescribed derivative order.

\begin{lemma}[Backward estimate from small curvature energy]
\label{lem:globalbackwardenergy}
Let $g(t)$, $t\in[0,\infty)$, be a global Calabi flow
satisfying the hypotheses of Theorem~\ref{thm:radiuslower}.
For every $t_0>0$, there exist constants
\[
 0<r_*<\min\{V^{1/4},(t_0/2)^{1/4}\},
 \qquad \eps>0,\qquad 0<\theta\le1,
\]
and, for every integer $k\ge0$, a constant $C_k<\infty$
with the following property. For every $x\in M$,
$t\ge t_0$, and $0<r<r_*$,
\begin{equation}\label{eq:Rtrigger}
 \begin{gathered}
 \int_{B_{g(t)}(x,2r)}
 |\Rm_{g(t)}|_{g(t)}^2\,dV_{g(t)}<\eps\quad
 \Longrightarrow
 \sup_{\substack{y\in B_{g(t)}(x,r)\\
                 t-\theta r^4\le u\le t}}
 r^{2+k}
 |\nabla^k\Rm_{g(u)}|_{g(u)}(y)
 \le C_k.
 \end{gathered}
\end{equation}
\end{lemma}
\begin{proof}
By Theorem~\ref{thm:radiuslower}, there is
$r_0>0$ such that
$
 \scr(y,u)\ge r_0$
for all $y\in M$, $u\ge t_0.$
Write $\theta_R$ and $C_m^{\rm R}$, $0\le m\le N$,
for the constants in test \textup{(R)}.
Set $\eps:=\eps_R$, and choose
\[
 0<r_*<
 \min\{r_0,V^{1/4},(t_0/2)^{1/4}\}.
\]
Fix $t\ge t_0$, $x\in M$, and $0<r<r_*$. Then
$t-r^4>t-t_0/2>0$. Suppose that
\[
 \int_{B_{g(t)}(x,2r)}
 |\Rm_{g(t)}|^2\,dV_{g(t)}<\eps_R.
\]
Let $q\in B_{g(t)}(x,r)$ and put $s:=r/4$. Then
\[
 B_{g(t)}(q,2s)
 \subset B_{g(t)}(x,3r/2)
 \subset B_{g(t)}(x,2r).
\]
Thus the energy hypothesis of test \textup{(R)} holds
at $(q,t,s)$. 
Introduce the normalized flow
and the fixed spatial domain
\[k_q(u):=s^{-2}g(t+s^4u),
 -\theta_R\le u\le0,\qquad
 \Omega_q:=B_{k_q(0)}(q,1)
          =B_{g(t)}(q,s).
\]
Test \textup{(R)} gives
\begin{align}\label{higher
curvature-derivative bounds}
 \sup_{\Omega_q\times[-\theta_R,0]}
 |\nabla^m\Rm_{k_q(u)}|_{k_q(u)}
 \le C_m^{\rm R},
 \qquad 0\le m\le N=12.
\end{align}

By Lemma~\ref{lem:noncollapse},
the normalized terminal metric satisfies noncollapsing.
Together with the terminal curvature bounds, this gives
a uniform positive lower bound for the injectivity
radius at $q$. The curvature-derivative bounds therefore
provide controlled terminal coordinate charts on a
uniformly sized neighborhood of $q$.
These charts can be kept fixed throughout the backward
interval by using the Calabi flow equation and \eqref{higher
curvature-derivative bounds}.
The connection variation formula and the higher
curvature-derivative bounds similarly control 
uniform spatial $C^8$ bounds on a smaller fixed
coordinate neighborhood.

Using the fixed complex structure and the local
holomorphic-chart construction in Step~1 of
Theorem~\ref{thm:selected}, we obtain time-independent
holomorphic charts
$\Psi_q:B_2\subset\C^2\longrightarrow M,
\Psi_q(0)=q,
$
whose images lie in $\Omega_q$, and constants
$\lambda>0$ and $C<\infty$ such that
\[
 \lambda e\le\Psi_q^*k_q(u)\le\lambda^{-1}e,
 \qquad
 \sup_{u\in[-\theta_R,0]}
 \|\Psi_q^*k_q(u)\|_{C^8(B_2)}\le C.
\]

Apply the interior smoothing estimate
\eqref{eq:interiorsmoothing} in these fixed holomorphic
charts, after a fixed Calabi rescaling of the time
interval. On smaller cylinders, it gives every spatial derivative
of the metric. Uniform ellipticity then converts these
coordinate estimates into covariant curvature estimates.
Consequently, for every integer $k\ge0$, there is
$A_k<\infty$, independent of $q,t,r$, such that
\[
 \sup_{v\in[t-\theta_Rs^4/2,t]}
s^{2+k} |\nabla^k\Rm_{g(v)}|_{g(v)}(q)=
 \sup_{u\in[-\theta_R/2,0]}
 |\nabla^k\Rm_{k_q(u)}|_{k_q(u)}(q)
 \le A_k.
\]
Since $s=r/4$, set
$
 \theta:=\frac{\theta_R}{2\cdot4^4},
 C_k:=4^{2+k}A_k.
$
The preceding estimate becomes
\[
 \sup_{v\in[t-\theta r^4,t]}
 r^{2+k}|\nabla^k\Rm_{g(v)}|_{g(v)}(q)
 \le C_k.
\]
The point $q$ was arbitrary in $B_{g(t)}(x,r)$, and
the constants are uniform in $q$. Taking the supremum
over this fixed terminal ball proves
\eqref{eq:Rtrigger}. The $C_k$ depend on
the fixed flow, $t_0$, and $k$.
\end{proof}

\begin{corollary}[Geometric properties of the time translates]
\label{cor:geometricinputs}
Under Theorem~\ref{thm:radiuslower}, let $t_j\to\infty$ and
write $g_j(\tau):=g(t_j+\tau)$. The following short-time properties hold for these translates,
including at arbitrary moving times in bounded intervals:
\begin{enumerate}[label=\textup{(B\arabic*)}]
\item Subsequential slice limits are completed compact
K\"ahler orbifolds. Convergence is global in the GH sense and
smooth on the entire geometric regular locus. For each
$K\Subset U\Subset U^+\Subset \mathcal R(X)$, comparison
embeddings are defined on $U^+$, and their images of $U$
contain all points of the underlying manifold corresponding to $K$ for large $j$.
The natural volume measures converge, with total limiting volume $V$.
\item For fixed terminal domains
$K\Subset U\Subset U^+\Subset \mathcal R(X)$, terminal
embeddings on $U^+$ give smooth local Calabi-flow limits on
$U\times[-\sigma,0]$ for some $\sigma>0$, after translating
the terminal time to zero. The terminal metric separations
and $\sigma$ have positive lower bounds independent of $j$.
The same choices work for isometric copies of these nested
domains in terminal limits.
\item If $a_j-b_j\to0$ in a bounded translated interval,
then for any sequences of points $x_j,y_j$ of the underlying manifold,
held fixed as time varies,
\[
 \left|d_{g_j(a_j)}(x_j,y_j)
       -d_{g_j(b_j)}(x_j,y_j)\right|\longrightarrow0.
\]
\end{enumerate}
Thus these short-time properties follow from the flow and Sobolev
assumptions. Lemma~\ref{lem:shortdrop} gives the corresponding statement
for finite flow segments with vanishing energy drop.
\end{corollary}
\begin{proof}
Fix a bounded translated interval. For large translates,
the positive canonical-radius bound is available throughout
this interval and a larger past collar. The fixed volume,
diameter bound, noncollapse, curvature-energy bound, and
vanishing translated dissipation come from
Section~\ref{sec:preliminaries}; (V) supplies the local upper
volume bound.

For \textup{(B1)}, apply Proposition~\ref{prop:slice}
globally on the compact sources, using
Lemma~\ref{lem:globalbackwardenergy} for the backward estimate. Here the Sobolev lower-order coefficient is $V^{-1/2}$.
The local upper volume bound supplied by (V), together with the
local absorption statement in
Lemma~\ref{lem:sobolev-passage}, gives the homogeneous inequality
with constant $2A$ on sufficiently small balls. Since
$2A<A_*$, the branchwise removal and both model gaps apply with
the constants already chosen. As in Steps 1--2 of
Theorem~\ref{thm:selected}, they exclude branching and
concentration over geometrically regular points.
This gives compact extremal K\"ahler orbifold
limits, smooth convergence on the entire regular locus,
and coverage. Lemma~\ref{lem:volume-passage} gives natural volume convergence;
since the source manifolds have total volume $V$, the completed
limit also has total volume $V$.

For \textup{(B2)}, fix
$K\Subset U\Subset U^+\Subset \mathcal R(X)$ in the terminal
limit $(X,h)$. Choose a common scale $a>0$ using the curvature
and injectivity bounds near $\overline{U^+}$, the positive
canonical radius, and
$4a<d_h(\overline U,X\setminus U^+)$.
Cover $\overline U$ by finitely many terminal balls of radius
$a/2$ whose concentric $2a$-balls lie in $U^+$.
Choose $a$ small enough that terminal smooth convergence makes
the larger-ball energies strictly less than $\eps_R$.
Apply (R) on the fixed terminal $a$-balls for an interval of
length $\theta a^4$, and apply the interior estimate~\eqref{eq:interiorsmoothing} on
their $a/2$-balls during the latter half of that interval.
The scale $a$ and the resulting lifetime depend only on these
intrinsic terminal data and the positive canonical radius.
They are independent of $j$ and unchanged when the nested
terminal domains are carried by isometries.

For \textup{(B3)}, repeat the fixed-core argument
\eqref{eq:selectedflux}--\eqref{eq:selectedIDC} at
arbitrary moving terminal times. The regular exterior
and finitely many small cores now exhaust the compact
source, so no outer boundary is needed. The terminal core
volumes are $O(\rho^4)$ by (V), and their collars are
controlled by \textup{(B2)}. Bounded boundary flux
preserves this bound when the time difference tends to
zero. Noncollapse gives the intrinsic diameter bound
\eqref{eq:selectedcorediam}, and comparison of the
quotient length metrics gives the claimed assertion
for the same points of the underlying manifold. First let $j\to\infty$
with $\rho$ fixed, then let $\rho\downarrow0$.
\end{proof}
\begin{corollary}[Terminal-ball backward pseudolocality]
\label{cor:backwardpseudolocality}
Let $g(t)$, $t\in[0,\infty)$, be a global Calabi flow
satisfying the hypotheses of Theorem~\ref{thm:radiuslower}.
Fix $t_0>0$, and let $N=12$ and $C_m$, $0\le m\le N$,
be as in test \textup{(R)}.
There exist constants
\[
 0<a<1/2,\qquad 0<\theta\le1,\qquad
 0<r_0<t_0^{1/4},
\]
such that the following holds. Suppose $x\in M$,
$t\ge t_0$, and $0<r<r_0/(2a)$ satisfy
\[
 |\Rm_{g(t)}|_{g(t)}\le r^{-2}
 \qquad\text{on }B_{g(t)}(x,2r).
\]
Then
\[
 \sup_{B_{g(t)}(x,r)\times[t-\theta a^4r^4,t]}
 |\nabla^m\Rm_{g(u)}|_{g(u)}(y)
 \le C_m(ar)^{-2-m},
 \qquad 0\le m\le N.
\]
Moreover, for every integer $k\ge0$, there is a
constant $D_k<\infty$ such that
\[
 \sup_{B_{g(t)}(x,r/2)\times
                 [t-\theta a^4r^4/2, t]}
 |\nabla^k\Rm_{g(u)}|_{g(u)}
 \le D_k(ar)^{-2-k}.
\]
The constants $a,\theta,r_0$ may depend on the fixed
flow and $t_0$.
\end{corollary}
\begin{proof}
Fix $t_0>0$. By Theorem~\ref{thm:radiuslower}, there is
$r_{\rm can}>0$ such that
$
 \scr(q,t)\ge r_{\rm can}$
for every $q\in M$, $t\ge t_0$.
Write $\theta_R$ for the time constant in test
\textup{(R)}. Let
$
 r_*,\eps_{\rm b},\theta_{\rm b},\{D_k\}_{k\ge0}
$
be the constants supplied by
Lemma~\ref{lem:globalbackwardenergy} for this flow
and this value of $t_0$.
Choose $a$ sufficiently small that
\[
 0<a<1/2,\qquad
 16B_0a^4<\min\{\eps_R,\eps_{\rm b}\}.
\]
Next choose
\[
 0<r_0<
 \min\{r_{\rm can},2r_*,t_0^{1/4}\},
 \qquad
 \theta:=\min\{\theta_R,2\theta_{\rm b}\}.
\]
Now put
$s:=ar$. Since $r<r_0/(2a)$, we get
$
 2s<r_0<r_{\rm can}$ and $s<r_*$.
The testing intervals are contained in the domain
of the flow, because
$t-(2s)^4>t_0-r_0^4>0.$

Let $q\in B_{g(t)}(x,r)$. Since $2a<1$, the triangle
inequality gives
\[
 B_{g(t)}(q,2s)
 \subset B_{g(t)}(x,(1+2a)r)
 \subset B_{g(t)}(x,2r).
\]
Thus the assumed terminal curvature bound holds
throughout $B_{g(t)}(q,2s)$.

The canonical-radius lower bound makes test
\textup{(V)} available at $(q,t)$ at scale $2s$.
Rescaling the volume bound in that test gives
$
 \Vol_{g(t)}B_{g(t)}(q,2s)<B_0(2s)^4.
$
Consequently,
\[
 \begin{aligned}
 \int_{B_{g(t)}(q,2s)}
 |\Rm_{g(t)}|_{g(t)}^2\,dV_{g(t)}
 &\le r^{-4}\Vol_{g(t)}B_{g(t)}(q,2s)\\
 &<r^{-4}B_0(2ar)^4
 =16B_0a^4
 <\min\{\eps_R,\eps_{\rm b}\}.
 \end{aligned}
\]
This estimate holds uniformly for every
$q\in B_{g(t)}(x,r)$.

Since the energy is less than $\eps_R$, apply
test \textup{(R)} at $(q,t)$ at scale $s$.
In particular, at the center $q$ it yields
\[
 \sup_{u\in[t-\theta_Rs^4,t]}
 |\nabla^m\Rm_{g(u)}|_{g(u)}(q)
 \le C_ms^{-2-m},
 \qquad 0\le m\le N.
\]
Because $\theta\le\theta_R$, this estimate holds
on $[t-\theta s^4,t]$.
Taking the supremum over
$q\in B_{g(t)}(x,r)$ and using $s=ar$ proves
the first conclusion.

For the higher derivatives, the same energy estimate
is less than $\eps_{\rm b}$, and $s<r_*$.
Lemma~\ref{lem:globalbackwardenergy}, applied at
$(q,t)$ with radius $s$, therefore gives
\[
 \sup_{u\in[t-\theta_{\rm b}s^4,t]}
 |\nabla^k\Rm_{g(u)}|_{g(u)}(q)
 \le D_ks^{-2-k}
 \qquad\text{for every }k\ge0.
\]
Our choice of $\theta$ ensures
$\theta/2\le\theta_{\rm b}$, so
$
 [t-\theta s^4/2,t]
 \subset[t-\theta_{\rm b}s^4,t].
$
Restricting to $q\in B_{g(t)}(x,r/2)$ and taking
the supremum proves the second conclusion.
All constants are uniform in $q,x,t,r$, as required.
\end{proof}

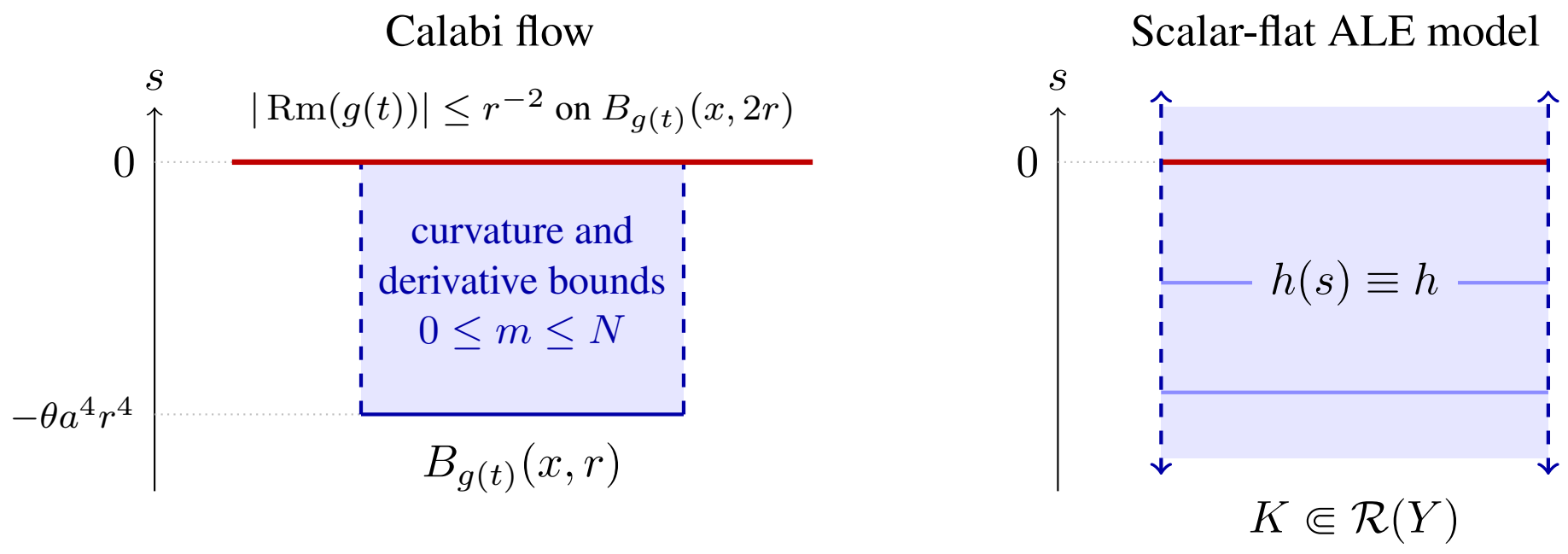
\begin{figure}[htbp]
\centering
\begin{tikzpicture}[x=1cm,y=0.85cm,font=\small]
 \node at (3.5,3.8) {Calabi flow};
 \draw[->] (0.9,-0.4) -- (0.9,3.1) node[above] {$s$};
 \draw[black!25,densely dotted] (0.9,2.6) -- (6.0,2.6);
 \draw[black!25,densely dotted] (0.9,0.3) -- (2.5,0.3);
 \node[left,inner sep=4pt] at (0.9,2.6) {$0$};
 \node[left,inner sep=4pt,font=\scriptsize] at (0.9,0.3)
  {$-\theta a^4r^4$};
 \fill[blue!10] (2.5,0.3) rectangle (5.0,2.6);
 \draw[blue!65!black,dashed,thick] (2.5,0.3) -- (2.5,2.6);
 \draw[blue!65!black,dashed,thick] (5.0,0.3) -- (5.0,2.6);
 \draw[blue!65!black,thick] (2.5,0.3) -- (5.0,0.3);
 \draw[red!75!black,very thick] (1.5,2.6) -- (6.0,2.6);
 \node[above,inner sep=6pt,font=\scriptsize] at (3.75,2.6)
  {$|\Rm(g(t))|\le r^{-2}$ on $B_{g(t)}(x,2r)$};
 \node[blue!65!black,align=center,font=\footnotesize] at (3.75,1.5)
  {curvature and\\derivative bounds\\$0\le m\le N$};
 \node[below,inner sep=6pt] at (3.75,0.3) {$B_{g(t)}(x,r)$};

 \begin{scope}[xshift=7.7cm]
  \node at (2.35,3.8) {Scalar-flat ALE model};
  \draw[->] (0.2,-0.4) -- (0.2,3.1) node[above] {$s$};
  \node[left,inner sep=4pt] at (0.2,2.6) {$0$};
  \draw[black!25,densely dotted] (0.2,2.6) -- (4.0,2.6);
  \fill[blue!10] (1.0,-0.1) rectangle (4.0,3.1);
  \foreach \y in {0.5,1.5} {
   \draw[blue!45,thick] (1.0,\y) -- (4.0,\y);
  }
  \draw[red!75!black,very thick] (1.0,2.6) -- (4.0,2.6);
  \draw[blue!65!black,dashed,thick,<->] (1.0,-0.25) -- (1.0,3.25);
  \draw[blue!65!black,dashed,thick,<->] (4.0,-0.25) -- (4.0,3.25);
  \node[fill=blue!10,inner sep=4pt] at (2.5,1.5) {$h(s)\equiv h$};
  \node[below,inner sep=5pt] at (2.5,-0.25) {$K\Subset\mathcal R(Y)$};
 \end{scope}
\end{tikzpicture}
\caption{The red lines mark the common zero-time slice in the two panels.}
\label{fig:backwardcomparison}
\end{figure}

Figure~\ref{fig:backwardcomparison} compares this corollary with
the stationary model. Under the hypotheses of the corollary,
curvature control on the larger terminal ball gives backward
estimates on the smaller fixed ball. For a scalar-flat ALE model,
$h(\tau)\equiv h$, so bounds on a compact regular set
$K\Subset\mathcal R(Y)$ persist for all time.
In both panels we use the shifted time $s:=\tau-t$, so the terminal
slice is at $s=0$; the spatial sets are fixed as time varies.

%%%%%%%%%%%%%%%%%%%%%%%%%%%%%%%%%%%%%%%%%%%%%%%%%%%%%%%%%%%%%%%%%%%%%%%%%%%%%%%%%%%%%%%%%%%%%%%%%%%%%%%%%%%%%%%%%%%%%%%%%%%%%%%%%%%%%%%%%%%%%%%%%%%%%%%%%%%%%%%%%%%%%%%%%%%%%%%%%%%%%%%%%%%%%%%%%%%%%%%%%%%%%%%%%%%%%%%%%%%%%%

\section{Long-time two-sided pseudo-locality}
\label{sec:upgrade}

We prove a curvature estimate under a small drop of the Calabi energy,
before extracting a limit flow. Following Chen--Wang
\cite{ChenWangB}, we use slice limits and short backward
cylinders to extend regularity in both time directions.

Throughout this section, fix a compact complex surface $(M,J)$,
a K\"ahler class $\Omega$, and bounds $A,E_0$ for the Sobolev inequality
and curvature energy. The volume $V$ is determined by $\Omega$.
All canonical radii use the model class and defining constants of
Section~\ref{sec:scr}. For the dependence statements below, these
choices are fixed once as functions of $(M,J,\Omega,A,E_0)$.
Theorem~\ref{thm:radiuslower} supplies the
required lower radius bound along the original global flow.

\subsection{Analyticity and exterior extension}

We write the extremal equation in a local K\"ahler potential and
apply Morrey's analytic elliptic regularity theorem \cite{Morrey}.

\begin{lemma}[Analyticity of extremal metrics]\label{lem:analytic}
A smooth extremal K\"ahler orbifold metric is real analytic in holomorphic uniformizing charts.
\end{lemma}
\begin{proof}
We choose a smooth real potential $\varphi$ and write $g_{i\bar j}=\varphi_{i\bar j}$. 
Since $g$ is extremal, \(Z=\nabla^{1,0}S\) is holomorphic. So \(\partial_{\bar j}Z^i=0\) and letting $f:=S-Z\varphi$, we have
\[
\begin{aligned}
\partial_{\bar j} f
=S_{\bar j}
  -(\partial_{\bar j}Z^i)\varphi_i
  -Z^i\varphi_{i\bar j}
=S_{\bar j}-Z^i g_{i\bar j}
=0,
\end{aligned}
\]
i.e. \(f\) is holomorphic.
We regard the holomorphic functions $Z,f$ as analytic data. The real potential solves
\[
-g^{i\bar j}\partial_i\partial_{\bar j}
\log\det(\varphi_{k\bar l})=\operatorname{Re}(Z\varphi+f).
\]
This operator is analytic where the Hermitian Hessian is positive definite. Its linearization has principal symbol $-c|\xi|_g^4$, for a fixed $c>0$ depending only on conventions. Thus it is locally uniformly elliptic. Nonlinear analytic interior regularity \cite{Morrey} gives analyticity of $\varphi$ and $g$. We apply the same argument on each smooth orbifold cover.
\end{proof}

Lemma~\ref{lem:endextension} uses nullness and Hartogs extension on
ALE ends. The compact case below instead uses analytic continuation
on complete frame manifolds. The regular exterior carries the
fundamental group, and equality of total volumes rules out a
nontrivial covering.

\begin{figure}[htbp]
\centering
\begin{tikzpicture}[x=1cm,y=1cm,font=\small]
 \node at (3.2,4.55) {The exterior $U$};
 \node at (10.7,4.55) {A regular collar near $p_a$};
 \path[fill=blue!9,draw=black!50]
 (0.35,3.3) .. controls (0.7,4.15) and (2.1,4.25) .. (3.15,3.95)
 .. controls (4.5,4.35) and (5.85,3.9) .. (6.15,3.05)
 .. controls (6.5,1.95) and (6.0,0.65) .. (5.1,0.4)
 .. controls (4.3,0.1) and (3.65,0.35) .. (3.2,0.4)
 .. controls (1.9,0.0) and (0.6,0.45) .. (0.3,1.25)
 .. controls (0.05,1.9) and (0.05,2.65) .. (0.35,3.3) -- cycle;
 \node at (0.5,3.85) {$X$};
 \node at (3.25,2.8) {$U$};
 \foreach \x/\y/\a in {1.65/2.15/1,4.65/1.55/2} {
   \fill[white] (\x,\y) circle (0.62);
   \fill[black!10] (\x,\y) circle (0.31);
   \draw[dashed,black!45] (\x,\y) circle (0.31);
   \path[pattern=north east lines,pattern color=orange!70!black,even odd rule]
     (\x,\y) circle (0.775) (\x,\y) circle (0.465);
   \draw (\x,\y) circle (0.62);
   \fill (\x,\y) circle (0.025);
   \node[above right,inner sep=2pt] at (\x,\y) {$p_{\a}$};
 }
 \begin{scope}[shift={(10.55,2.15)}]
   \fill[blue!9] (-2.05,-1.95) rectangle (2.45,1.95);
   \fill[white] (0,0) circle (1.4);
   \fill[black!10] (0,0) circle (0.7);
   \draw[dashed,black!45] (0,0) circle (0.7);
   \path[pattern=north east lines,pattern color=orange!70!black,even odd rule]
     (0,0) circle (1.75) (0,0) circle (1.05);
   \draw (0,0) circle (1.4);
   \fill (0,0) circle (0.025);
   \node[above left,inner sep=3pt] at (0,0) {$p_a$};
   \draw[->,black!60] (0.04,0) -- (2.2,0) node[above] {$r_a$};
   \draw (0.7,-0.05) -- (0.7,0.05);
   \draw (1.4,-0.05) -- (1.4,0.05);
   \node[below,inner sep=5pt] at (0.7,0) {$\rho$};
   \node[below,inner sep=5pt] at (1.4,0) {$2\rho$};
   \node at (2.05,1.5) {$U$};
   \node at (-1.72,1.7) {$L_a$};
   \draw[->,black!60] (-1.5,1.58) -- (-0.99,0.99);
   \node at (-1.68,-2.2) {collar};
   \draw[->,black!60] (-1.4,-2.05) -- (-0.95,-1.25);
 \end{scope}
\end{tikzpicture}
\caption{The exterior $U$ is blue, and the regular collars are hatched.}
\label{fig:compactexterior}
\end{figure}
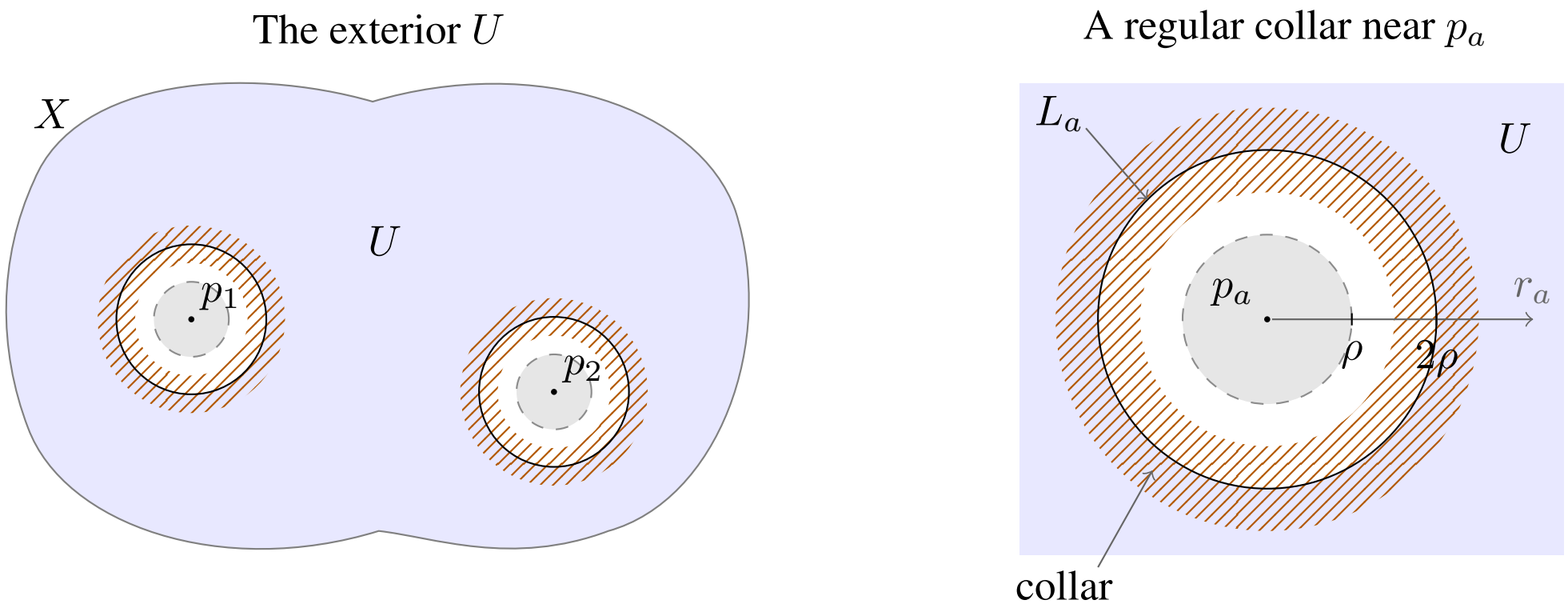

We describe the exterior and its collars explicitly. Let $X$ be a
compact $n$-dimensional orbifold with isolated singular points $p_a$.
Choose pairwise disjoint orbifold coordinate balls
$B_a(3\rho)\cong B_{\R^n}(0,3\rho)/\Gamma_a$, and write $r_a$ for
the radial coordinates. Set
\[
 U:=X\setminus\bigcup_a\overline{B_a(2\rho)},\qquad
 U^+:=X\setminus\bigcup_a\overline{B_a(\rho)}.
\]
Then $\overline U\subset U^+\Subset \mathcal R(X)$. The closure
$\overline U$ is a compact smooth manifold with boundary components
$L_a:=\{r_a=2\rho\}\cong S^{n-1}/\Gamma_a$. Each $L_a$ has a
two-sided regular collar
\[
 \{3\rho/2<r_a<5\rho/2\}
 \cong L_a\times(-\rho/2,\rho/2),
 \qquad s:=r_a-2\rho.
\]
Its outer half $s>0$ lies in $U$; its inner half $s<0$ lies in
the removed ball but is still regular. The portion $s\ge0$ is a
boundary collar of $\overline U$. The whole collar lies in $U^+$,
which leaves a regular neighborhood around $\overline U$ for local
estimates; see Figure~\ref{fig:compactexterior}. These are coordinate
balls, not necessarily geodesic balls. If $X$ is smooth, take
$U=U^+=X$.
In Figure~\ref{fig:compactexterior}, the solid inner boundaries
represent $L_a$, and the gray cores represent $\overline{B_a(\rho)}$.
Thus $U^+$ is the region outside the gray cores. The outer outline
is schematic; $X$ has no boundary.

\begin{lemma}[Compact exterior extension]\label{lem:compactextension}
Let $X,Y$ be connected compact oriented effective real-analytic
Riemannian orbifolds of dimension $n\ge3$. Assume that their
singularities are isolated and their volumes are equal. Let
$U\subset \mathcal R(X)$ be the connected open exterior of small
pairwise disjoint closed orbifold coordinate balls around all
singular points, with regular collars. If $X$ is smooth, take $U=X$.
Every orientation-preserving local isometry
$f:U\to \mathcal R(Y)$ extends uniquely to an orbifold isometry
$F:X\to Y$. If the metrics are K\"ahler and $f$ is holomorphic,
so is $F$.
\end{lemma}
\begin{proof}
The oriented orthonormal frame bundles $FX,FY$ are connected compact
analytic smooth manifolds. Each effective local group acts freely
on frames, so this construction does not require a manifold cover
of either orbifold \cite[Theorem~1.2]{Kank}. Equip them with the
analytic connection metrics
\[
 G:=\langle\theta,\theta\rangle+\langle\eta,\eta\rangle,
\]
where $\theta$ is the tautological form and $\eta$ is the
Levi-Civita connection form. Use the same fixed
$\operatorname{Ad}$-invariant inner product on $\mathfrak{so}(n)$
for both metrics. They are complete by compactness. The map $f$
lifts to a local isometry $\widehat f:FU\to FY$ preserving these
forms.

The frame fiber over an isolated quotient point with group
$\Gamma$ is an embedded copy of $SO(n)/\Gamma$, of codimension
$n$ in $FX$. Removing these finitely many fibers does not change
$\pi_1$: loops and their two-dimensional homotopies can be moved
off them since $n\ge3$. Radial deformation in the punctured
coordinate balls, lifted to frames, identifies $FU$ with this
complement up to homotopy. Consequently
\[
 \pi_1(FU)\longrightarrow\pi_1(FX)
\]
is an isomorphism.

Lift a germ of $\widehat f$ to the complete simply connected
covers of $FX,FY$. The analytic continuation theorem
\cite[Corollary~6.4, p.~256]{KN} extends it to an isometry of
these covers. Thus the original germ continues along every path
in $FX$, with homotopy-invariant continuation. Every loop based
in $FU$ is homotopic to a loop contained in $FU$. Along such a
loop the continuation follows the already defined map
$\widehat f$ and returns to the original germ. It is therefore
single-valued on $FX$, giving a local isometry
$\widehat F:FX\to FY$. This use of an open set carrying the
fundamental group is also explained in
\cite[Section~4.2, p.~2900]{GGSU}.

Compactness and connectedness make $\widehat F$ a finite covering.
With the common vertical normalization,
\[
 \Vol(FX,G_X)=\Vol(SO(n))\Vol(X),\qquad
 \Vol(FY,G_Y)=\Vol(SO(n))\Vol(Y).
\]
These identities follow by integration over the regular locus;
the exceptional fibers have measure zero. Since a local isometry
preserves volume, the covering degree is
\[
 \deg\widehat F
 =\frac{\Vol(FX,G_X)}{\Vol(FY,G_Y)}=1.
\]
Hence $\widehat F$ is a global isometry.

For each $A\in SO(n)$, let $R_A$ denote the right action on either
frame bundle. The maps $\widehat F\circ R_A$ and
$R_A\circ\widehat F$ agree on $FU$. Analytic continuation makes
them equal everywhere. Thus $\widehat F$ descends to an orbifold
isometry $F:X\to Y$ extending $f$. Uniqueness follows from the
same continuation argument. In the K\"ahler case
$F^*J_Y-J_X$ is parallel and vanishes on $U$, hence vanishes on
the connected regular locus and on its orbifold completion.
\end{proof}

Only the surjectivity of $\pi_1(FU)\to\pi_1(FX)$ is needed.
An arbitrary regular patch need not have this property. The same
exterior $U$ remains admissible for a pullback metric on $X$;
its deleted coordinate balls need not be geodesic balls for that
metric.

\subsection{Short-time limits}

\begin{lemma}[Short-time limits with vanishing energy drop]
\label{lem:shortdrop}
Fix $T,r_0>0$. Let $g_j(t)$, $-3T\le t\le3T$, be smooth Calabi
flows on $(M,J)$ in $\Omega$, satisfying \eqref{eq:Sob} and
\begin{equation}\label{eq:segmentbounds}
 \sup_{|t|\le3T}\int_M|\Rm(g_j(t))|^2\,dV\le E_0,
 \qquad
 \inf_{M\times[-2T,2T]}\scr_j\ge r_0.
\end{equation}
Suppose $\Ca_j(-3T)-\Ca_j(3T)\to0$.
On compact subintervals of $(-2T,2T)$, the following properties hold.
Slice and local cylinder convergence is subsequential, also at arbitrary
sequences of moving times:
\begin{enumerate}[label=\textup{(B\arabic*)}]
\item The completed slice limits are compact extremal K\"ahler
orbifolds. Convergence is global in the GH sense and smooth on
the entire geometric regular locus. For each
$K\Subset U\Subset U^+\Subset \mathcal R(X)$, comparison
embeddings are defined on $U^+$, and their images of $U$
contain all points of the underlying manifold corresponding to $K$ for large $j$.
The natural volume measures converge, with total limiting volume $V$.
\item For fixed terminal domains
$K\Subset U\Subset U^+\Subset \mathcal R(X)$, terminal
embeddings on $U^+$ give smooth local Calabi-flow limits on
$U\times[-\sigma,0]$ for some $\sigma>0$, after translating
the terminal time to zero. The terminal metric separations
and $\sigma$ are bounded below independently of $j$, and
uniformly for isometric copies of these nested domains with
terminal times in the chosen compact time interval.
\item If $a_j,b_j$ lie in the chosen compact subinterval and
$a_j-b_j\to0$, then for any sequences of points $x_j,y_j$ of the underlying manifold,
held fixed as time varies,
\[
 \left|d_{g_j(a_j)}(x_j,y_j)-d_{g_j(b_j)}(x_j,y_j)\right|\to0.
\]
\end{enumerate}
\end{lemma}
\begin{proof}
The dissipation identity gives
\[
 \int_{-3T}^{3T}D_j(t)\,dt
 =\tfrac12\bigl(\Ca_j(-3T)-\Ca_j(3T)\bigr)\longrightarrow0.
\]
Fix $I_0\Subset I_1\Subset(-2T,2T)$ and choose testing scales
$s<r_0/4$ with $s^4<\operatorname{dist}(I_0,\partial I_1)$.
Thus every backward testing interval based in $I_0$ lies in $I_1$. Tests (R) and (V) then hold with the
same constants for every $j$. On balls satisfying the energy trigger,
noncollapse and (R) give finite regular charts. The interior estimate~\eqref{eq:interiorsmoothing} supplies every prescribed derivative
order on smaller cylinders, still with uniform constants.

Proposition~\ref{prop:slice} now gives extremal multifold slice limits.
The model comparisons (P) and the two gaps in
Lemma~\ref{lem:modelgaps} exclude branching and curvature concentration
over geometrically regular points, exactly as in the proof of
Corollary~\ref{cor:geometricinputs}. The upper volume bound and
Lemma~\ref{lem:volume-passage} give the natural limiting measure and
total volume $V$. This proves (B1). For (B2), fix
$K\Subset U\Subset U^+\Subset \mathcal R(X)$ and cover
$\overline U$ by terminal balls whose concentric larger balls
lie in $U^+$. We use the argument of
Corollary~\ref{cor:backwardpseudolocality}: the terminal curvature
bounds, (V), and the uniform lower bound $r_0$ give a common
scale at which (R) applies. Its proof therefore gives backward
estimates uniformly in $j$. The bootstrap on the smaller balls
during the latter half of each backward interval gives all
derivative bounds. The sizes depend on the intrinsic
regular geometry near $\overline{U^+}$, the terminal metric
separations, $r_0$, and $\operatorname{dist}(I_0,\partial I_1)$.
They are therefore uniform in $j$ and for isometric copies
of these nested domains.

For (B3), extract the later of the two moving slices. Its regular
exterior and the collars of its finitely many singular cores have
short backward estimates from (B2). The boundary-flux argument
\eqref{eq:selectedflux} preserves the small core volumes when the
time difference tends to zero. Noncollapse gives the core diameter
bound \eqref{eq:selectedcorediam}. Collapsing the cores and comparing
the exterior metrics gives \eqref{eq:selectedIDC}, with no outer
boundary in this compact setting. First let $j\to\infty$, then let
the core radii tend to zero. This proves (B3).
\end{proof}

\subsection{Long-time curvature estimates}

The forward estimate requires an identification of the same points
at different times. We obtain it by starting at a later terminal
slice and using backward regularity. The flow is assumed to exist
on the entire displayed interval.

\begin{lemma}[Long-time forward pseudo-locality]\label{lem:forwardpersistence}
Fix $T,r_0,r>0$. There are constants
\[
 \eps=\eps(M,J,\Omega,A,E_0,T,r_0,r)>0,\qquad
 C_m=C(m,M,J,\Omega,A,E_0,T,r_0,r)<\infty,
\]
for $m\ge0$, with the following property. Let $g(t)$,
$-3T\le t\le3T$, be a smooth
Calabi flow in $\Omega$ satisfying \eqref{eq:Sob} and the bounds
\eqref{eq:segmentbounds}, with the subscript $j$ omitted. If
\begin{equation}\label{eq:smalldrop}
 \Ca(-3T)-\Ca(3T)<\eps,
 \qquad
 \sup_{B_{g(0)}(x,2r)}|\Rm(g(0))|\le r^{-2},
\end{equation}
then
\begin{equation}\label{eq:forwardpersistence}
 \sup_{B_{g(0)}(x,r)\times[0,T]}
 r^{2+m}|\nabla^m\Rm(g(t))|\le C_m,
 \qquad m\ge0.
\end{equation}
The spatial ball is fixed at time zero. The same $\eps$ works for
every derivative order.
\end{lemma}
\begin{proof}
\noindent\textit{Step 1: Backward construction from a terminal slice.}
We first prove that some $\eps_0,B>0$ give
\begin{equation}\label{eq:finiteforwardpersistence}
 \max_{0\le k\le N}
 \sup_{B_{g(0)}(x,3r/2)\times[-T/4,5T/4]}
 r^{2+k}|\nabla^k\Rm|\le B,
 \qquad N=12,
\end{equation}
under the hypotheses with energy drop less than $\eps_0$.
If this fails, choose flows $g_j$, centers $p_j$, points $q_j$ of the underlying manifold,
held fixed as time varies, with $q_j\in B_{g_j(0)}(p_j,3r/2)$, times $s_j\in[-T/4,5T/4]$, and
orders $k_j\le N$, such that
\[
 \Ca_j(-3T)-\Ca_j(3T)\longrightarrow0,
 \qquad
 r^{2+k_j}|\nabla^{k_j}\Rm(g_j(s_j))|(q_j)\longrightarrow\infty.
\]
Lemma~\ref{lem:shortdrop} supplies (B1)--(B3).
Before using the selected points and times, we establish the following
consequence of these properties on any compact interval
$[a,b]\Subset(-2T,2T)$. After passing to a subsequence, the completed
slice limits can be identified with the extremal family generated by
the terminal metric at $b$. Under these identifications, every fixed
sequence of points of the underlying manifolds has a constant limiting
point, after further extraction if necessary. Convergence on compact
regular domains is smooth in embeddings independent of time.
We shall also use this conclusion in the backward part of
Theorem~\ref{thm:upgrade}. The argument below uses (B1)--(B3) and
vanishing dissipation, without invoking spacetime weak compactness.

Translate the terminal time $b$ to zero, so $[a,b]$ becomes
$[a-b,0]$, and write $(X,h_0,J)$ for a completed terminal limit and
$\mathcal S:=\Sing X$. We work only on backward intervals whose original
times lie in $(-2T,2T)$. By (B1) and Lemma~\ref{lem:analytic}, the
completed slice limits are analytic extremal orbifolds.

Choose the exterior $U\Subset U^+$ as above, and a further
regular domain $U^+\Subset U^{++}\Subset \mathcal R(X)$;
see Figure~\ref{fig:compactexterior} and
Lemma~\ref{lem:compactextension}. If $\mathcal S=\varnothing$,
take $U=U^+=U^{++}=X$, with no spatial boundary.
Choose terminal embeddings $\Phi_j$ on an exhaustion of
$\mathcal R(X)$, initially restricted to $U^{++}$. A finite
cover and (B2) give smooth limits on $U^+\times I$ for one
backward interval $I:=[-\sigma,0]$.
Proposition~\ref{prop:local} and vanishing dissipation show that
every slice of this local limit is extremal.
Let $F_\tau$ be the complete holomorphic flow of
$W:=\tfrac12\nabla_{h_0}S_{h_0}$. Decrease $\sigma$ so that
characteristics starting in $\overline U$ remain in $U^+$.
Lemma~\ref{lem:transport} identifies the local limit on $U$ with
\begin{equation}\label{eq:shortsoliton}
h_\tau=F_\tau^*h_0.
\end{equation}

Lemma~\ref{lem:regularorbits} shows that $F_\tau$ preserves
$\mathcal R(X)$ and fixes $\mathcal S$. On any finite closed interval,
the trajectories of a compact regular set remain a positive
distance from $\mathcal S$ when $\mathcal S\ne\varnothing$. The lemma
also bounds the intrinsic curvature derivatives of $h_\tau$
uniformly in time.

Figure~\ref{fig:regularorbits} shows this preservation for a compact
regular region $K$. The blue sections represent $F_\tau(K)$ in
the fixed reference orbifold $(X,h_0)$; the red lines represent
fixed singular points. Their separation is positive on each finite
time interval. The identification of fixed points in the approximating
flows will be proved below.

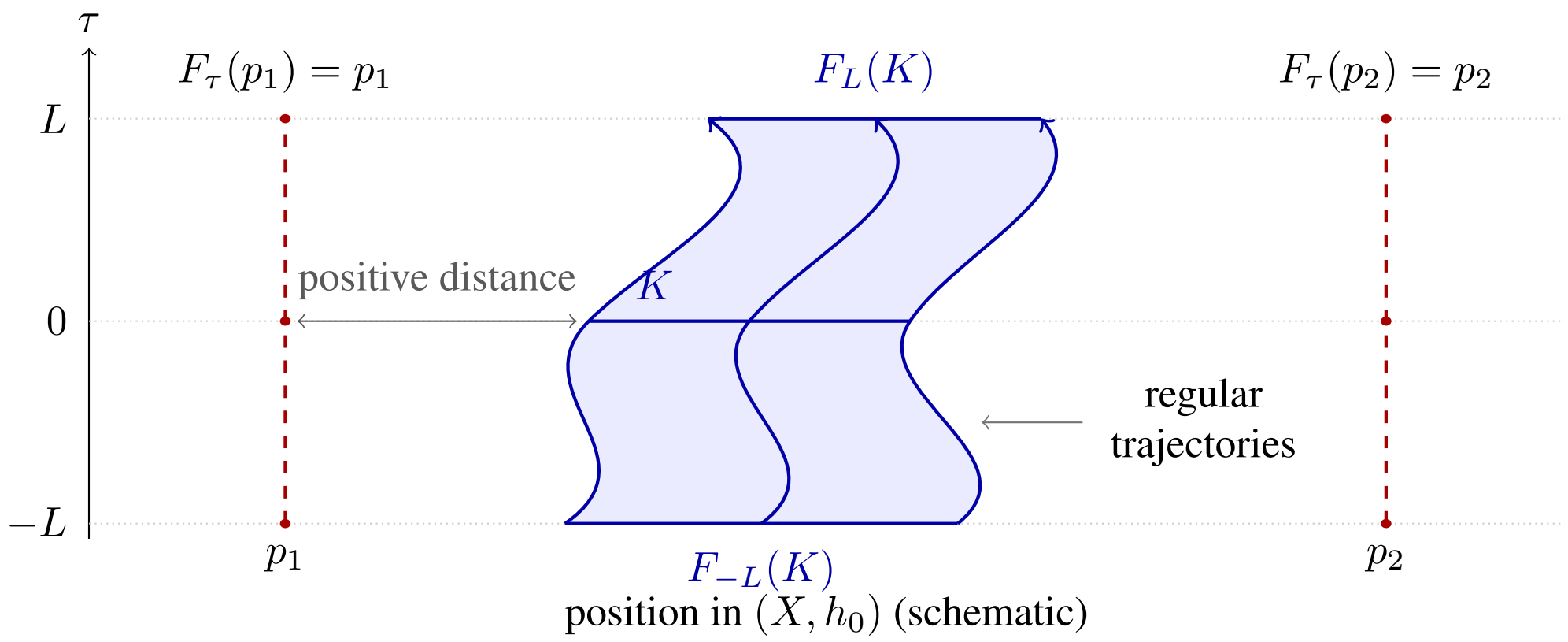
\begin{figure}[htbp]
\centering
\begin{tikzpicture}[x=1cm,y=0.85cm,font=\small]
 \draw[->] (0,-0.15) -- (0,4.7) node[above] {$\tau$};
 \foreach \y/\t in {0/{-L},2/{0},4/{L}} {
   \draw[black!20,densely dotted] (0,\y) -- (12.4,\y);
   \node[left,inner sep=5pt] at (0,\y) {$\t$};
 }
 \fill[blue!8]
 (4.0,0) .. controls (4.8,0.65) and (3.6,1.3) .. (4.2,2)
 .. controls (4.8,2.7) and (6.0,3.3) .. (5.2,4)
 -- (8.0,4)
 .. controls (8.5,3.3) and (7.3,2.7) .. (6.9,2)
 .. controls (6.5,1.3) and (8.0,0.65) .. (7.3,0) -- cycle;
 \foreach \x/\a in {1.65/1,10.9/2} {
   \draw[red!65!black,dashed,thick] (\x,0) -- (\x,4);
   \foreach \y in {0,2,4} {
     \fill[red!65!black] (\x,\y) circle (0.045);
   }
   \node[below,inner sep=5pt] at (\x,0) {$p_{\a}$};
   \node[above,inner sep=6pt] at (\x,4) {$F_\tau(p_{\a})=p_{\a}$};
 }
 \draw[blue!65!black,thick,->]
 (4.0,0) .. controls (4.8,0.65) and (3.6,1.3) .. (4.2,2)
 .. controls (4.8,2.7) and (6.0,3.3) .. (5.2,4);
 \draw[blue!65!black,thick,->]
 (7.3,0) .. controls (8.0,0.65) and (6.5,1.3) .. (6.9,2)
 .. controls (7.3,2.7) and (8.5,3.3) .. (8.0,4);
 \draw[blue!65!black,thick,->]
 (5.65,0) .. controls (6.4,0.65) and (5.05,1.3) .. (5.55,2)
 .. controls (6.05,2.7) and (7.25,3.3) .. (6.6,4);
 \draw[blue!65!black,thick] (4.0,0) -- (7.3,0);
 \draw[blue!65!black,thick] (4.2,2) -- (6.9,2);
 \draw[blue!65!black,thick] (5.2,4) -- (8.0,4);
 \node[blue!65!black,below,inner sep=6pt] at (5.65,0) {$F_{-L}(K)$};
 \node[blue!65!black,above,inner sep=5pt] at (4.75,2) {$K$};
 \node[blue!65!black,above,inner sep=6pt] at (6.6,4) {$F_L(K)$};
 \draw[<->,black!65] (1.75,2) -- (4.1,2)
   node[midway,above,inner sep=5pt] {positive distance};
 \node[align=center,anchor=west] at (8.45,1.0) {regular\\trajectories};
 \draw[->,black!55] (8.35,1.0) -- (7.5,1.0);
 \node at (6.2,-0.9) {position in $(X,h_0)$ (schematic)};
\end{tikzpicture}
\caption{On finite time intervals, $F_\tau(K)$ stays separated from the fixed singular points.}
\label{fig:regularorbits}
\end{figure}

At arbitrary moving times $\tau_j\to\tau\in I$, (B1) gives a completed
slice $(Y_\tau,k_\tau)$. Compose its GH approximation maps with the
fixed embeddings on $U^+$. After a subsequence these maps converge
locally uniformly on $U$. The inclusion
$\overline U\subset U^+$ has positive $h_0$-distance from
$X\setminus U^+$. Uniform metric equivalence on
$\overline U\times I$ and on a regular neighborhood of
$\overline U$ in $U^+$ makes this separation positive in the
pullback metrics, independently of $j$ and $\tau\in I$.
First exit therefore keeps sufficiently short ambient minimizing
curves in the chart image and places a small ambient ball inside it. The limit is therefore
a local isometry $(U,h_\tau)\to(Y_\tau,k_\tau)$ with regular image.
Transition maps to the regular charts in (B1), together with convergence
of the complex structures, show that it is holomorphic.
The natural volume convergence in (B1) gives
$\Vol(Y_\tau,k_\tau)=V$, while
$\Vol(X,h_\tau)=\Vol(X,h_0)=V$ by pullback invariance.
Lemma~\ref{lem:compactextension} therefore extends this map uniquely to
\[
 I_\tau:(X,h_\tau,J)\longrightarrow(Y_\tau,k_\tau,J_\tau).
\]
We use $I_\tau$ to identify each completed slice with $X$, agreeing
with the fixed charts on $U$. Moreover, $I_\tau\circ F_{-\tau}$ is
a holomorphic isometry from $(X,h_0)$ onto $(Y_\tau,k_\tau)$.
The limits of fixed points of the underlying manifold at different
times have not yet been identified. We prove this next.

\medskip
\noindent\textit{Step 2: Identification of points through time.}
We take a countable dense set $\{a_l\}\subset U$. For any sequence of points $x_j$ of the underlying manifold,
held fixed as time varies, the functions
\[
f_{j,l}(\tau):=d_{g_j(\tau)}\bigl(x_j,\Phi_j(a_l)\bigr)
\]
are bounded by the diameter estimate and asymptotically equicontinuous by (B3). A diagonal subsequence converges uniformly, for each $l$, to a continuous $f_l$. Moving-slice compactness and the maps $I_\tau$ represent these values as distances from a point $\gamma(\tau)\in X$ to $a_l$, measured by $h_\tau$.

Distances to a nonempty regular open set distinguish points. Indeed, if two distance functions agree there, we choose a point where both are differentiable. Equality of the gradients and lengths gives the same initial direction and length for a minimizing geodesic, hence the same endpoint. A minimizing geodesic cannot traverse an isolated nontrivial quotient point in its interior, because the spherical quotient link has diameter less than $\pi$ and the corner can be shortened. Singular endpoints cause no ambiguity. Thus these distances determine $\gamma(\tau)$ uniquely. The continuous injection
\[
I\times X\longrightarrow I\times\R^{\mathbb N},\qquad
(\tau,x)\longmapsto
\bigl(\tau,(d_{h_\tau}(x,a_l))_l\bigr)
\]
is an embedding by compactness. It follows that $\gamma$ is continuous.
For any $\tau_j\to\tau$, uniform convergence gives
$f_{j,l}(\tau_j)\to f_l(\tau)$ for every $l$. Hence every
normalized moving-time limit of $x_j$ is $\gamma(\tau)$.
Further subsequences used to obtain fresh charts do not change
this path.

Suppose $x:=\gamma(\tau_0)$ is regular. Choose connected regular
domains
$x\in D\Subset D^+\Subset D^{++}\Subset \mathcal R(X)$
and nonempty open sets $O\Subset O^+\Subset D\cap U$.
At $\tau_0$, use $I_{\tau_0}$ to identify the domains of fresh
terminal embeddings $\Psi_j$ on $D^{++}$. Condition (B2)
gives smooth local limits on $D^+$ during a fixed backward
interval. Vanishing dissipation and Lemma~\ref{lem:transport}
identify the limit on $D$ with $h_\tau$, after shortening the
interval so characteristics from $\overline D$ remain in $D^+$.
On $O$,
$\Phi_j^{-1}\Psi_j$ is defined by terminal coverage and converges
smoothly to the identity by Lemma~\ref{lem:retainmaps}. This
transition is independent of $\tau$, so the two charts give the
same limiting identification on $O$ throughout their common cylinder.

At each earlier time, the $\Psi_j$ charts induce a local isometry
from $(D,h_\tau)$ into $(X,h_\tau)$ under the global identification.
The estimates on $D^+$ and
$d_{h_{\tau_0}}(\overline D,X\setminus D^+)>0$ give a
positive lower bound, independent of $j$ and time in this
short interval, for the length of a path starting in
$\Psi_j(D)$ and leaving $\Psi_j(D^+)$. Thus small
ambient balls lie in its image, so its limiting image is regular.
This uses only the fresh charts. The local isometry is the identity
on $O$, hence on connected $D$, by uniqueness of local isometries.
Terminal coverage places $x_j$ in $\Psi_j(D)$, with
$\Psi_j^{-1}(x_j)\to x$. These coordinates are independent of
time. Consequently $\gamma(\tau)=x$ on a backward neighborhood
of $\tau_0$. These local charts establish constancy before persistence of the
original $\Phi_j$ away from $O$ is known.

No uniform size is required for these neighborhoods near $\mathcal S$.
On each connected component of $\gamma^{-1}(X\setminus\mathcal S)$,
the path is constant:
otherwise the first time at which its distance from one of its values
reaches an intermediate positive value contradicts local backward
constancy. At an endpoint of the component, continuity gives the same
regular value. Hence the component cannot end in the interior of $I$ and must equal
$I$. If the path has no regular value, its image lies in the finite
set $\mathcal S$, so connectedness again makes it constant.
Thus every fixed sequence of points has the same limit throughout $I$.
In particular, regularity at an earlier time prevents its limit from
becoming singular at a later time. This conclusion uses no forward
regularity estimate for the approximating flows.

\medskip
\noindent\textit{Step 3: Smooth convergence in the original embeddings.}
Fix $K\Subset U_0\Subset U_1\Subset \mathcal R(X)$ and retain
the original $\Phi_j$ on a neighborhood of $\overline{U_1}$.
Every sequence $\Phi_j(z_j)$ with $z_j\in\overline{U_1}$ has
regular limit points at arbitrary moving times, by the preceding
identification of fixed points of the underlying manifold. If an intrinsic curvature
derivative were unbounded there, (B1), normalized by the exterior
extension, would give smooth covered charts at the offending
moving points, a contradiction. Thus all intrinsic derivatives
are uniformly bounded on this fixed region of the underlying manifold throughout $I$.

Put $G_j:=\Phi_j^*g_j$ and $T_j:=\partial_\tau G_j$. The Calabi
equation bounds every $|\nabla^mT_j|_{G_j}$, and integration of
the zeroth bound gives uniform equivalence to $G_j(0)$. The
connection induction in Step~6 of Theorem~\ref{thm:selected}
then bounds all spatial derivatives in these same coordinates.
Since $\Phi_j^*J_j$ is fixed in time and converges smoothly,
successive differentiation of the Calabi equation bounds every
mixed derivative. Hence $G_j$ is smoothly precompact on
$\overline{U_0}\times I$.

Let $G$ be a subsequential limit. For $z\in K$, choose
$B_{h_0}(z,a)\subset B_{h_0}(z,2a)\Subset U_0$.
First exit from the smaller ball, using metric equivalence on
the larger one throughout $I$, ensures that sufficiently short
source minimizers remain in its image. The preceding identification of fixed point sequences, applied also
at moving times, identifies the resulting local distances of $G(\tau)$
with those of $h_\tau$ under the identity map. Therefore $G(\tau)=h_\tau$. Every smooth
subsequential limit is the same, proving smooth convergence in
the original embeddings, including mixed derivatives.
The inclusions $K\Subset U_0\Subset U_1$ and the same
first-exit argument give coverage on $K$.

\medskip
\noindent\textit{Step 4: Continuation and the forward estimate.}
At a restart time $c$ in the interior of the preceding interval,
measured relative to the terminal time $b$ (so the original time is $b+c$),
take $U_c:=F_c^{-1}(U)$, $U_c^+:=F_c^{-1}(U^+)$, and
$U_c^{++}:=F_c^{-1}(U^{++})$. The
holomorphic isometry $F_c:(X,h_c)\to(X,h_0)$ preserves the
intrinsic exterior data and carries $W$ to itself. By (B2),
both the regular-chart lifetime and the characteristic lifetime
have the same positive lower bound at every restart. Reduce this
bound below half the distance of $[a,b]$ from the boundary of
$(-2T,2T)$.
Fresh embeddings are identified at $c$ using the preceding identification
of the completed slice limits. On fixed domains
$V\Subset V^+$ in the common regular region, where both
embeddings are defined on $V^+$, their transitions on $V$
are independent of time and converge smoothly to the identity
by Lemma~\ref{lem:retainmaps}.

Apply the preceding argument on the new interval with this
normalization. On the overlap, the limits of each fixed point
sequence agree at $c$ and are constant in both constructions;
hence they agree throughout the overlap. The limiting local
flows are the same $h_\tau$ by Lemma~\ref{lem:transport} and uniqueness of
exterior extension. Lemma~\ref{lem:retainmaps} transfers smooth
convergence back to the original $\Phi_j$ on every fixed regular
compact set. Finitely many overlapping intervals cover the prescribed
backward interval $[a,b]$. Only auxiliary charts are replaced; the retained
embeddings remain fixed. Lemma~\ref{lem:regularorbits} places all
$F_{-c}(\overline{U^{++}})$ in one compact subset of $\mathcal R(X)$
as $c$ ranges over the prescribed finite interval. Thus all nested exterior domains used at
the restarts lie in one member of the original exhaustion.
Their separations in $h_c$ equal those of the original domains
in $h_0$, so they do not decrease with the number of restarts.

This proves the announced identification on $[a,b]$.
For the forward estimate, choose the original terminal time $b=3T/2$
and continue backward to $a=-T/2$. The terminal limits of the selected
points $q_j$ need not be assumed regular. Instead, at the original
reference time zero these points have regular limit points. Indeed, their distance from the boundary of
$B_{g_j(0)}(p_j,2r)$ is at least $r/2$. Choose a fixed
$0<\ell<\min\{r/8,r_0/4,V^{1/4}/8\}$ so small that
$B_0(2\ell)^4r^{-4}<\eps_R$ and $\ell^4<T/4$, so the
backward testing interval lies inside $(-2T,2T)$.
Then (V) and the reference curvature bound give
\[
 \int_{B_{g_j(0)}(q_j,2\ell)}|\Rm|^2<\eps_R.
\]
Test (R) and the bootstrap give all terminal derivative bounds on
smaller balls about $q_j$. Noncollapse gives a positive local
injectivity radius. These charts, with ambient ball coverage by first exit, show that
the corresponding points in every completed time-zero limit are regular.

The identification just proved keeps the limit of each fixed point
$q_j$ unchanged throughout $[-T/2,3T/2]$. After taking $s_j\to s$,
its limit at the moving time $s_j$ is therefore regular as well.
Property (B1), including coverage, bounds every curvature derivative
at $q_j$ at that time. Since $k_j$ takes only finitely many values,
this contradicts its selection and proves \eqref{eq:finiteforwardpersistence}.

Set $\eps:=\eps_0$. The order-two bound controls
$|\partial_tg|\le C|\nabla^2S|$, so the metrics are uniformly equivalent
on the larger fixed ball throughout $[-T/4,5T/4]$.
The reference curvature bound, (R), and noncollapse give regular
holomorphic charts at time zero on coordinate balls
$B_2\subset\mathbb C^2$, with uniform finite bounds.
Choose them inside the fixed source ball
$B_{g(0)}(x,3r/2)$ so that their concentric $B_1$-images
cover $B_{g(0)}(x,r)$. Their sizes and terminal metric
separations have positive lower bounds depending only on the
fixed data. Retain these charts at all times.
The connection-variation induction turns the bounds through
order $N=12$ into uniform spatial $C^8$ metric bounds on
these $B_2$-charts. Metric equivalence and first exit preserve
the positive separation of the $B_1$-images from the boundaries
of the $B_2$-images throughout $[-T/4,5T/4]$.
Apply the interior estimate~\eqref{eq:interiorsmoothing} on the inner coordinate balls
and on $[0,T]$, using the estimates on the larger balls and
the additional earlier interval $[-T/4,0]$ of length $T/4$.
This gives every prescribed derivative bound on
$B_{g(0)}(x,r)\times[0,T]$. This proves
\eqref{eq:forwardpersistence}. The finite-order constants $\eps_0,B$
depend only on $(M,J,\Omega,A,E_0,T,r_0,r)$ with the canonical-radius
choices fixed above; the bootstrap adds only dependence on $m$.
Thus one energy threshold works for every derivative order.
\end{proof}

\begin{samepage}
\begin{theorem}[Long-time two-sided pseudo-locality]\label{thm:upgrade}
Under the geometric hypotheses of Lemma~\ref{lem:forwardpersistence},
there are constants
\[
 \eps=\eps(M,J,\Omega,A,E_0,T,r_0,r)>0,\qquad
 C_m=C(m,M,J,\Omega,A,E_0,T,r_0,r)<\infty,
\]
such that \eqref{eq:smalldrop} implies
\begin{equation}\label{eq:persistence}
 \sup_{B_{g(0)}(x,r)\times[-T,T]}
 r^{2+m}|\nabla^m\Rm(g(t))|\le C_m,
 \qquad m\ge0.
\end{equation}
The spatial ball is fixed at time zero, and the same $\eps$ works
for every derivative order.
\end{theorem}
\end{samepage}
\begin{proof}
Lemma~\ref{lem:forwardpersistence} gives the estimate on $[0,T]$.
For the backward estimate, first obtain the finite-order bound
\begin{equation}\label{eq:finitepersistence}
 \max_{0\le k\le N}
 \sup_{B_{g(0)}(x,3r/2)\times[-5T/4,0]}
 r^{2+k}|\nabla^k\Rm|\le B_-,
 \qquad N=12,
\end{equation}
under a sufficiently small energy-drop threshold $\eps_->0$.
If no such constants exist, take a sequence with energy drop tending
to zero and points $q_j$ in the displayed fixed ball at which one of
these derivatives becomes unbounded at times $s_j\in[-5T/4,0]$.
The reference curvature bound, (R), and noncollapse show that every
subsequential time-zero limit of $q_j$ is regular, as in the proof of
Lemma~\ref{lem:forwardpersistence}.

Apply the backward construction proved there with terminal time
$b=0$ on $[-3T/2,0]$. It identifies the completed slice limits with
the terminal extremal family and identifies each fixed sequence
$q_j$ with the same regular point throughout this interval.
After taking $s_j\to s$, (B1), including coverage at moving times,
therefore bounds its curvature derivatives at $s_j$. This contradiction
proves \eqref{eq:finitepersistence}.

The order-two bound gives metric equivalence on the fixed ball.
Using the time-zero charts and the connection-variation argument in
Lemma~\ref{lem:forwardpersistence}, the finite-order bound gives
uniform spatial $C^8$ metric bounds. The interior bootstrap on
the smaller ball, with the earlier interval $[-5T/4,-T]$, gives all
curvature derivatives on $[-T,0]$. Its energy threshold is independent
of the derivative order. Taking the minimum of $\eps_-$ and the
forward threshold, and the maximum of the corresponding constants,
proves \eqref{eq:persistence} with the stated dependence.
\end{proof}

For the original global flow, Theorem~\ref{thm:radiuslower} supplies
$r_0$ on every sufficiently late translated interval. For each fixed
$T$, Lemma~\ref{lem:tail} makes the energy drop on
$[t_j-3T,t_j+3T]$ smaller than the threshold above. Thus the estimate
is available before extracting a spacetime limit.

% END proofs/compact_spacetime_upgrade_section.tex

\section{Weak compactness}
\label{sec:mainproof}
\label{sec:topology}
The long-time two-sided pseudo-locality theorem supplies curvature
estimates on fixed regular regions throughout bounded translated
intervals. We first extract a limit flow and prove
Theorem~\ref{thm:main}\textup{(I)--(III)}. We then prove
continuity of the Calabi functional, completing \textup{(IV)},
and finish with the finite-group equivariant version.

The motivation for the topology below is the pointed
$\hat C^\infty$-CG convergence recalled in
Section~\ref{sec:pointedconvergence}, Definition~\ref{def:modelhat}.
It combines GH convergence with smooth convergence on the regular
part through compatible comparisons. For flows, these comparisons
must also identify the same points of the underlying manifolds at
different times. We therefore require time-independent relations and
regular embeddings, with uniform distance control and smooth
convergence on every compact time interval. This gives the spacetime
version, abbreviated as \emph{spacetime $\hat C^\infty$-CG convergence},
and preserves the identification of points through time as in
Chen--Wang \cite{ChenWangA,ChenWangB}.

\begin{definition}[Spacetime pointed-$\hat C^\infty$-Cheeger--Gromov convergence]
\label{def:topology}
\leavevmode\par
Let $(X_j,p_j,g_j(\tau),J_j)$ and $(X,p,h(\tau),J_\infty)$ be
families of complete K\"ahler orbifolds on a time interval $\mathcal I$.
Each underlying orbifold and complex structure is fixed in time.
The metrics are smooth in space and time in orbifold charts, and
their length distances $d_j(\tau)$ and $d(\tau)$ are jointly
continuous. Basepoints may be singular.

We say that the former converge to the latter in the spacetime
pointed-$\hat C^\infty$-Cheeger--Gromov topology if the following
holds on every compact interval $I\subset\mathcal I$. There are
$R_j\to\infty$, $\epsilon_j\to0$, time-independent relations
$\mathscr C_j\subset X_j\times X$ containing $(p_j,p)$, and a
regular exhaustion $U_j\Subset\mathcal R(X)$, with these properties:
\begin{enumerate}[label=\textup{(\roman*)}]
\item For every $\tau\in I$, the two projections of $\mathscr C_j$
contain respectively $B_{d_j(\tau)}(p_j,R_j)$ and
$B_{d(\tau)}(p,R_j)$, and
\begin{equation}\label{eq:distortion}
 \sup_{\tau\in I}\operatorname{dis}_{\tau}(\mathscr C_j)
 \le\epsilon_j,\qquad  \operatorname{dis}_{\tau}(\mathscr C_j)
 :=\sup_{(x,a),(y,b)\in\mathscr C_j}
 |d_j(\tau)(x,y)-d(\tau)(a,b)|.
\end{equation}

\item There are smooth embeddings
$\Phi_j:U_j\to\mathcal R(X_j)$, independent of $\tau$, such that
$(\Phi_j(x),x)\in\mathscr C_j$ for every $x\in U_j$.
On every $K\Subset\mathcal R(X)$,
\begin{equation}\label{eq:smoothconv}
 \Phi_j^*g_j\longrightarrow h\quad\hbox{in }C^\infty(K\times I),
 \qquad
 \Phi_j^*J_j\longrightarrow J_\infty\quad\hbox{in }C^\infty(K).
\end{equation}
The first convergence includes all mixed space and time derivatives.
\end{enumerate}
The data are chosen along one sequence, by exhaustion in space and
time. When all spaces are compact, we use full correspondences,
surjective onto both spaces, and omit basepoints if none are specified.
\end{definition}

The pointed formulation compares balls at every time. The mutual ball
containment in Theorem~\ref{thm:selected} verifies this requirement for
its noncompact limits. On a single time
slice the definition reduces to Definition~\ref{def:modelhat}:
selecting a point in each fiber of a relation gives a GH
approximation, and on a regular compact set the selection is
$\Phi_j$. Conversely, the graph of a GH approximation can be
enlarged by pairs at distance tending to zero to give a relation.
Exhaustion gives the formulation above.

The spacetime embeddings are $(x,\tau)\mapsto(\Phi_j(x),\tau)$;
in particular they preserve $\partial_\tau$. No metric on the
product spacetime is being introduced. A fixed singular set in this
definition is justified in the proof of Theorem~\ref{thm:main}; the
definition itself does not prove that singularities persist.

\begin{lemma}[Regular neighborhoods and moving time slices]
\label{lem:spacetimecoverage}
In Definition~\ref{def:topology}, let
$K\Subset U\Subset\mathcal R(X)$, with $U$ a smooth relatively
compact domain. For all sufficiently large $j$,
\begin{equation}\label{eq:spacetimecoverage}
 \{x\in X_j:(x,a)\in\mathscr C_j\text{ for some }a\in K\}
 \subset\Phi_j(U).
\end{equation}
This holds uniformly on $I$. If $\tau_j\to\tau\in I$, the time
slices converge in pointed-$\hat C^\infty$-Cheeger--Gromov topology
to $(X,p,h(\tau),J_\infty)$, using the same regular embeddings.
\end{lemma}
\begin{proof}
Choose a smooth domain $U'$ with $K\Subset U'\Subset U$.
On $\overline{U'}\times I$, smooth convergence bounds
$\Phi_j^*g_j(\tau)$ below by a fixed positive multiple of one
background metric. Consequently a curve starting in $\Phi_j(K)$
and first leaving $\Phi_j(U')$ has length at least $b>0$,
independently of $j$ and $\tau$. The embedding is defined past
$\overline{U'}$, so a first exit meets $\Phi_j(\partial U')$.
If $(x,a)\in\mathscr C_j$ and $a\in K$, the graph condition and
\eqref{eq:distortion} give
$d_j(\tau)(x,\Phi_j(a))\le\epsilon_j<b$. Thus $x\in\Phi_j(U')$.

For moving times, \eqref{eq:distortion} compares the source slice
with $(X,d(\tau_j),p)$. Choose relatively compact domains
$V\Subset V^+\subset X$ containing the comparison region, with
$h$ smooth in orbifold charts on a neighborhood of $\overline{V^+}$.
Smoothness and completeness give pointed GH convergence
$(X,d(\tau_j),p)\to(X,d(\tau),p)$. Indeed, a sufficiently large
compact domain has boundary farther from $p$ than the ball under
consideration; this remains true for nearby times, by the same
first-exit argument. Distances on the resulting compact comparison
region converge uniformly. Equation~\eqref{eq:smoothconv} then
gives smooth convergence on its regular part at the moving times.
\end{proof}

\begin{lemma}[Extension from a dense set of times]
\label{lem:densetimes}
Fix a compact interval $I$ and a dense subset $D\subset I$.
On a fixed spatial comparison region, suppose time-independent
relations $\mathscr C_j$ satisfy
$\operatorname{dis}_\tau(\mathscr C_j)\to0$ for every $\tau\in D$.
Suppose also that both distance families have a common modulus of
continuity in time on the projections of these relations:
\begin{align*}
 |d_j(s)(x,y)-d_j(t)(x,y)|&\le\omega(|s-t|),\\
 |d(s)(a,b)-d(t)(a,b)|&\le\omega(|s-t|),
 \qquad \omega(\delta)\longrightarrow0.
\end{align*}
Then their distortions tend to zero uniformly on $I$.
The same conclusion holds with an additional error $\eta_j\to0$
in the first inequality.

Suppose, in addition, that fixed regular embeddings compatible with
these relations give the prescribed smooth limits at times in $D$.
For each $U\Subset U^+\Subset\mathcal R(X)$, assume their
pullback metrics have bounds, uniform in $j$, for every mixed
derivative on $U^+\times I$. Then the smooth convergence is
uniform on $\overline U\times I$.
\end{lemma}
\begin{proof}
Choose a finite $\delta$-net $\{t_1,\ldots,t_N\}\subset D$ of $I$.
The triangle inequality gives
\[
 \sup_{t\in I}\operatorname{dis}_t(\mathscr C_j)
 \le \max_{1\le k\le N}
       \operatorname{dis}_{t_k}(\mathscr C_j)
       +2\omega(\delta)+\eta_j.
\]
First let $j\to\infty$, then $\delta\to0$.
For the smooth statement, the derivative bounds give precompactness
in every finite differentiability order on smaller domains.
Every subsequential limit agrees with the prescribed family on
$D$, hence everywhere by continuity. Uniqueness of this limit
gives smooth convergence of the sequence.
\end{proof}

The criterion requires common point identifications and fixed
regular charts; separate slice limits do not supply them.
For example, a nonconstant family $F_\tau^*h$ is slice-wise
isometric to $h$ but cannot be identified with it by one fixed
isometry. Below we instead extract smooth limits on whole compact
time intervals in fixed charts, then prove uniform GH convergence
by identifying fixed points also at moving times.

The metric-space convergence also determines a topology. For compact spaces with jointly continuous distance
families on $I$, set
\begin{equation}\label{eq:familyGH}
 \delta_I(X,Y):=\frac12\inf_{\mathscr C}
 \sup_{\tau\in I}\operatorname{dis}_\tau(\mathscr C),
\end{equation}
where $\mathscr C$ ranges over full correspondences independent of
time. This is the correspondence description of GH distance
\cite{BBI}, with a single correspondence for the whole interval.

\begin{lemma}\label{lem:familymetric}
The function $\delta_I$ is a metric on compact metric families
modulo a single time-independent isometry preserving all distances.
Definition~\ref{def:topology} adds compatible smooth convergence
on the regular part.
\end{lemma}
\begin{proof}
Reversing and composing correspondences proves symmetry and the
triangle inequality. If $\delta_I(X,Y)=0$, close correspondences
of distortion tending to zero in the compact space $X\times Y$,
using reference-time metrics. A Hausdorff-convergent subsequence
has a limit projecting onto both factors. Joint continuity makes
its distortion zero at every time. At the reference time it is
the graph of a bijection; at every other time this same bijection
is an isometry. Smooth neighborhoods are specified by finite
derivative bounds on finitely many pairs of regular domains
$U\Subset U^+$. The embeddings are defined on $U^+$ and are
tested on $\overline U$. Restricting to an intermediate domain
before composing the embeddings preserves a positive spatial
margin and composes these tests. For pointed families one tests on
bounded balls as in Definition~\ref{def:topology}.
\end{proof}

Time-independent pullback commutes with $\partial_\tau$, so the
Calabi equation passes to the limit under \eqref{eq:smoothconv}
without an additional Lie derivative. In the proof below,
extremality and Lemma~\ref{lem:transport} then identify the limit
by \eqref{eq:mainfamily} in the original comparison maps.

Volume convergence is a separate conclusion of the geometric
estimates, as in Lemma~\ref{lem:volume-passage}.
\begin{definition}\label{volume measures converge}
In the compact case, weak convergence of volume measures uniformly on
compact time intervals means the following. For every compact interval $I$,
there are Borel GH approximation maps $f_j:X_j\to X$. They are independent
of $\tau$, compatible with Definition~\ref{def:topology}, and satisfy,
for every $\varphi\in C(X)$,
\[
 \lim_{j\to\infty}\sup_{\tau\in I}
 \left|\int_{X_j}\varphi\circ f_j\,dV_{g_j(\tau)}
       -\int_X\varphi\,dV_{h(\tau)}\right|=0.
\]
\end{definition}

\subsection{Proof of Theorem~\ref{thm:main}\textup{(I)--(III)}}

\begin{proof}
Fix $t_j\to\infty$ and define
$ g_j(\tau):=g(t_j+\tau), \tau\in[-t_j,\infty).$

\smallskip
\noindent\emph{Step 1: Estimates on bounded translated intervals.}
Lemma~\ref{lem:energy} gives
bounded curvature energy.
We assume the Sobolev inequality \eqref{eq:Sob} holds with the
constants $A,V$, and
Lemma~\ref{lem:noncollapse} gives the same volume
lower bounds and diameter bound at every translated
time.
For each fixed $T\ge1$, the translated flows are
defined on $[-3T,3T]$ for all sufficiently large $j$.
Lemma~\ref{lem:tail} provides vanishing dissipation on these intervals.

We apply Theorem~\ref{thm:radiuslower} to the original
flow at times with $t_0=1$ to get $r_{\rm can}>0$ such that
$\scr(x,t)\ge r_{\rm can}$, $x\in M,\ t\ge1$.
Set $r_0:=\min\{1,r_{\rm can}/2\}$. 
Then, for every fixed $T\ge1$ and large $j$, we have a uniform lower bound for the spacetime canonical radius
\[
 \inf_{M\times[-2T,2T]}\scr_j\ge r_0.
\]
In particular, the hypotheses
\eqref{eq:segmentbounds} required in
Theorem~\ref{thm:upgrade} hold on these translated
segments. 

\smallskip
\noindent\emph{Step 2: Choice of the time-zero comparison maps.}
Corollary~\ref{cor:geometricinputs}\textup{(B1)},
applied at time zero, gives a subsequence converging
in the GH sense to a compact K\"ahler orbifold
$(X,h_0,J_\infty)$. The metric convergence is smooth
on $\mathcal R(X)$, and $h_0$ extends smoothly in
orbifold charts across the finite set
$
 \mathcal S:=\operatorname{Sing}X.
$
Choose increasing open domains
$
 U_1\Subset U_2\Subset\cdots\Subset\mathcal R(X),
\bigcup_jU_j=\mathcal R(X)$.
After a diagonal subsequence, choose comparison
embeddings $\Phi_j:U_j\to M$ with
 $\Phi_j^*g_j(0)\longrightarrow h_0$,
 $\Phi_j^*J\longrightarrow J_\infty$
smoothly on compact regular subsets.
We choose these maps compatibly with the GH
comparisons, in the sense of
Definition~\ref{def:modelhat}.
Using the relation/map conversion explained after
Definition~\ref{def:topology}, choose full
correspondences
$
 \mathscr C_j^0\subset M\times X
$
whose two projections are surjective and which
satisfy
\[
 \operatorname{dis}_0(\mathscr C_j^0)\longrightarrow0,
 \qquad
 (\Phi_j(z),z)\in\mathscr C_j^0,
 \quad z\in U_j.
\]
Each $\Phi_j$ is used at every
time and does not depend on $\tau$.

\smallskip
\noindent\emph{Step 3: Smooth convergence in the retained embeddings.}
Fix $T\ge1$ and nested regular domains
\[
 K\Subset V_0\Subset V_1\Subset V_2
 \Subset\mathcal R(X).
\]
For large $j$, $\Phi_j$ is defined on a neighborhood
of $\overline{V_2}$.
Terminal smooth convergence gives uniform curvature
bounds on $\Phi_j(V_2)$ and a uniform lower bound
for $\Phi_j^*g_j(0)$ by a positive multiple of $h_0$.

The first-exit length comparison used to prove
the capture condition \eqref{eq:scrcapture}
therefore gives a uniform positive distance between
$\Phi_j(\overline{V_1})$ and
$M\setminus\Phi_j(V_2)$.
Choose $r>0$ below this spatial margin and small
enough that the terminal curvature bound is at most
$r^{-2}$.
A finite cover of $\overline{V_0}$ then supplies
points $z_1,\ldots,z_L\in V_1$ such that
$
 \Phi_j(\overline{V_0})
 \subset
 \bigcup_{\ell=1}^L
 B_{g_j(0)}(\Phi_j(z_\ell),r),
$
while every concentric $2r$-ball lies in
$\Phi_j(V_2)$ and satisfies
\[
 \sup_{B_{g_j(0)}(\Phi_j(z_\ell),2r)}
 |\Rm_{g_j(0)}|\le r^{-2}.
\]

The hypotheses of Theorem~\ref{thm:upgrade} have
been verified in Step~1, and
$$
 \Ca(t_j-3T)-\Ca(t_j+3T)\longrightarrow0.
$$
Apply that theorem at the finitely many centers
$\Phi_j(z_\ell)$. It gives bounds for every
intrinsic curvature derivative on the fixed set
$\Phi_j(\overline{V_0})$ throughout $[-T,T]$,
uniformly in $j$.

Starting from these intrinsic bounds, apply the
metric-variation and connection-variation estimates in
Step~3 of the proof of
Lemma~\ref{lem:forwardpersistence}.
The order-two curvature bound controls the metric
velocity, and integration gives
$
 c_{K,T}h_0\le\Phi_j^*g_j(\tau)\le C_{K,T}h_0.
$
The connection estimates, the terminal smooth
bounds, and the Calabi flow equation give
\[
 \sup_{K\times[-T,T]}
 \left|
 (\nabla^{h_0})^a\partial_\tau^c
       (\Phi_j^*g_j)
 \right|_{h_0}
 \le C_{K,T,a,c},
 \qquad a,c\ge0.
\]
Here the pulled-back complex structures
$\Phi_j^*J$ are independent of time.

Diagonal extraction over the regular exhaustion
and $T=1,2,\ldots$ gives
\[
 \Phi_j^*g_j(\tau)\longrightarrow h(\tau)
 \quad\text{in }C^\infty_{\mathrm{loc}}
       (\mathcal R(X)\times\R),
 \qquad h(0)=h_0.
\]
Lemma~\ref{lem:parallel-complex} gives the K\"ahler
limit with complex structure $J_\infty$.
The time independence of $\Phi_j$ allows the
Calabi flow equation to pass to the limit.

Apply Proposition~\ref{prop:local} on regular
backward cylinders, shifting each chosen terminal
time to zero. The vanishing dissipation in Step~1
shows that every $h(\tau)$ is extremal.

By Corollary~\ref{cor:geometricinputs}\textup{(B1)},
$h_0$ extends smoothly to the compact orbifold $X$.
Its extremality extends across $\mathcal S$ by
continuity. Lemma~\ref{lem:regularorbits} therefore
gives the complete real holomorphic field
$
 W_0:=\frac12\nabla_{h_0}S_{h_0}
$
and its flow $F_\tau$. This flow preserves
$\mathcal R(X)$ and fixes every point of $\mathcal S$.
The global case of Lemma~\ref{lem:transport} gives
\[
 h(\tau)=F_\tau^*h_0
 \qquad\text{on }\mathcal R(X),\quad \tau\in\R.
\]
This formula extends the family smoothly to $X$.

The quotient-link argument at the beginning of
Section~\ref{sec:extremal-limits} shows that every
minimizing $h_0$-geodesic with regular endpoints
stays in $\mathcal R(X)$.
Since $F_\tau:(X,h(\tau))\to(X,h_0)$ is an
isometry preserving the regular locus, the same
statement holds for every $h(\tau)$.
So, we have geodesic convexity of
$\mathcal R$ for every $h(\tau)$.

\smallskip
\noindent\emph{Step 4: Uniform metric comparisons using the same pairs.}
\begin{claim}
Fix $T\ge1$. With distortion defined by
\eqref{eq:distortion}, we prove
\begin{equation}\label{eq:fixedrelationuniform}
 \sup_{|\tau|\le T}
 \operatorname{dis}_\tau(\mathscr C_j^0)
 \longrightarrow0.
\end{equation}
\end{claim}
\begin{proof}
Suppose \eqref{eq:fixedrelationuniform} fails.
After passing to a subsequence, choose
$\varepsilon>0$, 
$\tau_j\to\tau\in[-T,T]$, and pairs
\[
 (x_j,a_j),(y_j,b_j)\in\mathscr C_j^0,
 \qquad a_j\to a,\quad b_j\to b,
\]
such that
$
 \left|
 d_{g_j(\tau_j)}(x_j,y_j)
 -d_{h(\tau_j)}(a_j,b_j)
 \right|\ge\varepsilon.
$

Use the construction established in Steps~1--4
of the proof of
Lemma~\ref{lem:forwardpersistence} on
$
 [-3T/2,3T/2]\Subset(-2T,2T).
$
That construction uses the three properties in Corollary~\ref{cor:geometricinputs}\textup{(B1)--(B3)} and vanishing dissipation.

Starting from the later terminal time $3T/2$,
the construction gives, after further extraction,
a compact complex orbifold $Y$ and an extremal
family $k(\tau)$ on this interval.
Under its identifications, each sequence $x_j$ or $y_j$ of points
of the underlying manifolds, held fixed as time varies, has
a time-independent limiting point in $Y$, including at moving times.
Its regular comparison embeddings are independent
of time.

We identify this family with the family on $X$.
Let $\mathscr D_j^0\subset M\times Y$ be compatible
full correspondences for its time-zero limit.
Relate $v\in Y$ to $c\in X$ whenever some $x\in M$
satisfies
$
 (x,v)\in\mathscr D_j^0,
 (x,c)\in\mathscr C_j^0.
$
These composed relations are full and have
distortions tending to zero at time zero.
The compact-correspondence argument in the proof
of Lemma~\ref{lem:familymetric}, applied at the
single time zero, gives an isometry
\[
 I_0:(Y,k(0))\longrightarrow(X,h_0).
\]

We check that this identification is holomorphic.
The isometry preserves the regular loci, because
regularity is characterized by Euclidean tangent
cones, as explained in
Section~\ref{sec:extremal-limits}.
Choose a connected regular exterior
$U\subset\mathcal R(Y)$ as in
Lemma~\ref{lem:compactextension}; if $Y$ is smooth,
take $U=Y$.
The containment property in
Corollary~\ref{cor:geometricinputs}\textup{(B1)}
makes the transition maps between the two systems
of regular embeddings well-defined near
$\overline U$.
The transition-map estimates in the proof of
Lemma~\ref{lem:retainmaps} give smooth convergence
of these maps to $I_0|_U$.
Passing to the limit in the complex-structure
identities shows that $I_0|_U$ is holomorphic.

Both time-zero metrics are extremal and hence
analytic by Lemma~\ref{lem:analytic}.
Both have volume $V$ by
Corollary~\ref{cor:geometricinputs}\textup{(B1)}.
Lemma~\ref{lem:compactextension} therefore extends
$I_0|_U$ to a holomorphic orbifold isometry
\[
 I:(Y,k(0),J_Y)\longrightarrow(X,h_0,J_\infty).
\]
This extension equals $I_0$ everywhere:
the two isometries agree on $U$, so their images
of any point have identical distances to $I_0(U)$.
The distance-to-regular-points argument in Step~2
of the proof of
Lemma~\ref{lem:forwardpersistence} shows that those
distances determine the point uniquely.

By \eqref{eq:Wnaturality}, $I$ identifies the
extremal vector fields of $k(0)$ and $h_0$.
Lemma~\ref{lem:transport} consequently gives
\[
 k(\tau)=I^*h(\tau)
 \qquad\text{for }\tau\in[-3T/2,3T/2].
\]
Thus one time-independent map identifies the
entire auxiliary family with the family on $X$.

Choose partners $v_j,w_j\in Y$ such that $ (x_j,v_j),(y_j,w_j)\in\mathscr D_j^0.$
After extraction, $v_j\to v$ and $w_j\to w$.
The composed time-zero relations give
$
 I(v)=a, I(w)=b.
$
The identification of points of the underlying manifolds held fixed
as time varies, proved in Lemma~\ref{lem:forwardpersistence}
also for moving times, now implies
\[
 d_{g_j(\tau_j)}(x_j,y_j)
 \longrightarrow d_{k(\tau)}(v,w)
 =d_{h(\tau)}(a,b).
\]
The smooth orbifold family $h(\tau)$ has jointly
continuous distances, so
\[
 d_{h(\tau_j)}(a_j,b_j)
 \longrightarrow d_{h(\tau)}(a,b).
\]
This contradicts the chosen positive distortion
and proves \eqref{eq:fixedrelationuniform}.
\end{proof}
The correspondences $\mathscr C_j^0$ remain full
and contain the original graph pairs
$(\Phi_j(z),z)$. Together with Step~3, this verifies
Definition~\ref{def:topology} on every compact
time interval.
The singular set is $\mathcal S$ for every
limiting metric, and each $F_\tau$ fixes it.
Under these fixed comparisons, the singular set
of the limiting spacetime is therefore
$\mathcal S\times\R$.

\smallskip
\noindent\emph{Step 5: Uniform weak convergence of volume measures.}
We prove the conclusion in
Definition~\ref{volume measures converge}. The proof is similar to that of Lemma~\ref{lem:volume-passage}.
Close $\mathscr C_j^0$ in $M\times X$, retaining
the same notation. This preserves its projections,
graph pairs, and distance estimates.

For each $x\in M$, the set
$
 \mathscr C_j^0(x):=
 \{z\in X:(x,z)\in\mathscr C_j^0\}
$
is nonempty and compact. We can choose a Borel
map $f_j^0:M\to X$ with
$f_j^0(x)\in\mathscr C_j^0(x)$.
The image $A_j:=\Phi_j(U_j)$ is open, and
$\Phi_j^{-1}$ is continuous there. Define
\[
 f_j(x):=
 \begin{cases}
 \Phi_j^{-1}(x),&x\in A_j,\\
 f_j^0(x),&x\notin A_j.
 \end{cases}
\]
This is Borel and independent of time.
Its graph remains in $\mathscr C_j^0$, and
$
 f_j\circ\Phi_j=\mathrm{id}_{U_j}.
$
The fullness and vanishing uniform distortion of
the correspondences make these maps GH
approximations uniformly on each compact time
interval. Indeed, if $(x,z)\in\mathscr C_j^0$,
comparison with $(x,f_j(x))$ gives
$
 d_{h(\tau)}(z,f_j(x))
 \le\operatorname{dis}_\tau(\mathscr C_j^0).
$

Fix $I:=[-T,T]$ and put
$
 \epsilon_j:=
 \sup_{\tau\in I}
 \operatorname{dis}_\tau(\mathscr C_j^0)
 \longrightarrow0.
$
We control neighborhoods of the
singular points. Write
$\mathcal S=\{\sigma_1,\ldots,\sigma_\ell\}$ and
let $N_r(\mathcal S)$ be its $h_0$-neighborhood.
Smoothness of the orbifold family on $X\times I$
gives $A_T\ge1$ such that
$
 d_{h(\tau)}(z,z')\le A_Td_{h_0}(z,z'),
\tau\in I.
$
For each $\alpha$, choose $y_{j,\alpha}\in M$
with
$(y_{j,\alpha},\sigma_\alpha)\in\mathscr C_j^0$.
If $f_j(x)\in N_r(\sigma_\alpha)$, then
\[
 d_{g_j(\tau)}(x,y_{j,\alpha})
 \le d_{h(\tau)}(f_j(x),\sigma_\alpha)+\epsilon_j
 <A_Tr+\epsilon_j.
\]
Hence
\[
 f_j^{-1}(N_r(\mathcal S))
 \subset
 \bigcup_{\alpha=1}^{\ell}
 B_{g_j(\tau)}(y_{j,\alpha},A_Tr+\epsilon_j).
\]
By Theorem~\ref{thm:radiuslower},
test~\eqref{eq:scrV} is available uniformly at
these translated times:
\[
 \Vol_{g_j(\tau)}B_{g_j(\tau)}(y,\varrho)
 \le B_0\varrho^4,
 \qquad 0<\varrho<r_0,
\]
for every $y\in M$, $\tau\in I$, and large $j$.
Fix $r>0$ sufficiently small that
$(A_T+1)r<r_0$, and then take $j$ large enough
that $\epsilon_j<r$. The preceding inclusion gives
\[
 \sup_{\tau\in I}
 \Vol_{g_j(\tau)}
 f_j^{-1}(N_r(\mathcal S))
 \le C_Tr^4.
\]
The quotient charts described in
Section~\ref{sec:extremal-limits}, together with
uniform metric equivalence on $I$, imply
\[
 \sup_{\tau\in I}
 \Vol_{h(\tau)}N_r(\mathcal S)\le C_Tr^4.
\]

Choose a continuous function $\chi_r$ with
$0\le\chi_r\le1$, equal to zero on
$N_r(\mathcal S)$ and equal to one outside
$N_{2r}(\mathcal S)$.
For $\varphi\in C(X)$, the function
$\chi_r\varphi$ has compact regular support,
so the convergence follows by using Lemma~\ref{lem:spacetimecoverage}, specifically
\eqref{eq:spacetimecoverage}.
As in the proof of
Lemma~\ref{lem:volume-passage}, the omitted
integrals are bounded on both spaces by
$C_T\|\varphi\|_\infty r^4$. Therefore
\[
 \limsup_{j\to\infty}\sup_{\tau\in I}
 \left|
 \int_M\varphi\circ f_j\,dV_{g_j(\tau)}
 -\int_X\varphi\,dV_{h(\tau)}
 \right|
 \le C_T\|\varphi\|_\infty r^4\longrightarrow 0, \quad r\downarrow0.
\]
This proves the uniform weak
convergence required in
Definition~\ref{volume measures converge}.
Taking $\varphi\equiv1$ gives
$
 \Vol(X,h(\tau))=V.
$

\smallskip
\noindent\emph{Step 6: Total scalar curvature.}
The volume and average complex scalar curvature
of the source metrics are fixed by $\Omega$.
By \eqref{eq:energy} and \eqref{eq:dissipation},
\[
 \int_M S_{g_j(\tau)}^2\,dV_{g_j(\tau)}
 =\Ca(t_j+\tau)+V\bar S^2
 \le\Ca(0)+V\bar S^2.
\]
Cauchy--Schwarz and the exceptional-volume bound
from Step~5 give
\[
 \sup_{\tau\in I}
 \int_{f_j^{-1}(N_{2r}(\mathcal S))}
 |S_{g_j(\tau)}|\,dV_{g_j(\tau)}
 \le C_Tr^2.
\]
The limiting scalar curvatures are uniformly
bounded on the smooth compact orbifold family
$X\times I$, so their integrals over
$N_{2r}(\mathcal S)$ also tend uniformly to zero.

On the support of $\chi_r$, smooth convergence
of the metrics gives smooth convergence of
scalar curvature. Using the same change of
variables as in Step~5,
\[
 \sup_{\tau\in I}
 \left|
 \int_M(\chi_r\circ f_j)S_{g_j(\tau)}
                   \,dV_{g_j(\tau)}
 -\int_X\chi_r S_{h(\tau)}\,dV_{h(\tau)}
 \right|\longrightarrow0.
\]
Since
$
 \int_M S_{g_j(\tau)}\,dV_{g_j(\tau)}=V\bar S,
$
letting $j\to\infty$ and then $r\downarrow0$
gives
$
 \int_X S_{h(\tau)}\,dV_{h(\tau)}=V\bar S.
$
The interval $I$ was arbitrary, so these conclusions
hold for every $\tau\in\R$.

This proves \textup{(I)--(III)}.
Part \textup{(IV)} is proved in
Subsection~\ref{sec:calabicontinuity}.
\end{proof}

\begin{corollary}\label{cor:noaut}
If the limiting complex orbifold has no nonzero
holomorphic vector field, then $h(0)$ has constant
scalar curvature and $h(\tau)$ is stationary.
\end{corollary}
\begin{proof}
Since $h_0$ is extremal,
$\nabla^{1,0}_{h_0}S_{h_0}$ is an orbifold
holomorphic vector field. It therefore vanishes.
Thus $dS_{h_0}=0$ on the connected regular locus,
and smooth orbifold extension makes $S_{h_0}$
constant on $X$. Consequently $W_0=0$,
$F_\tau=\mathrm{id}_X$, and
$h(\tau)=F_\tau^*h_0=h_0$.
\end{proof}

\subsection{Continuity of the Calabi functional and proof of Theorem~\ref{thm:main}\textup{(IV)}}
\label{sec:calabicontinuity}

We prove Theorem~\ref{thm:main}\textup{(IV)} by excluding scalar-energy
concentration at times of small dissipation, then using monotonicity.
The shrinking-scale limits are scalar-flat; the following estimate
controls the remaining conical regions.

\begin{lemma}[A length-independent scalar estimate on a conical neck]
\label{lem:weightedDneck}
Let $e$ and $J_e$ be the flat K\"ahler metric and complex structure on
$C_\Gamma:=(\C^2\setminus\{0\})/\Gamma$, with radius $r$.
There exist constants $C>0$ and $\delta_*>0$ with the following
properties.
\begin{enumerate}[label=\textup{(\roman*)}]
\item \emph{Flat-cone estimate.}
Every real function $f\in C_c^\infty(C_\Gamma)$ satisfies
\begin{equation}\label{eq:flatweightedD}
 \int_{C_\Gamma}f^2\,dV_e
 \le C\int_{C_\Gamma}r^4|\Dcal_e f|^2\,dV_e.
\end{equation}

\item \emph{Stability under perturbation.}
Let $(U,g,J)$ be a K\"ahler annulus admitting a cone chart in which
\begin{equation}\label{eq:conegauge}
 \sum_{k=0}^2 r^k|\nabla_e^k(g-e)|_e
 +\sum_{k=0}^2 r^k|\nabla_e^k(J-J_e)|_e
 <\delta_*.
\end{equation}
Then every real function $f\in C_c^\infty(U)$ satisfies
\[
 \int_U f^2\,dV_g
 \le C\int_U r^4|\Dcal_g f|^2\,dV_g.
\]

\item \emph{Scalar-curvature estimate on a neck.}
Let $A\Subset\widetilde A$ be concentric necks, with
$\widetilde A$ satisfying the hypotheses of \textup{(ii)}.
Assume that each boundary component of $A$ is separated from
the corresponding boundary component of $\widetilde A$
by a fixed radius ratio.

Choose a smooth cutoff $\chi$ equal to one on $A$ and supported
compactly in $\widetilde A$, whose derivatives are supported in
transition collars
\[
 C_\partial\subset\widetilde A\setminus A
\]
of fixed radius ratio.
Then the scalar curvature $S:=S_g$ satisfies
\begin{equation}\label{eq:scalarneckestimate}
 \int_A S^2\,dV_g
 \le C r_+^4\int_{\widetilde A}|\Dcal_g S|^2\,dV_g
 +C\int_{C_\partial}
       \bigl(S^2+r^2|\nabla S|^2\bigr)\,dV_g,
\end{equation}
where $r_+$ is the outer radius of $\widetilde A$.
\end{enumerate}
The constants are independent of the ratio of the outer and
inner radii of the annuli; in \textup{(iii)}, they may depend
on the fixed collar ratios.
\end{lemma}
\begin{proof}
We lift to $\C^2\setminus\{0\}$ and work with $\Gamma$-invariant functions.
We put $t:=\log r$, $u:=r^2f$ and
\[
 Z:=\sum_{a=1}^2z^a\partial_{z^a},\qquad
 T:=i(Z-\overline Z).
\]
The field $T$ generates the scalar circle action on $S^3$, which
commutes with $\Gamma$. We have
$Z=(\partial_t-iT)/2$, $\overline Z=(\partial_t+iT)/2$ and
\[
 \int f^2\,dV_e=\int|u|^2\,dt\,d\sigma,
 \qquad
 (Z^2-Z)f=r^{-2}(Z-1)(Z-2)u.
\]
The analogous identity holds for $\overline Z$. We decompose into the
circle modes $T=im$, $m\in\Z$, and take the Fourier transform in $t$.
The two multipliers are
\[
 (a_++ib)(a_+-1+ib),\qquad
 (a_-+ib)(a_--1+ib),
 \quad a_\pm:=\pm m/2-1,\quad b:=\xi/2.
\]
Since $a_++a_-=-2$, one is at most $-1$, and the modulus of its
multiplier is at least two. Plancherel therefore bounds $\|u\|_2^2$
by the sum of the two squared output norms. But
\[
 (Z^2-Z)f=\sum_{a,b}z^az^b f_{ab},\qquad
 (\overline Z^2-\overline Z)f
       =\sum_{a,b}\overline z^a\overline z^b f_{\overline a\overline b}.
\]
For real $f$, these are conjugate contractions of its pure complex
Hessians, bounded by $Cr^2|\Dcal_e f|$. This proves
\eqref{eq:flatweightedD}. Dividing the lifted integrals by $|\Gamma|$
gives the quotient statement.

We next prove stability with constants independent of neck length.
Scaled interior estimates for the elliptic real operator $\Dcal_e$
on overlapping dyadic annuli, summed with bounded overlap, give
\[
 \int(r^4|\nabla_e^2f|^2+r^2|\nabla_e f|^2)
 \le C\left(\int r^4|\Dcal_e f|^2+\int f^2\right).
\]
The principal symbol is injective on real functions. The same estimate
holds on the quotient by lifting. Under \eqref{eq:conegauge}, the
coefficient difference satisfies
\[
 |(\Dcal_g-\Dcal_e)f|
 \le C\delta_*(|\nabla_e^2 f|+r^{-1}|\nabla_e f|).
\]
The target tensor bundles are identified by the natural projections
between the nearby complex structures. The metric norms and volume
forms are uniformly equivalent.
Combining the last two estimates with \eqref{eq:flatweightedD} and
choosing $\delta_*$ small absorbs the error. For functions supported
in an annulus all flat estimates apply after extension by zero. No
boundary condition on the metric is required.

We choose $\chi=1$ on $A$, supported in $\widetilde A$, with
$|\nabla\chi|\le C/r$ and $|\nabla^2\chi|\le C/r^2$ on
$C_\partial$. We apply the perturbed inequality to $f=\chi S$ and use
the product rule for $\Dcal_g$. Its cutoff terms are bounded in the
weighted square integral by
$C\int_{C_\partial}(S^2+r^2|\nabla S|^2)$, while $r^4\le r_+^4$.
This gives \eqref{eq:scalarneckestimate}.
\end{proof}
\begin{lemma}[Scalar estimate from a cone atlas]
\label{lem:atlasDneck}
Let $A\Subset\widetilde A$ be K\"ahler necks, with a cutoff
$\chi$ and transition collars $C_\partial$ satisfying the
boundary-separation and cutoff assumptions of
Lemma~\ref{lem:weightedDneck}\textup{(iii)}.

Choose a sufficiently large logarithmic width $L$, and then a
sufficiently small tolerance $0<\delta(L)\le\delta_*$.
Suppose that $\widetilde A$ carries a positive continuous radius function $r$,
not necessarily smooth, and a \emph{cone atlas} satisfying:
\begin{enumerate}[label=\textup{(\roman*)}]
\item \emph{Chart bands and coverage.}
The $j$th chart is defined on
$
 \{\,|\log r-jL|<2L+2\,\}
$
and is used on the inner band
$
 \{\,|\log r-jL|\le L\,\}.
$
These inner bands cover $\operatorname{supp}\chi$.

\item \emph{Bounded overlap.}
The larger chart bands have bounded overlap, independently of
$L$ and of the number of charts.

\item \emph{Perturbation bounds.}
The cone radius $\rho_j$ in the $j$th chart satisfies
$
 \frac12r\le\rho_j\le2r.
$
In that chart, the metric and complex structure satisfy
\eqref{eq:conegauge} from
Lemma~\ref{lem:weightedDneck}\textup{(ii)}, with $r$ replaced by
$\rho_j$ and $\delta_*$ replaced by $\delta(L)$.
\end{enumerate}
Then the scalar-curvature estimate \eqref{eq:scalarneckestimate}
of Lemma~\ref{lem:weightedDneck}\textup{(iii)} holds, where $r_+$
is the outer radius of $\widetilde A$ with respect to $r$.
The constant is independent of the ratio of the outer and inner
radii of $\widetilde A$.
\end{lemma}
\begin{proof}
Fix a tolerance ceiling $\bar\delta\le\delta_*$ sufficiently small
that the perturbed estimate and the local elliptic estimates
from Lemma~\ref{lem:weightedDneck} have uniform constants.
We shall choose $L$ below and then take
$\delta(L)\le\bar\delta$.

Let $f\in C_c^\infty(\widetilde A)$ be real, with
$\operatorname{supp}f\subset\operatorname{supp}\chi$.
We first establish the weighted estimate for $f$.
Choose $\eta\in C_c^\infty((-2,2))$ with
$0\le\eta\le1$ and $\eta=1$ on $[-3/2,3/2]$.
Writing $\rho_j$ for the cone radius in the $j$th chart, set
\[
 \phi_j
 :=\eta\!\left(\frac{\log\rho_j-jL}{L}\right).
\]
Since $|\log\rho_j-\log r|\le\log2$, for $L\ge2\log2$
the function $\phi_j$ equals one on the corresponding inner band.
Consequently,
$\sum_j\phi_j^2\ge1$
on $\operatorname{supp}f.$
Moreover, its support lies in
\[
 |\log r-jL|<2L+\log2<2L+2.
\]
Thus $\phi_jf$ is supported compactly inside its chart and
extends smoothly by zero.

The chart bounds and radius comparison give
$|\nabla\phi_j|\le\frac{C}{Lr}$ and
$|\nabla^2\phi_j|\le\frac{C}{Lr^2}$.
Here the second derivative contains both
$\eta''\,d\log\rho_j\otimes d\log\rho_j/L^2$
and $\eta'\nabla^2\log\rho_j/L$.

Apply Lemma~\ref{lem:weightedDneck}\textup{(ii)}
to $\phi_jf$ in each chart. The product rule gives
\[
 |\Dcal_g(\phi_jf)-\phi_j\Dcal_gf|
 \le C\bigl(
 |\nabla\phi_j|\,|\nabla f|
 +|\nabla^2\phi_j|\,|f|
 \bigr).
\]
Comparing $\rho_j$ with $r$ and summing with bounded overlap,
we obtain
\[
 F\le C_1E+\frac{C_2}{L^2}(G+F),
\]
where
\[
 F:=\int_{\widetilde A}f^2\,dV_g,\qquad
 G:=\int_{\widetilde A}r^2|\nabla f|^2\,dV_g,\qquad
 E:=\int_{\widetilde A}r^4|\Dcal_gf|^2\,dV_g.
\]
The constants $C_1,C_2$ are independent of $L$, the number of
charts, and the neck length.

We next control $G$. Subdivide each inner chart band into
annuli of fixed logarithmic width, retaining fixed larger
subannuli inside the chart domains.
The scaled interior elliptic estimate used in the proof of
Lemma~\ref{lem:weightedDneck} applies on these nested pairs.
Their normalized separation and coefficient bounds are uniform.
Their overlap is bounded within each chart, and the atlas
overlap controls multiplicity between charts.
Summing the local estimates and comparing the local radii
with $r$ therefore gives
\[
 G\le C_0(E+F),
\]
with $C_0$ independent of $L$ and neck length.

Substitution yields
\[
 F\le
 \left(C_1+\frac{C_2C_0}{L^2}\right)E
 +\frac{C_2(C_0+1)}{L^2}F.
\]
Choose $L$ sufficiently large that
$
 \frac{C_2(C_0+1)}{L^2}<\frac12.
$
Absorbing the last term proves
\[
 \int_{\widetilde A}f^2\,dV_g
 \le C\int_{\widetilde A}r^4|\Dcal_gf|^2\,dV_g,
\]
with a constant independent of the ratio of the end radii.

Finally, take $f=\chi S$, using the standard end cutoff from
the proof of Lemma~\ref{lem:weightedDneck}\textup{(iii)}, with
$|\nabla\chi|\le C/r$ and $|\nabla^2\chi|\le C/r^2.
$
The cutoff computation established there applies to the
weighted estimate just proved and gives
\eqref{eq:scalarneckestimate}.
Its error terms are supported in $C_\partial$, and
$r^4\le r_+^4$ bounds the dissipation term.
All constants depend only on the fixed atlas and collar bounds,
not on the neck length.
\end{proof}

Zhang--Zhou \cite[Proposition~3.6]{ZZ} prove a curvature-decay
estimate on necks for cscK metrics. Here Lemmas~\ref{lem:weightedDneck}
and~\ref{lem:atlasDneck} control scalar energy using dissipation
and end collars. We construct the required covering by regular
regions and necks below.

A sequence of \emph{late-time slices} means $(M,g(u_i),J)$
with $u_i\to\infty$.

\begin{lemma}[Regular cores and conical necks on late-time slices]
\label{lem:scalarcover}
Let $g_i:=g(u_i)$, where $u_i\to\infty$, and suppose the slices
converge to $(Y_0,k_0):=(X,h)$ as in
Theorem~\ref{thm:main}\textup{(I)--(III)}. After passing to a
subsequence, the following hold.
\begin{enumerate}[label=\textup{(\roman*)}]
\item \emph{Regular cores.}
There are finitely many nonflat, complete scalar-flat K\"ahler ALE
orbifolds $(Y_a,k_a,J_a)$ obtained at length scales $s_{i,a}\to0$,
$1\le a\le N$. Set $s_{i,0}=1$. For $0\le a\le N$, compact
domains $K_a(R)$ exhaust $\mathcal R(Y_a)$ as $R\to\infty$.
On neighborhoods of these domains the comparison embeddings give
smooth convergence of $s_{i,a}^{-2}g_i$ and $J$ to the corresponding
limit. Write $U_{i,a}(R)$ for the image of $K_a(R)$.

\item \emph{Necks.}
Fix a logarithmic chart width $L$ as in Lemma~\ref{lem:atlasDneck}.
For $R$ sufficiently large and then $i$ sufficiently large,
there are $m\le N$ pairs $A_{i,b}(R)\Subset\widetilde A_{i,b}(R)$
such that
\[
 M=\bigcup_{a=0}^N U_{i,a}(R)\;\cup\;\bigcup_{b=1}^m A_{i,b}(R).
\]
Each neck joins the end of a child to a puncture of its parent.
Its radius $r$ is the $g_i$-distance to the child center, and its
outer radius is uniformly bounded in $g_i$. Each larger neck admits
the cone atlas of Lemma~\ref{lem:atlasDneck} on the support of the
cutoff in \textup{(iii)}. The chart radii are comparable to $r$, the overlap
is bounded independently of the neck length, and the largest
scale-invariant $C^2$ chart error tends to zero in the iterated
limit $i\to\infty$, then $R\to\infty$.

\item \emph{End collars.}
Each neck has a smooth cutoff supported compactly in
$\widetilde A_{i,b}(R)$ and equal to one on $A_{i,b}(R)$. Uniformly in $i$ and $R$,
\[
 |\nabla\chi_{i,b}|\le C/r,\qquad |\nabla^2\chi_{i,b}|\le C/r^2.
\]
Its derivatives are supported in two fixed-ratio collars in the
adjacent cores, at radii comparable to $R$ in the child metric
and to $R^{-1}$ in the parent metric.
\end{enumerate}
If the translates $g(t_i+\tau)$ converge to $h(\tau)$ as in
Theorem~\ref{thm:main}\textup{(I)--(III)} and
$u_i\in[t_i,t_i+1]$ with $u_i-t_i\to\tau_*$, the base is
$(X,h(\tau_*))$.
\end{lemma}
\begin{proof}
We verify the inputs to the selection and covering argument of
Lemma~\ref{lem:modelhierarchy} for Calabi-flow slices.

For any $p_i\in M$ and $\ell_i\downarrow0$, consider
$\widehat g_i(\tau):=\ell_i^{-2}g(u_i+\ell_i^4\tau)$.
Its backward time interval has length $u_i/\ell_i^4\to\infty$.
Theorem~\ref{thm:radiuslower}, the energy bound and the scaling
identities give, on every fixed backward interval $[-L,0]$,
\[
 \widehat{\scr}_i\ge r_0/\ell_i\to\infty,\qquad
 \widehat\beta_i:=\ell_i^2V^{-1/2}\to0,\qquad
 \int|\Rm_{\widehat g_i}|^2\le E_0,
\]
\[
 \int_{-L}^0\widehat D_i(\tau)\,d\tau
 =\int_{u_i-L\ell_i^4}^{u_i}D(t)\,dt\longrightarrow0.
\]
With $\lambda_i:=\ell_i^{-2}$ and the identity maps
of $M$, the complex structure remains $J$, and
\[
 [\lambda_i^{-1}\omega_{\widehat g_i(\tau)}]
 =[\omega(u_i+\ell_i^4\tau)]=\Omega.
\]
Thus Theorem~\ref{thm:selected} applies at every shrinking scale,
also after any further subsequence. The limit is a complete
stationary null scalar-flat ALE orbifold, and its smooth regular
charts cover compact subsets of its regular part.

Every quotient point of the base or of a rescaled limit has a
nonflat limit at a smaller scale. Indeed, bounded curvature near
that point, together with noncollapse and clause \textup{(R)},
would give a positive injectivity radius and a smooth limit of
the smooth source manifolds. To extract a nonflat limit, maximize
$|\Rm|$ times the squared distance to the boundary in a shrinking
neighborhood of the quotient point. Choose the neighborhood to
shrink slowly enough that its outer collar already converges
smoothly and the weighted maximum tends to infinity. Curvature
normalization then gives
\[
 |\Rm|(q_i)=1,\qquad
 \sup_{B(q_i,L_i)}|\Rm|\le4,\qquad L_i\to\infty.
\]
On a sufficiently small fixed ball, the volume upper bound makes
the curvature energy less than the threshold in \textup{(R)}.
That clause and noncollapse give smooth convergence at $q_i$,
so the rescaled limit is nonflat.

The energy-gap argument in Step~4 of the proof of
Lemma~\ref{lem:modelhierarchy}, applied with Sobolev
constant $A_*$, gives $\varepsilon_{\rm gap}>0$
such that every nonflat limit has energy at least
$\varepsilon_{\rm gap}$.
Two selections represent the same limit when their scale ratio
stays bounded above and below and the distance between their
centers is bounded in those scales. For a finite inequivalent
family, extract limits of all relative scales and center
distances. Comparable scales at diverging normalized distance
have disjoint compact cores. At separated scales a bounded
smaller-scale center must approach a singular point of the
larger limit: regular coverage and bounded larger-scale curvature
exclude a nonflat limit at a regular point. Remove these finitely
many attachment points from the larger cores. The regular
curvature measures have no atoms, so the cores can be chosen
disjoint with limiting energy greater than
$\varepsilon_{\rm gap}/2$ each. Their source images each have
energy greater than $\varepsilon_{\rm gap}/3$ for large $i$.
Hence every such family has at most $3E_0/\varepsilon_{\rm gap}$
members. Choose a family of maximal cardinality, allowing further
subsequences in this choice.

Group its finitely many centers at their mutual distance scales
as in Lemma~\ref{lem:modelhierarchy}. Every nonflat limit at a
grouping scale is already represented, by maximality. A flat
complete model with the inherited homogeneous Sobolev inequality
is $\mathbb R^4/\Gamma$ and has at most one quotient point.
Regular coverage implies that every finer nonflat selection
attaches at a quotient point. Consequently a flat level cannot
separate two nonempty groups and can be omitted. Every quotient
point of a retained limit has a smaller nonflat selection by the
preceding point-selection argument. Finiteness therefore makes
the last retained limits smooth. All selections lie above
singular points of the base. If none are needed, the compact
smooth convergence map onto the connected source is surjective,
and there is no remainder.

Consider two consecutive retained scales $a_i\ll b_i$, centered
at the smaller selection. At an intermediate scale
$a_i\ll\rho_i\ll b_i$, the limit must be flat; a nonflat limit
would enlarge the maximal family. Nor can it have a quotient
point away from the center. Such a point would have a nonflat
smaller-scale selection, introducing an additional grouping
scale between these consecutive ones. Finally, the center
itself cannot be regular, since it contains the smaller
nonflat selection. Thus the intermediate limit is a flat
quotient centered at its vertex. By the regular-coverage
conclusion of Theorem~\ref{thm:selected}, comparison charts
defined on larger limit annuli contain every point of the
underlying manifold in each smaller distance annulus.

This conclusion is uniform over the intermediate scales. More
precisely, on $Ra_i<d(p_i,\cdot)<b_i/R$, the largest optimal
normalized chart error on any annulus of prescribed fixed
logarithmic width tends to zero in the iterated limit
$i\to\infty$, $R\to\infty$. Otherwise choose a bad annulus
along $R_j\to\infty$ and then $i_j$ large enough that
$R_j^2a_{i_j}/b_{i_j}\to0$. Its scale is intermediate, contrary
to the preceding covered convergence. The same conclusion
holds with any prescribed finite number of derivatives.
At the two ends, it follows from ALE decay and smooth orbifold
extension in the retained limits.

For each intervening region use the positive radius function
$r:=d(p_i,\cdot)$. Smoothness of this function is unnecessary:
the cutoffs in Lemma~\ref{lem:atlasDneck} use the smooth cone
radius in each local chart. Covered annular convergence makes
each cone radius comparable to $r$ by a factor at most two.
Choose chart centers $L$ apart in $\log r$, with the larger
bands and smaller bands prescribed in that lemma. For its
fixed $L$, the preceding uniform convergence supplies these
charts with arbitrarily small normalized errors. Their
multiplicity is bounded by a universal constant, independently
of the number of bands. The finite number of intervening
regions also bounds the overlap between different regions.

Choose $K_a(R)$ with smooth radial boundaries in their ALE
ends and puncture charts. For example, retain the end up to
radius $2R$ and remove puncture balls of radius $1/(2R)$.
The intervening covered regions overlap these cores in collars
of radii $[R,2R]$ at a smaller-scale end and
$[1/(2R),1/R]$ at a larger-scale puncture, measured in the
respective normalized metrics. The compact cores and slightly
larger collar charts are fixed before $i\to\infty$. Define
each end cutoff in its smooth radial chart, make it constant
one toward the intervening region and constant zero toward
the core, and extend by these constants. Whole-annulus
coverage and the radial $C^0$ comparison with distance make
the chart formula agree with the constant extensions on
overlapping distance bands. This gives the stated cutoffs and their
scale-invariant derivative bounds. The fixed logarithmic
width needed for the cone atlas is obtained by retaining
larger collar charts; for fixed $L$ this changes only the
fixed radial factors.

These sets cover the source. To verify this, take an allegedly
uncovered sequence and locate it in the base limit. A regular
limit point is covered by a base chart. Near a singular point,
compare its distance with the finitely many selected scales.
An intermediate position is covered by the annular charts;
a bounded regular position at a retained scale is covered by
that core. Otherwise descend to a quotient point at the next
retained scale. This process ends at a smooth last limit,
whose bounded region is covered. A further branch would
contradict maximality of the selected family. The number of
intervening regions is at most the number of nonflat selections,
and all outer radii are bounded by the original diameter.
Finally, for moving times, Definition~\ref{def:topology} gives
the same base at $\tau_*$ and the argument applies unchanged.
\end{proof}

\begin{proof}[Proof of Theorem~\ref{thm:main}\textup{(IV)}]
Write $\Ca_\infty:=\lim_{t\to\infty}\Ca(g(t))$. Choose
$u_i\in[t_i,t_i+1]$ such that
\[
 D(u_i)\le\int_{t_i}^{t_i+1}D(t)\,dt+1/i\longrightarrow0.
\]
After passing to a subsequence, $u_i-t_i\to\tau_*\in[0,1]$.
The moving slices $g_i:=g(u_i)$ have base limit $h(\tau_*)$,
by the spacetime convergence already proved. This metric is
isometric to $h(0)$.

Apply Lemma~\ref{lem:scalarcover} to these slices. Write $S_i$
for their complex scalar curvatures. Smooth convergence gives convergence of $\int S_i^2$ on base
cores. On shrinking-scale cores the limit is scalar-flat, and
\[
 S_{\lambda g}^2\,dV_{\lambda g}=S_g^2\,dV_g
 \qquad (\lambda>0)
\]
in real dimension four. Its scalar-energy integral therefore
tends to zero for each fixed core.

On each of the remaining conical regions, Lemma~\ref{lem:atlasDneck}
and its end cutoffs give
\[
 \int_{A_i(R)}S_i^2\,dV_{g_i}
 \le C r_+^4D(u_i)
   +C\int_{C_{\partial,i}(R)}
         (S_i^2+r^2|\nabla S_i|^2)\,dV_{g_i}.
\]
The constants are independent of the ratio of the end radii.
The number of regions and their outer radii $r_+$ are uniformly
bounded, so the sum of the first terms tends to zero.

For fixed $R$, every transition collar lies in a smoothly
converging regular domain. On a shrinking-scale collar its
limit is scalar-flat. Both terms in its integral are scale
invariant, and the local radius is comparable to the model
radius, so that integral tends to zero. On a base collar of
radius $\rho$, smooth orbifold extension bounds $S$ and
$\nabla S$; its contribution is $O(\rho^4)+O(\rho^6)$.
The base collars shrink to the finite singular set as
$R\to\infty$. Consequently
\[
 \lim_{R\to\infty}\limsup_i
 \sum_A\int_{A_i(R)}S_i^2\,dV_{g_i}=0.
\]

The cover bounds the total scalar energy above by the sum of
these terms, and below by the base-core integral. Letting
$i\to\infty$ and then $R\to\infty$ gives
\begin{equation}\label{eq:uncenteredcalabicontinuity}
 \lim_i\int_M S_i^2\,dV_{g_i}
 =\int_X S_{h(\tau_*)}^2\,dV_{h(\tau_*)}.
\end{equation}
The volume and average scalar curvature agree with those of
the original flow by Theorem~\ref{thm:main}\textup{(III)}.
Subtracting $V\bar S^2$ proves
\[
 \Ca_\infty=\lim_i\Ca(g(u_i))
 =\Ca(h(\tau_*))=\Ca(h(0)).
\]
Holomorphic pullbacks preserve the Calabi functional, so
$\Ca(h(\tau))=\Ca_\infty$ for every $\tau$. Finally, for
$|\tau|\le T$ and all sufficiently large $i$,
\[
 \Ca(g(t_i+T))\le\Ca(g(t_i+\tau))\le\Ca(g(t_i-T)).
\]
Both endpoints tend to $\Ca_\infty$. This proves
\eqref{eq:maincalabiconvergence}.
\end{proof}

\begin{remark}\label{rem:furtherbubbling}
Once the canonical radius has a uniform positive lower bound,
further blow-ups on shrinking scales have the same geometric
features as blow-ups of extremal metrics. Theorem~\ref{thm:selected}
provides null scalar-flat ALE limits and smooth convergence on
their regular loci, with coverage of the corresponding regions.
Together with the model compactness and energy gap in Section~\ref{sec:stationary},
these ingredients also allow a bubble-tree description, as for
extremal metrics. We omit this additional structure to keep the
paper focused on flow convergence and continuity of the Calabi
functional.
\end{remark}

\subsection{Finite-group equivariant compactness}
\label{sec:equivariantcompactness}

We conclude with the equivariant version for a fixed finite group.

\begin{theorem}[Finite-group equivariant compactness]
\label{thm:equivariantcompactness}
Suppose that smooth K\"ahler metric families
$(M_j,g_j(\tau),J_j)$ converge in the sense of
Definition~\ref{def:topology}, on every bounded time interval to
$(X,h(\tau),J_\infty)$, where $X$ is a connected compact effective
K\"ahler orbifold with isolated quotient singularities.
Let a fixed finite group $\Gamma$ act on each $M_j$ by
time-independent biholomorphisms $\gamma_j$ preserving every
$g_j(\tau)$. After a further subsequence, there is one
time-independent action $\gamma\mapsto\gamma_\infty$ on $X$ by
orbifold biholomorphisms such that
\begin{equation}\label{eq:equivariantlimit}
(\gamma\eta)_\infty=\gamma_\infty\eta_\infty,\qquad
\gamma_\infty^*h(\tau)=h(\tau)
\quad(\gamma,\eta\in\Gamma,\ \tau\in\R).
\end{equation}
The GH correspondences can be chosen $\Gamma$-equivariantly,
retaining the compatible regular convergence and coverage in
Definition~\ref{def:topology}. In the original comparison charts,
the local action maps converge smoothly on regular compact
subsets to $\gamma_\infty$.
\end{theorem}
\begin{proof}
Fix $I:=[-T,T]$ and take the time-independent correspondences
$\mathscr C_j$ of Definition~\ref{def:topology}, with distortion
at most $\varepsilon_j\to0$ uniformly on $I$.
For $\gamma\in\Gamma$, set
\[
 T_{j,\gamma}:=
 \overline{\{(x,y):(p,x),(\gamma_jp,y)\in\mathscr C_j
                         \text{ for some }p\in M_j\}}.
\]
These relations project onto both factors of $X\times X$, and
any two pairs satisfy
\[
 |d_{h(\tau)}(x,x')-d_{h(\tau)}(y,y')|
 \le2\varepsilon_j,\qquad \tau\in I.
\]
After simultaneous Hausdorff extraction for the finite group,
each limiting relation is the graph of an isometry
$\gamma_\infty$ at time zero. Joint continuity and this inequality
make the same map an isometry at every time.
Following $p_j,\eta_jp_j,\gamma_j\eta_jp_j$ through
$\mathscr C_j$ proves the group law: the fibers have diameter
at most $\varepsilon_j$, so the last point has both limits
$\gamma_\infty\eta_\infty x$ and $(\gamma\eta)_\infty x$.

The correspondences
\[
 \mathscr C'_j:=\bigcup_{\gamma\in\Gamma}
       (\gamma_j\times\gamma_\infty)\mathscr C_j
\]
are exactly equivariant and contain the original chart graphs.
For each new pair $(p,x)$ choose an old pair $(p,y)$.
Convergence of $T_{j,\gamma}$ to the graph of $\gamma_\infty$
gives $d_{h(0)}(x,y)=o(1)$ uniformly. Uniform continuity of the
distance families therefore preserves vanishing distortion on $I$.
For $K\Subset K'\Subset U\Subset \mathcal R(X)$, such a $y$ lies
in $K'$ whenever $x\in K$ and $j$ is large. Applying
Lemma~\ref{lem:spacetimecoverage} to $K'\Subset U$ proves coverage.

Isometries preserve the regular locus, characterized by Euclidean
tangent cones. For $K\Subset \mathcal R(X)$, choose a larger regular
domain containing $K\cup\gamma_\infty K$. Graph convergence and
coverage make $\gamma_j^{\rm loc}:=\Phi_j^{-1}\gamma_j\Phi_j$
well-defined near $K$ and uniformly convergent to $\gamma_\infty$.
The identity
\[
 (\gamma_j^{\rm loc})^*(\Phi_j^*g_j(0))=\Phi_j^*g_j(0)
\]
bounds first derivatives; preservation of the Levi-Civita
connection and differentiation bound all higher derivatives.
Thus the convergence is smooth on $K$. Passing to the limit in
complex linearity proves holomorphicity there. By
\cite[Theorem 1.2 and Lemma 2.2]{Lange}, the distance isometry
lifts to smooth isometries of orbifold charts. Complex linearity
extends across their isolated quotient points by continuity.
Hence $\gamma_\infty$ is an orbifold biholomorphism.

Diagonal extraction over integer $T$ gives one action on all time
intervals: the regular local action maps have the same limits on
overlaps, and density identifies their extensions to $X$.
\end{proof}

\begin{corollary}[Equivariant limits of the Calabi flow]
\label{cor:equivariantcalabi}
Under the hypotheses of Theorem~\ref{thm:main}, suppose a finite
group $\Gamma$ acts on $M$ by biholomorphisms preserving
$\omega(0)$. Then every sequence $t_j\to\infty$ has a subsequence
with the convergence in that theorem enhanced by
Theorem~\ref{thm:equivariantcompactness}. The limiting action
preserves the singular set $\mathcal S$ and commutes with the
flow of the extremal vector field:
\[
\gamma_\infty\circ F_\tau=F_\tau\circ\gamma_\infty
\qquad(\gamma\in\Gamma,\ \tau\in\R).
\]
\end{corollary}
\begin{proof}
Uniqueness preserves the initial symmetries along the Calabi flow,
so Theorem~\ref{thm:equivariantcompactness} applies.
By \eqref{eq:Wnaturality} and
$\gamma_\infty^*h(0)=h(0)$, we have
$\gamma_\infty^*W_0=W_0$.
Hence $\gamma_\infty$ commutes with $F_\tau$.
Orbifold isomorphisms preserve $\mathcal S$.
\end{proof}

The GH correspondences are equivariant, although the smooth comparison
embeddings need not be. The limiting action is not necessarily effective.
Both results hold without assuming a polarization or a rational
K\"ahler class.

\vskip10pt
Haozhao Li, Institute of Geometry and Physics, and School of Mathematical
Sciences, University of Science and Technology of China, No. 96 Jinzhai
Road, Hefei, Anhui Province, 230026, China.
Email: \href{mailto:hzli@ustc.edu.cn}{hzli@ustc.edu.cn}\par\medskip

Bing Wang, Institute of Geometry and Physics, and School of Mathematical
Sciences, University of Science and Technology of China, Hefei 230026,
China, and Hefei National Laboratory, Hefei 230088, China.
Email: \href{mailto:topspin@ustc.edu.cn}{topspin@ustc.edu.cn}\par\medskip

Kai Zheng, University of Chinese Academy of Sciences, Beijing, 100190,
China.\\
Email: \href{mailto:KaiZheng@amss.ac.cn}{KaiZheng@amss.ac.cn}

\end{document}